\documentclass[12pt,reqno]{amsart}

\usepackage{amssymb,latexsym,graphicx}

\def\cal{\mathcal}

\def\Bbb{\mathbb}

\def\princa{\begin{enumerate}}

\newcommand\ej{\ifcase\value{potorro} \item  \or \end{list}%
 \or \end{list} \item \or \end{list} \item \or \end{list} \item%
 \or \end{list} \item \or \end{list} \item \or \end{list} \item%
 \or \end{list} \item \or \end{list} \item \or \end{list} \item%
 \or \end{list} \item \or \end{list} \item \or \end{list} \item  \fi}

\def\x#1{\def\@@eqncr{\let\@tempa\relax
  \ifcase\@eqcnt \def\@tempa{& & &}\or \def\@tempa{& &}
   \else \def\@tempa{&}\fi & #1\cr}}

\font\teneufm=eufm10
\font\seveneufm=eufm7
\font\fiveeufm=eufm5
\newfam\eufmfam
\textfont\eufmfam=\teneufm
\scriptfont\eufmfam=\seveneufm
\scriptscriptfont\eufmfam=\fiveeufm

\newfam\msbfam
\font\tenmsb=msbm10  scaled \magstep1  \textfont\msbfam=\tenmsb
\font\sevenmsb=msbm7  scaled \magstep1   \scriptfont\msbfam=\sevenmsb
\font\fivemsb=msbm5   scaled \magstep1  \scriptscriptfont\msbfam=\fivemsb

\def\RR{{\Bbb R}}
\def\CC{{\Bbb C}}

\def\ZZ{{\Bbb Z}}

 \def\HollowBox #1#2{{\dimen0=#1 \advance\dimen0 by -#2       
       \dimen1=#1 \advance\dimen1 by #2                       
        \vrule height #1 depth #2 width #2                    
        \vrule height 0pt depth #2 width #1                   
        \llap{\vrule height #1 depth -\dimen0 width \dimen1}%
       \hskip-#2                                           
       \vrule height #1 depth #2 width #2}}                   
 \def\BoxOpTwo{\mathord{\HollowBox{6pt}{.4pt}}\;}             

\def\crcrule{\hbox to
\textwidth{\rlap{\rule[-3.5pt]{84pt}{4pt}}\rule{\textwidth}{.5pt}}}
\def\thinrule{\vspace*{.13in}\hbox to \textwidth{\rule{\textwidth}{.5pt}}%
 \nopagebreak\vskip -1pc\nopagebreak}

\def\ra{\rightarrow}
\def\lra{\longrightarrow}
\def\Re{{\rm Re}\,}
\def\Im{{\rm Im}\,}

\def\wbar{\overline{w}} \def\ovw{\wbar}

\def\endpf{\medskip\hfill $\BoxOpTwo$}  

\def\sm{\setminus}

\def\btu{\bigtriangleup}

\def\arg{{\rm arg}\,}

\def\Res#1{\mathop{\rm Res}_{#1}\,}

\def\j2{\frac{j+1}{2}}
\def\nb{\nu_\beta}
\def\2bp{2\beta - \pi}
\def\sgn{{\rm sgn}}

\def\a{\alpha}
\def\b{\beta}

\def\d{\delta}
\def\l{\lambda}
\def\lj{\lambda^j}
\def\om{\omega} 
\def\p{\pi}
\def\x{\xi} 
\def\z{\zeta}

\def\x2p{x + 2\pi}
\def\ex2{e^{-x/2}}

\def\bmp{(\beta - \pi/2)}

\def\Jac{{\rm Jac\,}}

\newtheorem{theorem}{Theorem}[section]
\newtheorem{THM}{Theorem}

\newtheorem{proposition}[theorem]{Proposition}
\newtheorem{lemma}[theorem]{Lemma}

\newtheorem{remark}[theorem]{Remark}

\begin{document}

\parindent=17pt

\title[The Bergman kernel and 
projection on worm domains]{\large The Bergman kernel and projection on 
  non-smooth worm domains}

\author[S. Krantz]{Steven G. Krantz}
\address{American Institute of Mathematics, 360 Portage Avenue,
Palo Alto, California 94306-2244} \hfill \break
\email{{\tt skrantz@aimath.org}}
\vspace*{.15in}
\author[M. Peloso]{Marco M. Peloso}
\address{Dipartimento di Matematica, Corso Duca degli Abruzzi 24,
  Politecnico di Torino, 10129 
Torino, Italy}
\email{{\tt peloso@calvino.polito.it}}
\thanks{First author was supported
in part by NSF Grant DMS-9988854 and by a grant from the Dean of the Graduate
School at Washington University in St. Louis.  The second author was
supported in part by  
the 2002 Cofin project {\em Analisi Armonica}.}
									    
\begin{abstract}
We study the Bergman kernel and projection on
the worm domains 
$$
D_\b = \left \{\z \in \CC^2:\, 
\Re \bigl(\z_1 e^{-i\log |\z_2|^2} \bigr)>0 ,\, 
\big|\log |\z_2|^2\big| < \b - \frac{\pi}{2} \right \} 
$$
and
$$
D_\b' = \left \{z \in \CC^2 :\, \big|\Im z_1 - \log |z_2|^2\big| <
  \frac{\pi}{2},\, |\log |z_2|^2| < \b - \frac{\pi}{2} \right \} 
$$
for $\beta > \pi$.  
These two domains are biholomorphically 
equivalent via the mapping
$$
 D_\b' \ni (z_1,z_2) \mapsto (e^{z_1},z_2)\ni D_\b\, .
$$
We calculate the kernels explicitly, up to an error term that can be
controlled. 	    

As a result, we can determine the $L^p$-mapping properties of
the Bergman projections on these worm domains.  
Denote by $P$ the Bergman projection on $D_\b$ and by $P'$ the one
on $D_\b'$.  We calculate the sharp
range of $p$ for which the Bergman projection is
bounded on $L^p$.
More precisely we show
that 
$$
P': L^p(D_\b')\lra L^p(D_\b')
$$
boundedly
when $1<p<\infty$, while
$$
P: L^p(D_\b)\lra L^p(D_\b)
$$
if and only if $2/(1+\nu_\b)<p<2/(1-\nu_\b)$, where
$\nu_\b=\pi/(\2bp)$. 
Along the way, we give a new proof of the failure of Condition $R$ on
these worms.

Finally, we are able to show that the singularities of the
Bergman kernel on the boundary are not contained in the
boundary diagonal.
\end{abstract}

\maketitle

\thispagestyle{empty}

\section*{Introduction}

The (smooth) worm domain was first created by Diederich and Forn\ae ss
\cite{DF} to provide counterexamples to certain classical conjectures
about the geometry of pseudoconvex domains.  Chief among these
examples is that the smooth worm is bounded and pseudoconvex and
smooth yet its closure does not have a Stein neighborhood basis.  More
recently, thanks to work of Kiselman \cite{KIS}, Barrett \cite{BAR2}, 
and Christ
\cite{CHR}, the worm has played an important role in the study of the
$\overline{\partial}$-Neumann problem and Condition $R$ (see
\cite{CHS}).  Recall that 
Condition $R$ for a domain $\Omega$ is the assertion that the Bergman
projection maps $C^\infty(\overline{\Omega})$ to
$C^\infty(\overline{\Omega})$.  Work of Bell \cite{BEL1} has shown that
Condition $R$ is closely related to the boundary regularity of
biholomorphic mappings.  See also \cite{BOS2}.  Condition $R$ fails on the
worm domains.

In more detail, in the seminal paper \cite{CHR}, Christ shows that the
$\overline{\partial}$-Neumann problem is not globally hypoelliptic on
the smooth worm domain.  It follows then that Condition $R$ fails on
the smooth worm.  Christ's work provides considerable motivation for
attaining a deeper and more detailed understanding of the Bergman
kernel and projection on the worm domains.  That is what the present
paper achieves---for the particular worm domains\footnote{Throughout
this paper, we use notation for the different worm domains that
is recorded in \cite{CHS}.  This notation has become rather standard.}
$$
D_\b = \left \{\z \in \CC^2: 
\Re \bigl(\z_1 e^{-i\log |\z_2|^2}
\bigr)>0 ,\, \big|\log |\z_2|^2\big| < 
\b - \frac{\pi}{2} \right \} 
$$
and
$$
D'_\b = \left \{z 
\in \CC^2: \big|\Im z_1 - \log |z_2|^2\big| 
< \frac{\pi}{2},\,  |\log |z_2|^2| <
\b - \frac{\pi}{2} \right \} \, . 
$$

It should be noted that these two domains 
are biholomorphically equivalent via the mapping
$$
(z_1,z_2)\ni D_\b' \mapsto (e^{z_1},z_2)\ni D_\b\, .
$$
Moreover, these domains are {\it not} smoothly
bounded.   Each boundary is only Lipschitz, and, particularly, its
boundary is Levi flat, and clearly each is unbounded. 

These domains are rather different from the smooth worm
$\mathcal W_\b$ (see \cite{CHS} and the discussion below), which 
has all boundary points, except those on a singular annulus in the
boundary, strongly pseudoconvex. However our worm domain $D_\b$ 
is actually  a model for the smoothly bounded $\mathcal W_\b$ (see, for
instance, \cite{BAR2}), and it can be
expected that phenomena that are true on $D_\b$ will
in fact hold on $\mathcal W_\b$ as well.  We will say more about
this symbiotic relationship below.  We intend to explore the other
worms, particularly the smooth worm, in future papers.

The first author thanks Michael Christ for a very helpful
conversation. He also thanks the Mathematical Sciences
Research Institute and the American Institute of Mathematics
for generous and gracious support during a part of the writing
of this paper.			  

\section*{Table of Contents}
\vspace*{.15in}

\begin{description}   
\item[{\bf Section 1}] {\sl  
Statement of the Main Results}   
\item[{\bf Section 2}] {\sl Analysis of the Singularity of the Bergman
    Kernels}    
\item[{\bf Section 3}] {\sl Boundedness of the Bergman Projection on
    $D'_\b$}	   
\item[{\bf Section 4}] {\sl Boundedness of the Bergman Projection on $D_\b$}
\item[{\bf Section 5}]  {\sl  Decomposition of the Bergman
    Space} 
\item[{\bf Section 6}] {\sl The Sum of the $R_j$}	 
\item[{\bf Section 7}] {\sl Decomposition of the $J_j$} 
\item[{\bf Section 8}] {\sl The Bergman Kernels for $D_\b$  and
    $D'_\b$\! --\!    
Proof of Theorems \ref{THM1} and \ref{THM2}}
\item[{\bf Section 9}]  {\sl Foundational Steps in the Proofs of Propositions
    \ref{sum-of-Mj-proposition} and
\ref{sum-of-error-terms-proposition}}
\item[{\bf Section 10}]  {\sl Proof of Propositions 
\ref{sum-of-Mj-proposition} and
\ref{sum-of-error-terms-proposition}}
\item[{\bf Section 11}] {\sl Completion of the Proof of Proposition 
\ref{technical}}
\item[{\bf Section 12}] {\sl Concluding Remarks}    \bigskip
\end{description}

\section{Statement of the Main Results}

Throughout this paper we work on worm domains in $\CC^2$.   

If $r>0$ then we let $I_r = \{x \in \RR: -r \leq x \leq r\} = [-r,r]$.
Let $\b>\pi/2$ and let $\mathcal W$ denote the 
domain
$$
\mathcal W = \left \{(z_1, z_2)\in\CC^2:\,  
\big|z_1 - e^{i\log |z_2|^2} \big|^2 < 1 - \eta(\log|z_2|^2) \right \}
\, , 
$$
where
\begin{itemize}
\item[(i)]  $\eta \geq 0$, $\eta$ is even, $\eta$ is convex;\smallskip
\item[(ii)]  $\eta^{-1}(0) = I_{\b - \pi/2} \equiv [-\b + \frac{\pi}{2},
    \b - \frac{\pi}{2}]$; \smallskip
\item[(iii)]  there exists a number $a > 0$ such that $\eta(x) > 1$ if
    $x < - a$ or $x > a$; \smallskip
\item[(iv)]  $\eta'(x) \ne 0$ if $\eta(x) = 1$.
\end{itemize}
We will also write $\mathcal W=\mathcal W_\b$. 
Notice that the slices of $\mathcal W$ for $z_2$ fixed are discs
centered on the unit circle with centers that wind $(\2bp)/2\pi$ times
about that circle as $|z_2|$ traverses the range of values for which
$\eta(\log|z_2|^2) < 1$. 

It turns out that $\mathcal W$ is smoothly bounded and
pseudoconvex (see \cite{CHS}). It is the classical, smooth
worm studied by Diederich and Forn\ae ss (see \cite{DF}).
Moreover, its boundary is {\it strongly} pseudoconvex except
at the boundary points $(0,e^{i\log |z_2|^2})$ for $\big|\log
|z_2|^2 \big|<\b-\pi/2$ ---see \cite{DF} or \cite{CHS}. Notice
that these points constitute an annulus in $\partial\mathcal
W_\b$.

Starting from work by Kiselman \cite{KIS}, Barrett \cite{BAR2} showed
that the Bergman projection $P_{\mathcal W}$ of $\mathcal W$
does not map
the
$s$-order Sobolev space $W^s(\mathcal W)$ 
boundedly into itself
when $s\ge\nu_\b$, where
\begin{equation}\label{nu-beta}
\nu_\b= \frac{\pi}{\2bp}
\end{equation}
is half the reciprocal of the number of windings of $\mathcal W$.

In a later paper \cite{CHR}, Christ showed that the
$\bar\partial$-Neumann problem is not globally hypoelliptic on
$\mathcal W$, a fact that is equivalent to
the property that $P_{\mathcal W}:C^\infty(\overline{\mathcal W})
\not\ra C^\infty(\overline{\mathcal W})$.  In other words, $\mathcal W$
does not satisfy Condition $R$.\medskip

It is still an open question, of great interest, whether a
biholomorphic mapping of $\mathcal W$ onto another 
smoothly bounded pseudoconvex domain
extends smoothly to a diffeomorphism of the boundaries.
In this direction, it is worth mentioning that Chen \cite{CHE1} has
shown that the automorphism group of $\mathcal W$ reduces to the
rotations in the $z_2$-variable; hence all biholomorphic {\it self-maps} of
$\mathcal W$ do extend smoothly to the boundary.

It is also noteworthy that Chen \cite{CHE1} and Ligocka \cite{LIG}
have independently showed that the Bergman kernel of $\mathcal W$ cannot
lie in 
$C^\infty(\overline{\mathcal W}\times\overline{\mathcal W}\sm\btu)$, 
where $\btu$ is the boundary diagonal.  In fact, in \cite{CHE1} it is
shown that this phenomenon is a consequence of the presence of a
complex variety 
in the boundary of $\mathcal W$. 
\medskip
	 
There are other worm domains, which are neither bounded nor smooth, 
that are nonetheless of considerable interest.  These are:
the unbounded, non-smooth worm with
half-plane slices
\begin{equation}\label{D-beta}
D_\b = 
\left \{(\z_1, \z_2) \in \CC^2: \Re (\z_1 e^{-i \log |\z_2|^2}) > 0 \, ,
     \big|\log |\z_2|^2\big| < \b - \frac{\pi}{2} \right \}  \, 
\end{equation}
and the unbounded, non-smooth worm
with strip slices
\begin{equation}\label{D-beta'}
D'_\b = \left\{ (z_1, z_2) \in \CC^2: 
\big|\Im z_1 - \log |z_2|^2 \big| < \frac{\pi}{2} \, ,
     \big|\log |z_2|^2\big| < \b - \frac{\pi}{2} \right\} \, .  
\end{equation}
These  domains 
are biholomorphically equivalent via the mapping
$$
\Phi: D'_\beta \ni (z_1,z_2) \mapsto (e^{z_1},z_2) \equiv (\z_1,\z_2) \in D_\b\, .
$$
Kiselman introduced the domains $D_\b$, $D'_\b$ in \cite{KIS}.
In his fundamental work \cite{BAR2},
Barrett showed that it is possible to obtain information on the
Bergman kernel and projection on $\mathcal W$ 
from 
corresponding information on the non-smooth worm $D_\b$ by using
an exhaustion and limiting argument.  
\medskip

It is therefore of interest, in its own right and as a model
for the smooth case, to study the behavior of the Bergman
kernel and projection on the domain $D_\b$ and its
biholomorphic copy $D_\b'$. It is obvious from the
transformation rule of the Bergman kernel (see \cite{KRA1})
that it suffices to obtain the expression for the kernel in
one of the two domains. However, the $L^p$-mapping properties
of the Bergman projections of the two domains turn out to be
substantially different (just because $L^p$ spaces of
holomorphic functions do not transform canonically under
biholomorphic maps when $p\ne2$).

In this paper we determine, in
Theorems \ref{THM1} and \ref{THM2}, the explicit expression for the Bergman
kernels for $D_\b$ and $D_\b'$---up to controllable error terms.  
Once these are available we study the 
$L^p$-mapping properties of the corresponding Bergman
projections in Theorems \ref{THM3} and \ref{THM4}.  

More precisely we prove the following results. 

\begin{THM}\label{THM1}	 \sl
Let $c_0$ be a positive fixed constant.  Let
$\chi_1$ be a smooth cuf-off function on the real line,
supported on $\{x:\, |x|\le 2c_0\}$,
identically 1 for $|x|<c_0$. Set $\chi_2 =1-\chi_1$.

Let $\b>\pi$ and let $\nu_\b$ be defined as in (\ref{nu-beta}).
Let $h$ be fixed, with
\begin{equation}\label{(*)}
\nu_\b<h<\min(1,2\nu_\b) \, .  
\end{equation}
Then there exist
functions $F_1,F_2,\dots, F_8$ and $\tilde F_1,\tilde F_2,\dots, \tilde F_8$,
holomorphic in $z$ and anti-holomorphic in $w$, for $z=(z_1,z_2)$,
$w=(w_1,w_2)$ varying in a
neighborhood of $D'_\b$, and 
having size ${\cal O}(|\Re z_1-\Re w_1|)$, together with all their
derivatives, for $z,w \in \overline{D'_\b}$,
as $|\Re z_1-\Re w_1|\ra+\infty$.  
Moreover, there exist functions $E,\tilde E \in {\mathcal C}^\infty 
(\overline{D'_\b}\times \overline{D'_\b})$ such that
$$
D_{z_1}^\alpha D_{w_1}^\gamma E(z,w),
\, D_{z_1}^\alpha D_{w_1}^\gamma \tilde E(z,w)
=\mathcal O(|\Re z_1-\Re w_1|^{|\alpha|+|\gamma|})\ ,
$$
as $|\Re z_1-\Re w_1|\ra+\infty$. (Here, for $\lambda\in\CC$,
$D_\lambda$ denotes the 
partial derivative in $\l$ or $\bar\l$.)

Then the following holds.  Set
\begin{align}
K_b(z,w)
& =   \frac{F_1(z,w)
}{(i(z_1- \ovw_1)+2\b)^2 (e^{\bmp} - z_2 \ovw_2)^2}  
\notag \notag\\    
&\qquad\quad
+ \frac{F_2(z,w)}{(i(z_1- \ovw_1)+2\b)^2
(z_2 \ovw_2 - e^{-[i(z_1-  \ovw_1) +\pi]/2})^2 }
\notag \notag\\  
& \quad
+
\frac{F_3(z,w)}{(e^{[\pi-i(z_1- \ovw_1)]/2}-z_2 \ovw_2)^2
(e^{\bmp}-z_2 \ovw_2)^2}
\notag \notag\\
& \qquad +  
\frac{F_4(z,w)}{(i(z_1- \ovw_1)-2\b)^2 
(e^{[\pi - i(z_1- \ovw_1)]/2} - z_2 \ovw_2)^2}
\notag \notag\\
&\qquad\quad
+ \frac{F_5(z,w)}{(i(z_1- \ovw_1)-2\b)^2  
(z_2 \ovw_2 - e^{-\bmp})^2}  
\notag \notag\\
&\quad 
+ \frac{F_6(z,w)
}{(e^{-[i(z_1- \ovw_1)+\pi]/2}-z_2 \ovw_2)^2
(e^{-\bmp}-z_2 \ovw_2)^2} 
 \notag \notag\\
& \quad
+ \frac{F_7(z,w)}{
(i(z_1- \ovw_1)+2\b)^2(e^{\bmp} - z_2 \ovw_2)
(e^{-[i(z_1- \ovw_1) + \pi]/2}-z_2 \ovw_2)} 
\notag \notag\\
& \quad
+ \frac{F_8(z,w)}
{(i(z_1- \ovw_1)-2\b)^2  
(e^{-\bmp}-z_2 \ovw_2)(e^{[\pi-i(z_1- \ovw_1)]/2}-z_2 \ovw_2)}  
+E(z,w) \notag\\ 
& \equiv K_1(z,w)+\cdots+K_8(z,w) + E(z,w)\ .
\label{K1throughK8}
\end{align}
Define $K_{\tilde b}$ by replacing $F_1,\dots,F_8$ and $E$ by
$\tilde F_1,\dots,\tilde F_8$ and $\tilde E$ 
and thus $K_1,\dots,K_8$  by
$\tilde K_1,\dots,\tilde K_8$ respectively in formula
(\ref{K1throughK8}).  

Then there exist functions $\phi_1,\phi_2$ entire in $z$ and
$\ovw$ (that is, anti-holomorphic in 
$w$), which are of size ${\cal O}(|\Re z_1 -\Re w_1|)$, together with
all their 
derivatives, uniformly in all closed strips 
$\{|\Im z_1|+|\Im w_1|\le C\}$, such that  
the Bergman kernel $K_{D'_\b}$ on $D'_\b$ admits the asymptotic expansion
\begin{multline*}
K_{D'_\b}(z,w) =
\chi_1(\Re z_1 -\Re w_1) K_b(z,w) + 
\chi_2(\Re z_1 -\Re w_1)
\biggl\{ 
e^{-h\sgn (\Re z_1 -\Re w_1) \cdot (z_1- \ovw_1)}   
K_{\tilde b}(z,w) 
\\
+e^{-\nb \sgn (\Re z_1 -\Re w_1)\cdot  (z_1- \ovw_1)}
\biggl( 
 \frac{\phi_1 (z_1,w_1)}{(e^{[\pi- i(z_1- \ovw_1)]/2}- z_2 \ovw_2)^2} 
   + \frac{\phi_2 (z,w)}{
(e^{-[i(z_1- \ovw_1)+\pi]/2}-z_2 \ovw_2)^2} \biggr) \biggr\} \, .
\end{multline*}
Here $h$ is specified as in (\ref{(*)})  above.
\end{THM}

\begin{THM}\label{THM2}	 \sl
With the notation as in Theorem \ref{THM1}, 
there exist
functions $g_1,g_2$, 
$G_1,G_2,\dots, G_8$ and $\tilde G_1,\tilde G_2,\dots, \tilde G_8$,
holomorphic in $\z$ and anti-holomorphic in $\om$, for $\z=(\z_1,\z_2)$,
$\om=(\om_1,\om_2)$ varying in  $\overline{D'_\b}\setminus\{(0,z_2)\}$, 
such that
$$
\partial_{\z_1}^\alpha \partial_{\overline{\om}_1}^\gamma G(\z,\om)=
{\cal O} \bigl( |\z_1|^{-|\alpha|} |\om_1|^{-|\gamma|} \bigr)
\qquad\qquad\text{as}\quad |\z_1|,|\om_1|\ra 0\ ,
$$
where $G$ denotes any of the functions $g_j$, $G_j$, $\tilde G_j$.
Moreover, there exist functions $E,\tilde E \in {\mathcal C}^\infty
\bigl( \overline{D'_\b}\setminus\{(0,z_2)\} \times
\overline{D'_\b}\setminus\{(0,z_2)\} \bigr)$ such that
$$
D_{\z_1}^\alpha D_{\om_1}^\gamma E(\z,\om),
D_{\z_1}^\alpha D_{\om_1}^\gamma \tilde E(\z,\om) =
{\cal O} \bigl( |\z_1|^{-|\alpha|} |\om_1|^{-|\gamma|} \bigr)
\qquad\qquad\text{as}\quad |\z_1|,|\om_1|\ra 0\ .
$$
(Here $D_\lambda$, 
for $\lambda\in\CC$,
$D_\lambda$ denotes the 
partial derivative in $\l$ or $\bar\l$.)

Then the following holds.  Set
\begin{align} 
H_b(\z,w)
& =   \frac{G_1(\z,w)
}{(i\log(\z_1/\overline{\om}_1)+2\b)^2 (e^{\bmp} - \z_2 \overline{\om}_2)^2}  
\notag\\    
&\qquad\quad
+ \frac{G_2(\z,w)}{(i\log(\z_1/\overline{\om}_1)+2\b)^2
\bigl((\z_1/\overline{\om}_1)^{-i/2}e^{-\pi/2} 
-\z_2 \overline{\om}_2  
\bigr)^2 }
\notag\\  
& \quad
+
\frac{G_3(\z,w)}{\bigl((\z_1/\overline{\om}_1)^{-i/2}e^{\pi/2} 
-\z_2 \overline{\om}_2\bigr)^2
(e^{\bmp}-\z_2 \overline{\om}_2)^2}
\notag\\
& \qquad +  
\frac{G_4(\z,w)}{(i\log(\z_1/\overline{\om}_1)-2\b)^2 
\bigl((\z_1/\overline{\om}_1)^{-i/2}e^{\pi/2} 
- \z_2 \overline{\om}_2\bigr)^2}
\notag\\
&\qquad\quad
+ \frac{G_5(\z,w)}{(i\log(\z_1/\overline{\om}_1)-2\b)^2  
(e^{-\bmp}-\z_2 \overline{\om}_2)^2}  
\notag\\
&\quad 
+\frac{G_6(\z,w)
}{\bigl((\z_1/\overline{\om}_1)^{-i/2}e^{-\pi/2} 
-\z_2 \overline{\om}_2\bigr)^2
(e^{-\bmp}-\z_2 \overline{\om}_2)^2} 
 \notag\\
&\quad
+ \frac{G_7(\z,w)}{ 
(i\log(\z_1/\overline{\om}_1)+2\b)^2(e^{\bmp} - \z_2 \overline{\om}_2)
\bigl((\z_1/\overline{\om}_1)^{-i/2}e^{-\pi/2}
-\z_2 \overline{\om}_2\bigr)} 
\notag\\
&\quad
+ \frac{G_8(\z,w)}
{(i\log(\z_1/\overline{\om}_1)-2\b)^2  
(e^{-\bmp}-\z_2 \overline{\om}_2)\bigl((\z_1/\overline{\om}_1)^{-i/2}e^{\pi/2} 
-\z_2 \overline{\om}_2\bigr)}  
+E(\z,w)\notag\\
& \equiv H_1(\z,\om)+\cdots+H_8(\z,\om) + E(\z,\om)\ . 
\end{align}
Define $H_{\tilde b}$ by replacing $G_1,\dots,G_8$ and $E$ by
$\tilde G_1,\dots,\tilde G_8$ and $\tilde E$, and $H_1,\dots,H_8$  by
$\tilde H_1,\dots,\tilde H_8$,  respectively.

Then, setting $t=|\z_1|-|\om_1|$, we have this
asymptotic expansion for the Bergman kernel on $D_\b$:
\begin{align*}
& K_{D_\b} \bigl((\z_1,\z_2),(\om_1,\om_2)\bigr)\\
&= \chi_1(t)
\frac{H_b(\z,\om)
}{\z_1 \overline{\om}_1} +
\chi_2(t) 
\biggl\{  \biggl(\frac{|\z_1|}{|\om_1|}\biggr)^{-h\sgn t} 
e^{-h\sgn t \cdot (\arg\z_1+\arg\om_1)}   
\frac{H_{\tilde b}(\z,\om)}{ \z_1
  \overline{\om}_1}  
  \notag\\
& \quad
+ \biggl(\frac{|\z_1|}{|\om_1|}\biggr)^{-\nu_\b\sgn t} 
e^{-\nu_b\sgn t \cdot (\arg\z_1+\arg\om_1)}\biggl(    
\frac{g_1 (\z_1,\om_1)}{ \z_1\overline{\om}_1}\cdot
\frac{1}{\bigl(
( \z_1/\overline{\om_1})^{-i/2}
e^{\pi/2}- z_2 \overline{\om}_2 \bigr)^2} \notag\\
& \quad\qquad\qquad   
+\frac{ g_2 (\z,\om)}
{\z_1\overline{\om}_1} \cdot \frac{1}{\bigl( 
( \z_1/\overline{\om}_1)^{-i/2}
e^{-\pi/2}- z_2 \overline{\om}_2 \bigr)^2} \biggr)  \biggr\}
\, ,
\end{align*}
where $h$ is defined in (\ref{(*)}).
\end{THM}

\begin{THM}\label{THM3} \sl
Let $P$  denote the Bergman projection 
on the domain $D_\b$, $\b>\pi$. Then
$$
P: L^p(D_\b) \ra L^p(D_\b) 
$$
 is bounded
if and only if
$$
\frac{2}{1 + \nu_\b} < p < \frac{2}{1 - \nu_\b} \ .
$$
\end{THM}

\begin{THM}\label{THM4}	   \sl
Let $P'$  denote the Bergman projection 
on the domain $D'_\b$, $\b>\pi$. 
Then 
$$
P': L^p(D'_\b) \ra L^p(D'_\b) 
$$
is bounded for $1 < p < \infty$.
\medskip   \\
\end{THM}

{\bf Remark.} A comment about the parameter $\b$ is now in order. The
definition of $D_\b$ and $D'_\b$ (as well as the one of $\mathcal
W_\b$) requires that $\b>\pi/2$.  

However, for simplicity of the arguments, in this paper
we restrict ourselves to the
case $\b>\pi$.  This is not a serious constraint.  The most
interesting situations occur as $\b\ra+\infty$ (which means more twists in
the geometry of the worm), and we believe that
the expansion of the Bergman kernels for $D_\b$ and $D'_\b$ (Theorems
\ref{THM1} and \ref{THM2}) will 
extend to the case $\pi/2<\b\le\pi$. 

On the other hand, the result about the $L^p$ boundedness of the Bergman
projection on $D_\b$ (Theorem \ref{THM4}) requires that $\nu_b$,
defined in (\ref{nu-beta}), satisfies $\nu_\b<1$.  It is immediate to
check that $\nu_\b<1$ if and only if $\b>\pi$.  It would 
therefore be of interest to study the 
$L^p$ boundedness of the Bergman projection on $D_\b$ when 
$\pi/2<\b\le\pi$.
\medskip   \\

The paper is organized as follows. In Section
\ref{added-section} we discuss the singularities of the kernels whose
expansions are given in Theorems \ref{THM1} and \ref{THM2}.
Theorems \ref{THM3} and \ref{THM4} of that section treat the $L^p$
boundedness of 
the Bergman projection.  Section \ref{added-section} 
also contains a comment 
concerning the
failure of Condition $R$ on the worm domains.

In
Sections \ref{section-boundedness-D-prime-beta} and 
\ref{section-boundedness-D-beta}, assuming the validity of the
expansions 
in Theorems \ref{THM1} and \ref{THM2},  we prove 
Theorems \ref{THM4} and \ref{THM3}, respectively.

Section \ref{section-boundedness-D-beta} 
 records (Theorem \ref{straube-thm}) an explicit result about limitations
of $L^p$ boundedness of the Bergman projection on the smooth worm
${\mathcal W}_\beta$.  This result is proved by way of an exhaustion
procedure \`{a} la Barrett \cite{BAR2}.

Sections \ref{decomposition-B-space} through \ref{appendix} are
devoted to the proof of 
the 
expansions (i.e., Theorems \ref{THM1} and \ref{THM2}) for the Bergman
kernels on the  
worms $D'_\b$ and $D_\b$.
We do the bulk of our work in this paper directly with $D'_\b$.
At the end, we shall transfer the result to $D_\b$
via a biholomorphic mapping.

It is worth noting that the proofs of a number of our more technical results,
including Proposition \ref{technical}, 
Lemmas \ref{NEW-LEMMA-IANDIII}--\ref{NEW-LEMMA-II*} and Lemma
\ref{lemma-Q}, are relegated to  
Section \ref{appendix} of the paper.
\bigskip   \\

\section{Analysis of the Singularity of the Bergman Kernels}
\label{added-section}
\vspace*{.12in}

\subsection*{Singularities of the Kernel on \boldmath $D'_\b$}
We can draw a picture of the section of 
$D'_\b$ in the $(\Im z_1,\log|z_2|)$-plane.  

 \begin{figure}
      \includegraphics[height=3.15in, width=4.6in]{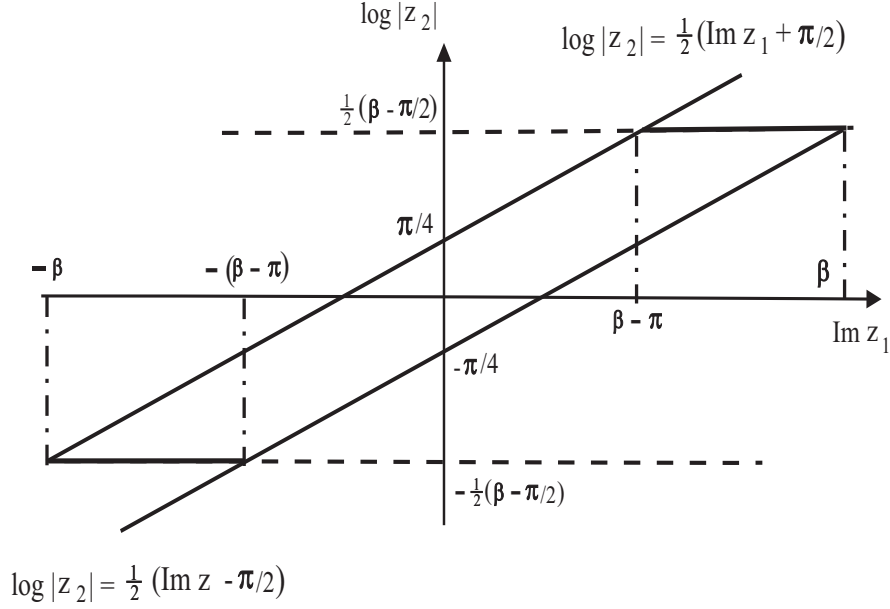}
      \caption{A representation of the domain $D'_\b$.}
      \label{figure1}
    \end{figure}

Recall that 
$$
D'_\b = \left\{ (z_1, z_2) \in \CC^2: 
\big|\Im z_1 - \log |z_2|^2 \big| < \frac{\pi}{2} \, ,
     \big|\log |z_2|^2\big| < \b - \frac{\pi}{2} \right\} \, ,  
$$
and recall
the expansion of $K_{D'_\b}$ given in
Theorem \ref{THM1}.  

We begin by analyzing the behavior at the ``bounded portion of the
boundary''. 
Then (using the notation from Theorem \ref{THM1}) we notice the
following facts: 

\smallskip
\subsection*{\boldmath $(K_1)$}
For $z,w\in D'_\b$ the term $K_1$ becomes singular (if
and) only if
$$
z_2\ovw_2\ra e^{\b-\pi/2}\qquad\text{or}\quad
\Im(z_1-\ovw_1)\ra 2\b\ .
$$
This can happen only if $\log|z_2|^2,\, \log|w_2|^2 \ra
\b-\pi/2$. Thus $K_1$ is singular only when both $z$ and $w$ tend to 
the top side of the
domain in Figure 1.

Notice that the top right corner of the domain in Figure 1
does not correspond to a single point in the domain $D'_\b$.  In fact,
we could rotate in the $z_2$-variable as well as translate the
$z_1$-variable by a real constant.  Thus the singularity of the
kernel $K_1$ is not contained in the boundary diagonal of
$\overline{D'_\b}\times\overline{D'_\b}$, since for $w\in
D'_\b$ fixed, the singular set
is
$\{ (z_1,z_2)\in\overline{D'_\b}:
|z_2|=e^{[\b-\pi/2]/2},\, -(\b-\pi)<\Im z_1<\b \}$.
The interior of this set has real dimension 3.
This phenomenon appears in the case of all the
other terms, and we shall not repeat this comment again.

\smallskip
\subsection*{\boldmath $(K_5)$} This term is symmetric to $K_1$ and it
is singular 
as  
$$
\log|z_2|^2,\, \log|w_2|^2 \ra
-(\b-\pi/2)\qquad
\text{and, in this case, also when}\quad
\Im z_1, \Im w_1 \ra -\b\ .
$$
Thus, $K_5$ is singular only 
when both $z$ and $w$ tend to 
the lower side
 of the domain in
Figure 1. 

\smallskip
\subsection*{\boldmath $(K_2)$} 

 For $z,w\in D'_\b$ the term $K_2$ becomes singular (if
and) only if
$$
z_2\ovw_2\ra e^{-[i(z_1-\ovw_1)+\p]/2}\qquad\text{or}\quad
\Im(z_1-\ovw_1)\ra 2\b\ .
$$
This can happen only if $\log|z_2|^2\ra\Im z_1-\frac\pi2$ and 
$\log|z_2|^2\ra\Im z_1-\frac\pi2$ and 
thus $K_2$ is singular only when both $z$ and $w$ tend to 
the right oblique line that bounds the
domain in Figure 1.

\smallskip
\subsection*{\boldmath $(K_4)$} 

This term is symmetric to $K_2$ and it is singular
on the left oblique line of the domain in Figure 1.

\smallskip
\subsection*{\boldmath $(K_3,\, K_6)$} 

The kernel $K_3$ is singular when 
$$
z_2\ovw_2\ra e^{[\p-i(z_1-\ovw_1)]/2}\qquad\text{or}\quad
z_2\ovw_2\ra e^{\b-\pi/2}\ .
$$
Therefore $K_3$ is singular when 
both $z$ and $w$ tend to either the
upper horizontal line, or the right oblique line on the boundary of
the domain in Figure 1.

Analogously, 
the kernel $K_6$ is singular when 
both $z$ and $w$ tend to either the
lower horizontal line, or the left oblique line on the boundary of
the domain in Figure 1.

\smallskip
\subsection*{\boldmath $(K_7,\, K_8)$} 

Finally, these kernels have singularities of types that have been
discussed already.  In fact, there exists a constant $C>0$ such that
$$
|K_7(z,w)| \le C \bigl( |K_1(z,w)|+|K_2(z,w)|\bigr)\, ,
$$
and 
$$
|K_8(z,w)| \le C \bigl( |K_4(z,w)|+|K_5(z,w)|\bigr)\, ,
$$

\smallskip

We finish this part by noticing that, 
as $|\Re z_1|,|\Re w_1| \ra +\infty$,
the principal term of the kernel  is $\approx 
e^{-\nu_\b|\Re z_1-\Re w_1|}$, as Theorem \ref{THM1} clearly indicates.
\medskip   \\

\subsection*{Singularities of the Kernel on \boldmath $D_\b$}
Now we turn to the Bergman kernel for the domain $D_\b$; we
recall that the latter is given by
$$
D_\b
= 
\left \{(\z_1, \z_2) \in \CC^2: \Re (\z_1 e^{-i \log |\z_2|^2}) > 0 \, ,
     \big|\log |\z_2|^2
\big| < \b - \frac{\pi}{2} \right \}  \, .
$$
The kernel is described in Theorem \ref{THM2}.

First of all, we comment on the definition of the function $\log\z_1$
on $D_\b$.  For $\z=(\z_1,\z_2)\in D_\b$ and $\zeta_2$ fixed, $\z_1$
varies in the half-plane $\Re(\z_1 e^{-i\log|\z_2|^2})>0$, that is,
\begin{equation}\label{inequality1-for-D-beta}
-\frac\pi2 + \log |\z_2|^2 < \arg \z_1 < \frac{\pi}{2}
    + \log |\z_2|^2 \ .
\end{equation}
Therefore $\log \z_1$ is well defined as a function of
$(\z_1,\z_2)\in D_\b$.

Notice also that $D_\b$ can be described  by inequality
(\ref{inequality1-for-D-beta}) together with the condition
$-(\b-\pi/2)<\log|\z_2|^2 < \b-\pi/2$. 
These inequalities also imply that $-\b<\arg\z_1<\b$.  
\medskip \\

The terms 
$i\log(\z_1/\overline{\om}_1)\pm2\b$ 
in the expansion for the kernel become singular 
as $\arg \z_1,\, \arg \om_1 \ra \pm \b$.  But this can happen only if 
$\log |\z_2|^2, \log|\om_2|^2 \ra 
\pm(\b - \pi/2)$, that is
$|\z_2\overline{\om}_2| \ra e^{\pm (\b - \pi/2)}$.
\medskip    \\

Next observe that
$|\z_2\overline{\om}_2| \ra e^{\pi/2} \big| \bigl( 
\z_1/\overline{\om}_1 \bigr)^{\pm i/2} \big|=  
e^{\pi/2 \mp ({\rm arg}\, \z_1 + {\rm arg}\, \om_1)/2}$, 
which occurs exactly when 
$$
\log |\z_2|^2 \ra  \frac\pi2 \mp\arg\z_1
\quad\text{and}\quad
\log|\om_2|^2 \ra \frac\pi2 \mp\arg\om_1 \, .
$$

\subsection*{Failure of Condition R}
From the expansion given by 
Theorem \ref{THM2} it is clear that, for any $\om\in D_\b$ fixed, the
$A^2(D_\b)$-function $K_{D_\b}(\cdot,\om)$ does not belong to
$L^2\bigl(D_\b,|\z_1|^{-2\nu_\b}dV\bigr)$.  This in particular implies
that $K_{D_\b}(\cdot,\om)\not\in W^{\nu_\b}(D_\b)$ and that, therefore,
condition R fails on $D_\b$.
This fact was already implicit on the work of Barrett \cite{BAR2}, and
it appears explicitely in \cite{CHS} Prop. 6.5.5.
\bigskip   \\

\section{Boundedness of the Bergman Projection on $D'_\b$} 
\label{section-boundedness-D-prime-beta}
\vspace*{.12in}

We now come to the proof of Theorem \ref{THM4}.  By the asymptotic
expansion of the Bergman kernel of $D'_\b$ given in 
Theorem \ref{THM1} and the discussion in the previous section we have
that
\begin{align}
& |K_{D'_\b}(z,w)| \notag\\
& \le C e^{-h|\Re z_1-\Re w_1|}
\Bigl( |K_1(z,w)|+ \cdots
+|K_6(z,w)| |\bigr) \notag\\
& \qquad+ C e^{-\nu_\b |\Re z_1-\Re w_1|}
\biggl(\big|e^{[\pi- i(z_1- \ovw_1)]/2}- z_2 \ovw_2\big|^{-2} 
   + \big|e^{-[i(z_1- \ovw_1)+\pi]/2}-z_2 \ovw_2\big|^{-2} 
\biggr) \notag\\
& \equiv B_1(z,w)+\cdots +B_8(z,w) \label{B1-B8} \, ,
\end{align}
where $h$ is as in (\ref{(*)})  of Theorem \ref{THM1}.
Therefore, in order to prove the boundedness of the Bergman
projection $P_{D'_\b}$, it suffices to prove the boundedness of the
operators $T_{B_1},\dots, T_{B_8}$ with positive kernels
$B_1,\dots,B_8$, respectively.

We will show that the operators $T_{B_1},T_{B_2}, T_{B_3}$ and
$T_{B_7}$ and $T_{B_8}$ are bounded on $L^p(D'_\b)$, for
$1<p<\infty$. The boundedness of the remaining operators will follow
by completely analogous arguments because of the symmetries discussed
in the previous section.

In the course of the proof of Theorem \ref{THM4} and 
of the next few lemmas, we will make use of the classical Schur's lemma
(for which see, e.g., \cite{Zhu}, or \cite{SAD}).  We will also use the
following estimates of ``Forelli-Rudin'' type (see \cite{FORR}).
\begin{lemma}\label{disk-estimates}
{\bf (i)} Let $0<R<Q$, $0<p<1$ and $q>0$. Let $\l,\tau$ denote (single)
complex variables.  Then there exists a constant $C>0$ such that,
for all $|\tau|<R$, we have
\begin{multline*}
\iint_{|\l|<R} 
\frac{1}{(Q^2-|\l|^2)^{q}}\frac{1}{(R^2-|\l|^2)^{p}}
\frac{1}{|Q^2-\tau \overline{\l}|^2} dV(\l) \\
\le C
\min \left(
\frac{1}{(Q^2-R^2)^q} \frac{1}{(1-|\tau|^2)^p},
\frac{1}{(1-|\tau|^2)^{p+q}} \right) 
\
.
\end{multline*}

{\bf (ii)} Moreover, if $a>0$, then there exists a constant $C>0$ such
that 
$$
\int_{-\infty}^{+\infty} \frac{e^{-h|x|}(1+|x|)}{x^2+a^2}\, dx \le C \cdot \frac1a
$$
as $a\ra0^+$.  Here $h$ is a fixed positive constant.
\end{lemma}
\proof Fact {\bf (i)} is an immediate consequence of the classical 
Forelli-Rudin type estimates
(e.g. see \cite{Zhu}), 
since $(Q^2-|\l|^2)^{-q} \le (Q^2-R^2)^{-q}$.

Fact {\bf (ii)} follows by the change of variables $x=ay$ and a simple
pointwise estimate.
We leave the elementary details to the reader. 
\endpf

\begin{proposition}\label{prop-B1}     \sl
For $z,w\in D'_\b$ let 
$$
B_1(z,w)
 =  \frac{e^{-h{|\Re z_1-\Re w_1|}}\bigl(1+|\Re z_1 -\Re w_1|\bigr)}{
\big|i(z_1-\ovw_1)+2\b\big|^2\big|e^{\b-\pi/2}  -z_2\ovw_2\big|^2} \, ,
$$
where $h$ is defined by (\ref{(*)}), 
and let $T=T_{B_1}$ be the integral operator
$$
T f(z) \equiv T_{B_1} f(z) =\int_{D'_\b} B_1(z,w) f(w) dV(w)\ .
$$
Then $T:L^p(D'_\b)\ra L^p(D'_\b)$ for $1<p<\infty$.
\end{proposition}
\proof
Let $a>0$ be a number to be specified later and define
\begin{equation}\label{phi-prop-B1}
\varphi (w) =
\bigl[ ( e^{\b-\pi/2}-|w_2|^2)( \b-\Im w_1)\bigr]^{-a} \ .
\end{equation}
We wish to show (following the paradigm of Schur's lemma) 
that there exists a constant $C>0$ such  that 
\begin{equation}\label{Schur1}
\int_{D'_\b} B_1(z,w) \varphi^p(w) dV(w) \le C \varphi^p (z)
\end{equation}
and
\begin{equation}\label{Schur2}
\int_{D'_\b} B_1(z,w) \varphi^{p'}(w) dV(w) \le C \varphi^{p'} (z)
\end{equation}
for all $z\in D'_\b$, where $p'$ is the exponent conjugate to $p$.

We write
$w_1 = t+iu$, 
and we break the region of integration into three parts:
\begin{align*}
E_1 & =D'_\b\cap \bigl\{ w:\, -\b< u<-\b+\pi \bigr\}\, , \\
E_2 & =D'_\b\cap \bigl\{ w:\, -\b+\pi < u< \b-\pi
\bigr\}\, , \\
E_3 & =D'_\b\cap \bigl\{ w:\, \b-\pi < u< \b \bigr\} \ .
\end{align*}

Let $I_1,\, I_2$, and $I_3$ respectively denote the integrals over
$E_1,E_2$ and 
$E_3$. 
Then, applying Lemma \ref{disk-estimates}, we see that
\begin{align*}
I_1 & = \int_{E_1} B_1(z,w) \varphi^p(w) dV(w) \\
& \le \int_{-\infty} ^{+\infty} \int_{-\b}^{-\b+\pi} 
\frac{e^{-h(|\Re z_1-t|)}
\bigl(1+|\Re z_1 -t|\bigr)}{
 (\Re z_1-t)^2+(2\b-\Im z_1 -u)^2 } \frac{1}{(\b-u)^{ap}} \\
& \qquad\times
 \iint_{e^{[-\b+\pi/2]/2}<|w_2|< e^{[u+\pi/2]/2}} 
\frac{1}{
\bigl( e^{\b-\pi/2}-|w_2|^2\bigr)^{ap} 
\big|e^{\b-\pi/2}  -z_2\ovw_2\big|^2}
\, dV(w_2)\, dudt \\
&\le C \frac{1}{(e^{\b-\pi/2}  -|z_2|^2)^{ap}}  
\int_{-\infty} ^{+\infty} \int_{-\b}^{-\b+\pi} 
\frac{e^{-h(|\Re z_1-t|)}
\bigl(1+|\Re z_1-t|\bigr)}{
 (\Re z_1-t)^2+(2\b-\Im z_1 -u)^2 } \frac{1}{(\b-u)^{ap}} \, 
dudt \\
& \le C \frac{1}{(e^{\b-\pi/2}  -|z_2|^2)^{ap}}  \\
& \le C \varphi^p (z)\ .
\end{align*}

The estimate for $I_2$ is similar. 
Applying Lemma \ref{disk-estimates} {\bf (ii)}
we have
\begin{align*}
I_2 & = \int_{E_2} B_1(z,w) \varphi^p(w) dV(w) \\
& \le \int_{-\infty} ^{+\infty} \int_{-\b+\pi}^{\b-\pi} 
\frac{e^{-h(|\Re z_1-t|)}
\bigl(1+|\Re z_1 -t|\bigr)}{
 (\Re z_1-t)^2+(2\b-\Im z_1 -u)^2 } \frac{1}{(\b-u)^{ap}} \\
& \qquad\times
 \iint_{e^{[-u+\pi/2]/2}< |w_2|< e^{[u+\pi/2]/2}} 
\frac{1}{
\bigl( e^{\b-\pi/2}-|w_2|^2\bigr)^{ap} 
\big|e^{\b-\pi/2}  -z_2\ovw_2\big|^2}
\, dV(w_2)\, dudt \\
&\le C \frac{1}{(e^{\b-\pi/2}  -|z_2|^2)^{ap}}  
\int_{-\infty} ^{+\infty} \int_{-\b+\pi}^{\b-\pi} 
\frac{e^{-h(|\Re z_1-t|)}
\bigl(1+|\Re z_1-t|\bigr)}{
 (\Re z_1-t)^2+(2\b-\Im z_1 -u)^2 } \frac{1}{(\b-u)^{ap}} 
\, dudt \\
& \le C \frac{1}{(e^{\b-\pi/2}  -|z_2|^2)^{ap}}  \\
& \le C \varphi^p (z)\ .
\end{align*}

Finally, we estimate $I_3$.   For $0<a<1/p$ we have
\begin{align}
I_3 & = \int_{E_3} B_1(z,w) \varphi^p(w) dV(w) \notag \\
& \le \int_{-\infty} ^{+\infty} \int_{\b-\pi}^{\b} 
\frac{e^{-h(|\Re z_1-t|)}
\bigl(1+|\Re z_1 -t|\bigr)}{
 (\Re z_1-t)^2+(2\b-\Im z_1 -u)^2 } \frac{1}{(\b-u)^{ap}} \notag \\
& \qquad\times
 \iint_{e^{[-u+\pi/2]/2}< |w_2|< e^{[\b-\pi/2]/2}} 
\frac{1}{
\bigl( e^{\b-\pi/2}-|w_2|^2\bigr)^{ap} 
\big|e^{\b-\pi/2}  -z_2\ovw_2\big|^2}
\, dV(w_2)\, dudt \notag \\
 &\le C \frac{1}{(e^{\b-\pi/2}-|z_2|^2) ^{ap} }
\int_{-\infty} ^{+\infty} \int_{\b-\pi}^{\b} 
\frac{e^{-h(|\Re z_1-t|)}
\bigl(1+|\Re z_1-t|\bigr)}{
 (\Re z_1-t)^2+(2\b-\Im z_1 -u)^2 } 
\frac{1}{(\b-u)^{ap}} \, dudt \notag \\
& \le C \frac{1}{(e^{\b-\pi/2}-|z_2|^2) ^{ap} }
\int_{\b-\pi}^{\b} 
\frac{1}{2\b-(\Im z_1 +u) } 
\frac{1}{(\b-u) ^{ap} }\, du \ .\label{line-integral}
\end{align}

Now, by a simple change of variables, we see that 
the last integral above equals
\begin{align*}
\int_0^\pi\frac{1}{v+\b-\Im z_1}\frac{1}{(1-e^{-v})^{ap}}\, dv
& \le C \int_0^\pi\frac{1}{v+\b-\Im z_1}\frac{1}{v^{ap}}\, dv\\
& \le C \frac{1}{(\b-\Im z_1)^{ap}} \ ,
\end{align*}
again, 
as long as $0<ap<1$.  Inserting this estimate into
(\ref{line-integral}) we obtain that 
$I_3 \le C \varphi^p(z)$, provided that $0<a<1/p$.

Now, we may repeat the previous argument, with $p$ replaced by $p'$,
to obtain that 
$$
\int_{D'_\b} B_1(z,w) \varphi^{p'}(w) dV(w) \le C \varphi^{p'} (z)
$$
provided that  $0<a<1/p'$.

Choosing $a$ such that $0<a<\min(1/p,1/p')$ we can find $\varphi$ 
so that (\ref{Schur1}) and (\ref{Schur2}) are satisfied.  This
 concludes the proof.
\endpf

\begin{proposition}\label{prop-B2}   \sl
For $z,w\in D'_\b$ let 
$$
B_2(z,w)
 =  \frac{e^{-h{|\Re z_1-\Re w_1|}}\bigl(1+|\Re z_1 -\Re w_1|\bigr)}{
\big|i(z_1-\ovw_1)+2\b\big|^2
\big|e^{[-i(z_1-\ovw_1)+\pi]/2}  -z_2\ovw_2\big|^2}  \, ,
$$
where $h$ is as specified in (\ref{(*)}),
and let $T=T_{B_2}$ be the integral operator
$$
T f(z) =\int_{D'_\b} B_2(z,w) f(w) dV(w)\ .
$$
Then $T:L^p(D'_\b)\ra L^p(D'_\b)$ for $1<p<\infty$.
\end{proposition}
\proof
The proof follows the same lines as the previous one.

For $0<a<\min(1/p,1/p')$ we now define
\begin{equation}\label{phi-prop-B2}
\varphi (w) =
\bigl[ ( e^{[\Im w_1 +\pi/2]}-|w_2|^2)( \b-\Im w_1)\bigr]^{-a} 
\end{equation}
and, again, we wish to show there exists a constant $C>0$ such  that 
$$
\int_{D'_\b} B_2(z,w) \varphi^p(w) dV(w) \le C \varphi^p (z)
$$
and
$$
\int_{D'_\b} B_2(z,w) \varphi^{p'}(w) dV(w) \le C \varphi^{p'} (z)
$$
for all $z\in D'_\b$.

We consider the former integral first.
Once more we write $w_1=t+iu$,  break the region of integration into
$E_1$, $E_2$ and $E_3$, defined as in the proof of
Proposition \ref{prop-B1}, and call the respective integrals $I_1$, $I_2$ and
$I_3$. 

We begin with $I_3$.
Using Lemma \ref{disk-estimates}, {\bf (i)} and {\bf (ii)}, we have
\begin{align*}
I_3 
& \le C 
\int_{-\infty} ^{+\infty} \int_{\b-\pi}^{\b} 
\frac{e^{-h(|\Re z_1-t|)}
\bigl(1+|\Re z_1 -t|\bigr)}{
 (\Re z_1-t)^2+(2\b-\Im z_1 -u)^2 } \frac{1}{(\b-u)^{ap}} \notag \\
& \qquad\times
 \iint_{e^{[-u+\pi/2]/2}< |w_2|< e^{[\b-\pi/2]/2}} 
\frac{1}{
\bigl( e^{[u+\pi/2]}-|w_2|^2\bigr)^{ap} 
\big|e^{\b-\pi/2}  -z_2\ovw_2\big|^2}
\, dV(w_2)\, dudt \notag \\
&  \le C
\frac{1}{
\bigl( e^{\b-\pi/2}-|z_2|^2\bigr)^{ap} }
\int_{-\infty} ^{+\infty} \int_{\b-\pi}^{\b} 
\frac{e^{-h(|\Re z_1-t|)}
\bigl(1+|\Re z_1 -t|\bigr)}{
 (\Re z_1-t)^2+(2\b-\Im z_1 -u)^2 } \frac{1}{(\b-u)^{ap}} \, dudt\\ 
&  \le C
\frac{1}{
\bigl( e^{\b-\pi/2}-|z_2|^2\bigr)^{ap} }
\int_{\b-\pi}^{\b} 
\frac{1}{(2\b-\Im z_1 -u)(\b-u)^{ap}} \, du\\
 &  \le C
\frac{1}{
\bigl( e^{\b-\pi/2}-|z_2|^2\bigr)^{ap} }\cdot
\frac{1}{(\b-\Im z_1)^{ap}} \\
& \le C \varphi^p(z)\ .
\end{align*}

Next we turn to $I_2$. We have that
\begin{align*}
I_2 
& \le \int_{-\infty} ^{+\infty} \int_{-\b+\pi}^{\b-\pi} 
\frac{e^{-h(|\Re z_1-t|)}
\bigl(1+|\Re z_1 -t|\bigr)}{
 (\Re z_1-t)^2+(2\b-\Im z_1 -u)^2 } \frac{1}{(\b-u)^{ap}} \\
& \qquad\times
 \iint_{e^{[-u+\pi/2]/2}< |w_2|< e^{[u+\pi/2]/2}} 
\frac{1}{
\bigl( e^{[u+\pi/2]}-|w_2|^2\bigr)^{ap} 
\big|e^{\b-\pi/2}  -z_2\ovw_2\big|^2}
\, dV(w_2)\, dudt \\
& \le C
\frac{1}{
\bigl( e^{\b-\pi/2}-|z_2|^2\bigr)^{ap} }
\int_{-\infty} ^{+\infty} \int_{-\b+\pi}^{\b-\pi} 
\frac{e^{-h(|\Re z_1-t|)}
\bigl(1+|\Re z_1 -t|\bigr)}{
 (\Re z_1-t)^2+(2\b-\Im z_1 -u)^2 } \frac{1}{(\b-u)^{ap}} \, dudt\\ 
& \le C
\frac{1}{
\bigl( e^{\b-\pi/2}-|z_2|^2\bigr)^{ap} } \\
& \le C \varphi^p(z)\ .
\end{align*}

Finally we estimate $I_1$, where again we use the fact that the
integrals in $u$ and $t$ are bounded uniformly in $z_1$:
\begin{align*}
I_1
& \le C 
\int_{-\infty} ^{+\infty} \int_{-\b}^{-\b+\pi} 
\frac{e^{-h(|\Re z_1-t|)} \bigl(1+|\Re z_1 -t|\bigr)}{
 (\Re z_1-t)^2+(2\b-\Im z_1 -u)^2 } \frac{1}{(\b-u)^{ap}} \notag \\
& \qquad\times
 \iint_{e^{[-\b+\pi/2]/2}<|w_2|< e^{[u+\pi/2]/2}} 
\frac{1}{
\bigl( e^{[u+\pi/2]}-|w_2|^2\bigr)^{ap} 
\big|e^{\b-\pi/2}  -z_2\ovw_2\big|^2}
\, dV(w_2)\, dudt \notag \\
& \le C\frac{1}{(e^{\b-\pi/2} -|z_1|^2)^{ap}}\\
& \le C \varphi^p(z)\ . 
\end{align*}

Again, we may repeat the argument with $p$ replaced by $p'$.
This concludes the proof.
\endpf

In the proof of the next proposition we shall again need the following 
``Forelli-Rudin type estimate'' (see also Lemma 31. above). 
\begin{lemma}\label{new-Forelli-Rudin-estimate}   \sl
{\bf (i)} Let $0<R<Q$, $\varrho>0$ and let
$m=\min(\varrho,Q)$. Let $\l, \tau$ denote (single) complex
variables with $|\tau|<m$ and
$|\l|<R$, and $\theta\in\mathbf R$.  Then, for
$0<a,b<1$, there exists 
a constant $C>0$ such that
\begin{multline*}
\iint_{|\l|<R}\frac{1}{|R\varrho-e^{i\theta}\tau\bar\l|^2}\cdot
\frac{1}{|Q^2-\tau\bar\l|^2}\cdot
\frac{1}{(R^2-|\l|^2)^{a}(Q^2-|\l|^2)^b}\, dV(\l) 
\\ \le C \frac{1}{(Q^2-R^2)^{b}}\cdot \frac{1}{(m^2-|\tau|^2)^a}
\cdot \frac{1}{\big|Q^2-e^{i\theta}\frac{R}{\varrho} |\tau|^2\big|^2}
\end{multline*}
as $|\l|\ra R^{-}$, $R\ra S^{-}$, and $|\tau|\ra\varrho^{-}$.

{\bf (ii)} Let  $\delta>0$ and 
$0<b<1$.  Then there exists a constant
$C>0$ such that
$$
\int_{-\infty}^{+\infty}
\frac{e^{-h|x|}(1+|x|)}{|1-be^{ix}|^{1+\delta}} \le
C\frac{1}{(1-b)^\delta} \ .
$$
\end{lemma}
\proof Although this is a fairly straightforward proof, we sketch it here
for completeness and for the reader's convenience.

In order to prove part {\bf (i)} we claim that, for
$0<a<1$, there exists a constant $C>0$
such that, for all $\tau,\tau'$ in the unit disk,
\begin{multline*}
\iint_{|\l|<1} \frac{1}{|1-\tau'\bar\l|^2\,|1-z\bar\l|^2\, (1-|\l|^2)^a}\,
dV(\l) \\
\le C \frac{1}{|1-\tau'\bar\tau|^2}
\cdot \min\bigl\{ (1-|\tau'|^2)^{-a},\,  (1-|\tau|^2)^{-a} \bigr\} \, .
\end{multline*}
Given the claim, part {\bf (i)} follows by noticing that 
$(Q^2-|\l|^2)^{-b} \le (Q^2-R^2)^{-b}$ and rescaling.

In order to prove the claim, one first integrates over the region
$E\equiv
\bigl\{ \l:\, |1-\tau'\bar\l|\le \gamma |1-\tau\bar\l| \bigr\}$, for
some $\gamma>0$ to be chosen later.  Then one integrates over the
symmetric region 
$E'\equiv
\bigl\{ \l:\, |1-\tau\bar\l|\le \gamma |1-\tau'\bar\l| \bigr\}$, and then
on the unit disk with $E\cup E'$ removed.

The proof of {\bf (ii)} is achieved by noticing that the integral is controlled
by a constant times
$$
\int_0^1 \frac{1}{(1-b)^{1+\delta} +x^{1+\delta}}\, dx\le C
\frac{1}{(1-b)^\delta} \ . 
$$
This concludes the (sketch of the) proof of the lemma.
\endpf	
\medskip \\

\begin{proposition}\label{prop-B3}   \sl
For $z,w\in D'_\b$ let 
$$
B_3(z,w)
 =  \frac{e^{-h{|\Re z_1-\Re w_1|}}\bigl(1+|\Re z_1 -\Re w_1|\bigr)}{
\big|e^{[\pi-i(z_1- \ovw_1)]/2}-z_2 \ovw_2\big|^2\,
\big|e^{\bmp}-z_2 \ovw_2\big|^2}  \, ,
$$				      
where $h$ is defined in (\ref{(*)}), 
and let $T=T_{B_3}$ be the integral operator
$$				       
T_{B_3} f(z) =\int_{D'_\b} B_3(z,w) f(w) dV(w)\ .
$$
Then $T:L^p(D'_\b)\ra L^p(D'_\b)$ for $1<p<\infty$.
\end{proposition}
\proof
For $0<a<\min(1/p,1\p')$ we now define
\begin{equation}\label{phi-prop-B3}
\varphi (w) =
\bigl[ ( m(w)^2
-|w_2|^2)(e^{\b-\pi/2}-|w_2|^2)\bigr]^{-a}
\end{equation}
where $m(w)\equiv \min\bigl( e^{[\b-\pi/2]/2}, e^{[\Im w_1+\pi/2]/2}\bigr)$.
Again, we wish to show there exists a constant $C>0$ such  that 
$$
\int_{D'_\b} B_3(z,w) \varphi^p(w) dV(w) \le C \varphi^p (z)
$$
and
$$
\int_{D'_\b} B_3(z,w) \varphi^{p'}(w) dV(w) \le C \varphi^{p'} (z)
$$
for all $z\in D'_\b$.

Once more we write $w_1=t+iu$,  break the region of integration into
$E_1$, $E_2$ and $E_3$, and call the respective integrals $I_1$, $I_2$
and $I_3$. 

Then we have,
\begin{multline*}
I_1 
 \le C 
\int_{-\infty} ^{+\infty} \int_{-\b}^{-\b+\pi} 
\iint_{e^{[-\b+\pi/2]/2}< |w_2|<  e^{[u+\pi/2]/2}}  
\frac{
e^{-h(|\Re z_1-t|)}
\bigl(1+|\Re z_1 -t|\bigr)}{
\big|e^{[\pi+\Im z_1+u-i(\Re z_1-t)]/2}-z_2 \ovw_2\big|^2}
\\
\times
 \frac{1}{\big|e^{\b-\pi/2}  -z_2\ovw_2\big|^2}
\cdot
 \frac{1}{
\bigl( e^{[u+\pi/2]}-|w_2|^2\bigr)^{ap}
\bigl(e^{\b-\pi/2}-|w_2|^2\bigr)^{ap} }
\, dV(w_2)\, dudt \ .
\end{multline*}

Now we change variables in the outer integral, setting $s=t-\Re z_1$,
and we apply Lemma \ref{new-Forelli-Rudin-estimate} {\bf (i)} with
$e^{[u+\pi/2]/2}=R$, $e^{[\b-\pi/2]/2}=Q$, $e^{[\Im z_1+\pi/2]/2}=\varrho$:
\begin{align*} 
I_1 & \le
\int_{-\infty} ^{+\infty} 
\int_{-\b}^{-\b+\pi} 
\iint_{e^{[-\b+\pi/2]/2}< |w_2|<  e^{[u+\pi/2]/2}}  
\frac{e^{-h(|s|)}(1+|s|)}{
\big|e^{[\pi+\Im z_1+u]/2}-e^{-is/2}z_2 \ovw_2\big|^2}
\\
& \qquad \times
 \frac{1}{\big|e^{\b-\pi/2}  -z_2\ovw_2\big|^2}
\cdot
 \frac{1}{
\bigl( e^{[u+\pi/2]}-|w_2|^2\bigr)^{ap}
\bigl(e^{\b-\pi/2}-|w_2|^2\bigr)^{ap} }
\, dV(w_2)\, duds \\
& \le C \frac{1}{(e^{[\Im z_1+\pi/2]}-|z_2|^2)^{ap} }
\int_{-\infty} ^{+\infty} \int_{-\b}^{-\b+\pi}  
\frac{e^{-h(|s|)}(1+|s|)}{ 
\big| e^{\b-\pi/2} -
e^{-is/2} e^{[u -\Im z_1]/2} |z_2|^2\big|^2} \\
& \qquad \times \frac{1}{
\bigl(e^{\b-\pi/2}- e^{[u+\pi/2]} \bigr)^{ap} } \, duds  \\
 & \le C \frac{1}{(e^{[\Im z_1+\pi/2]}-|z_2|^2)^{ap} } \\
& \qquad\times
\int_{-\b}^{-\b+\pi} 
\frac{1}{\bigl( e^{\b-\pi/2} e^{[\Im z_1+\pi/2]/2} - e^{[u+\pi/2]/2}
  |z_2|^2\bigr) 
(e^{\b-\pi/2}- e^{[u+\pi/2]})^{ap} } \, du \ .
\end{align*}
Here we have applied Lemma \ref{new-Forelli-Rudin-estimate} {\bf (ii)} 
to the outer integral in $ds$.  Now it is easy to see that the
integral in $du$ is uniformly bounded so that 
$$
I_1\le C  \frac{1}{(e^{[\Im z_1+\pi/2]}-|z_2|^2)^{ap} } 
\le C \varphi^p(z)\ .
$$

The proof of the estimate for $I_2$ is identical to the one for $I_1$
and we skip the details.  

Finally, for $I_3$ we begin by arguing as in
the previous cases to obtain that
\begin{multline*}
I_3
\le
\int_{-\infty} ^{+\infty} e^{-h(|s|)}(1+|s|)
\int_{\b-\pi}^{\b} 
\iint_{e^{[-\b+\pi/2]/2}< |w_2|<  e^{[\b-\pi/2]/2}}  
\frac{1}{
\big|e^{[\pi+\Im z_1+u]/2}-e^{-is/2}z_2 \ovw_2\big|^2}
\\
\times
 \frac{1}{\big|e^{\b-\pi/2}  -z_2\ovw_2\big|^2}
\cdot
 \frac{1}{
\bigl( e^{[u+\pi/2]}-|w_2|^2\bigr)^{ap}
\bigl(e^{\b-\pi/2}-|w_2|^2\bigr)^{ap} }
\, dV(w_2)\, duds \ .
\end{multline*}
Again we apply first
Lemma \ref{new-Forelli-Rudin-estimate} {\bf (i)}, this time
with 
$e^{[\b-\pi/2]/2}=R$, $e^{[u+\pi/2]/2}=Q$, 
and
$e^{[\Im z_1+\pi/2]/2}=\varrho$, and then part {\bf (ii)} of the
same lemma. We obtain
\begin{align*}
I_3
& \le C \frac{1}{(e^{\b-\pi/2}-|z_2|^2)^{ap} }
\int_{-\infty} ^{+\infty} \int_{\b-\pi}^{\b}  
\frac{e^{-h(|s|)}(1+|s|)}{ 
\big| e^{[u+\pi/2]}
-e^{-is/2} e^{-[\Im z_1+\pi/2]/2} e^{(\b-\pi/2)/2} 
  |z_2|^2\big|^2} \\
& \qquad \times \frac{1}{
\bigl(e^{\b-\pi/2}- e^{[u+\pi/2]} \bigr)^{ap} } \, duds  \\
 & \le C \frac{1}{(e^{\b-\pi/2}-|z_2|^2)^{ap} } \\
& \qquad\times
\int_{\b-\pi}^{\b} 
\frac{1}{\bigl( e^{[u+\pi/2]}
-e^{-[\Im z_1+\pi/2]/2} e^{(\b-\pi/2)/2} 
  |z_2|^2\bigr) 
(e^{[u+\pi/2]}-e^{\b-\pi/2})^{ap} } \, du \ .
\end{align*}

Now we recall that, for $y>0$, if $0<ap<1$,  
$$
\int_0^{+\infty}\frac{1}{x^{ap}}\cdot \frac{1}{x+y}\, dy \le
C\frac{1}{y^{ap}} 
$$
as $y\ra0^+$.  Then the last integral in the display above can be
easily seen to be less or equal to a constant times
$
\bigl( e^{[\Im z_1+\pi/2]/2} e^{(\b-\pi/2)/2} 
  -|z_2|^2 \bigr)^{-ap}$.
Therefore it follows that
$$
I_3 
\le C\frac{1}{(e^{\b-\pi/2}-|z_2|^2)^{ap} }
\cdot \frac{1}{
\bigl( e^{[\Im z_1+\pi/2]/2} e^{(\b-\pi/2)/2} 
  -|z_2|^2 \bigr)^{ap}} 
 \le C \varphi^p(z)\ ,
$$
as it is easy to check.

Again, we may repeat the argument with $p$ replaced by $p'$ and
then obtain that $T:L^p(D'_\b)\ra L^p(D'_\b)$ is bounded for
$1<p<\infty$. 
\endpf
\medskip  \\

\proof[End of the Proof of Theorem \ref{THM4}]
In order to conclude the proof of Theorem \ref{THM4} we only need to
show that the operators $T_{B_7}$ and $T_{B_8}$ are bounded on
$L^p(D'_\b)$, for $1<p<\infty$.

The estimates for these kernels are simpler than the ones in
\ref{prop-B1}--\ref{prop-B3}.
In these last two cases it
suffices to consider test functions $\varphi$ equal to
$$
(e^{[\Im w_1 +\pi/2]}- |w_2|^{2})^{-a}\qquad\text{and}\quad
   (|w_2|^2 -e^{[\Im w_1 -\pi/2]})^{-a}
$$
respectively. 

To avoid further repetitions,
we shall leave the simple details to the reader.
\endpf
\vspace*{.12in}

\section{Boundedness of the Bergman Projection on $D_\b$} 
\label{section-boundedness-D-beta} 
\vspace*{.12in}

We now turn our attention to the $L^p$ boundedness of the
Bergman projection operator on $D_\b$.  Even though
the domain $D_\b$ is biholomorphic to $D'_\b$,
the $L^p$ behavior of the Bergman projection on these two domains
is not {\it a priori} identical.  That is because the space
of $L^p$ holomorphic functions for $p \ne 2$ does not transform
canonically under the 
biholomorphic mapping. 
Put in other words, the Jacobian of the mapping plays a non-obvious
role in the calculations.  And in fact it turns out that
(in contrast to the situation on $D'_\b$) 
the projection $P_{D_\b}$ is bounded on $L^p(D_\b)$ only for the
restricted range $2/(1 + \nu_\b) < p < 2/(1 - \nu_\b)$.
 We shall now establish this last
assertion.

We shall use our asymptotic expansion for the Bergman
kernel to determine the $L^p$ boundedness of the Bergman projection
on the domain $D_\b$.
We remark that, in his review \cite{Str}
of Barrett's paper \cite{BAR2}, E. Straube
describes some results that are related to those that are presented
here. 

We begin with the negative result.
Let $1 < p < \infty$ and assume that $P:L^p(D_\b) \rightarrow
L^p(D_\b)$ is bounded. It follows that
for any $\z\in D_\b$ fixed, 
$K_{D_\b} (\cdot,\z) \in L^{p'}(D_\b)$, where
$p' = p/[p-1]$ is the exponent conjugate to $p$.  

For, if $P=P_{D_\b}$ is bounded, then
for all $f\in L^p(D_{\b})$ and all $\z\in D_\b$,
$$ 
|\langle f, K_{D_\b}(\cdot, \z) \rangle | =  |Pf(\z)| 
\leq  c_\z \|Pf\|_{L^p} \le C \|f\|_{L^p}\, .
$$

\begin{lemma}\label{unboundedness}	\sl
For any $\z\in D_\b$ it holds that
$K_{D_\b}(\cdot,\z)\in L^p(D_\b)$ only if
$2/(1 + \nu_\b) < p < 2/(1 - \nu_\b)$.
\end{lemma}
\proof
Fix $\z\in D_\b$ and define
$$
\Omega_\z = \bigl\{\om \in D_\b: |\om_1| < |\z_1| \ , \
1/4 \leq | e^{\pi/2} (\z_1/\overline{\om}_1 )^{\pm i/2} - 
\z_2 \overline{\om}_2 |
          \leq 1/2\bigr\} \, .
$$
Recall the expansion for the kernel $K_{D_\b}(\z,\om)$ given in
Theorem \ref{THM2}.  Then, for $\om\in\Omega_\z$,
we have that 
$$
|H_b(\z,\om)|,\, |H_{\tilde{b}}(\z,\om)|\le C_\z
$$
for some constant independent of $\om$, so that
\begin{equation}\label{estimate-from-below-K-D-beta}
|K_{D_\b}(\z,\om)|
\ge c_\z |\om_1|^{\nu_\b -1}
\end{equation}
for $\om\in\Omega_\z$.  

Therefore
\begin{align*}
\int_{D_{\b}} |K_{D_\b} (\cdot,\z)|^{p'} \, dV(\om) 
& \ge 
    \int_{\Omega_\z} |K_{D_\b} (\cdot,\z)|^{p'} \, dV(\om)\\
& \ge  c_\z 
\int_{\Omega_\z} \bigl(|\om_1|^{\nu_\b-1} 
\bigr)^{p'} \, dV(\om) 
 = c \int_0^{|\z_1|} r^{p'(\nu_\b-1) +1} \, dr \ .
\end{align*}

Obviously for convergence we need $p' (\nu_\b -1)+1>-1 $, that is
$p' < 2/[1 -\nu_\b]$.  Hence if
$p \geq 2/[1 - \nu_\b]$ then the integral diverges.  The other
result, for $p \leq 2/[1 + \nu_\b]$, follows by duality.
This proves the lemma.
\endpf

For the record, we record now the fact that we can use
Barrett's exhaustion procedure (see \cite{BAR2}) to obtain a
negative result with the same indices on the smooth worm
$\mathcal W$. Details of the proof appear in \cite{KP2}.

\begin{theorem}\label{straube-thm}	 \sl
Let $\mathcal P$ denote the Bergman projection on the smooth, bounded
worm $\mathcal W=\mathcal W_\b$, with $\b>\pi/2$.  Then, if $\mathcal P:
L^p(\mathcal W_\b)\ra L^p(\mathcal W_\b)$ is bounded, necessarily
$2/[1+\nu_\b]<p<2/[1-\nu_\b]$. 
\end{theorem}
\vspace*{.12in}

Next we prove the positive part of Theorem \ref{THM3}.

\begin{theorem}\label{positive-part-THM3}    \sl
Let 
$2/[1+\nu_\b]<p<2/[1-\nu_\b]$. Then
$$
P_{D_\b} : L^p(D_\b)\ra L^p(D_\b)
$$ 
is bounded.
\end{theorem}

In order to prove this fact we use the 
biholomorphism $\Phi: D'_\beta \ra D_\beta$
given by $\Phi(z_1, z_2) = (e^{z_1}, z_2)$, so that
$\det \Jac_\CC \Phi = e^{z_1}$.

Therefore, by the transformation rule for the Bergman kernels {\em
  via} 
biholomorphic mappings,
we  wish to show that the integral operator $T$
with kernel 
$L(z,w)\equiv
e^{-z_1}K_{D'_\b}(z,w)e^{\ovw_1}$
is bounded on $L^p(D'_\b)$ when $2/[1+\nu_\b]<p<2/[1-\nu_\b]$.

As in the proof of Theorem \ref{THM4}, we are going to use  Schur's
lemma.  With the notation of Section
\ref{section-boundedness-D-prime-beta}, recalling (\ref{B1-B8}), we
have that
\begin{align}
|L(z,w)|
& \le C e^{\Re w_1 -\Re z_1} |K_{D'_\b}(z,w)| \le C e^{\Re w_1 -\Re z_1} 
\bigl( B_1(z,w)+\cdots B_8(z,w) \bigr) \notag\\
& \equiv A_1(z,w)+\cdots A_8(z,w)\ .
\label{A1-A8}
\end{align}

For each $j=1,\dots,8$, we
wish to determine a positive function $\psi$ on $D'_\b$ such that
$$
\int_{D'_\b} A_j (z,w) \psi^p(w)\, dV(w)
\le C \psi^p(z)\ . 
$$

As in the proof of Theorem \ref{THM4}, by symmetry
it suffices to consider the
cases of $j=1,2,3$ and $j=7$.  

\begin{proposition}\label{prop-A1}   \sl
For $z,w\in D'_\b$ let 
$$
A_1(z,w)
 =  \frac{e^{-h{|\Re z_1-\Re w_1|}} e^{\Re w_1 -\Re z_1}
\bigl(1+|\Re z_1 -\Re w_1|\bigr)}{
\big|i(z_1-\ovw_1)+2\b\big|^2\big|e^{\b-\pi/2}  -z_2\ovw_2\big|^2} \, ,
$$
for $h$ as in (\ref{(*)}),
and let $T$ be the integral operator
$$
T f(z) =\int_{D'_\b} A_1(z,w) f(w) dV(w)\ .
$$
Then $T:L^p(D'_\b)\ra L^p(D'_\b)$ for 
$2/[1+\nu_\b]<p<2/[1-\nu_\b]$. 
\end{proposition}
\proof
Let $a,b>0$ and define
\begin{equation}\label{psi-prop-L1}
\psi (w) = e^{-b\Re w_1} 
\bigl[ ( e^{\b-\pi/2}-|w_2|^2)( \b-\Im w_1)\bigr]^{-a} \ .
\end{equation}

Arguing as in the proof of Proposition \ref{prop-B1}, 
dividing the region of integration into $E_1\cup E_2\cup E_3$ and 
performing the
integrations in $w_2$ and $\Im w_1$ (having set
$t=\Re w_1$), we find that 
\begin{align*}
&\int_{E_1\cup E_2} A_1(z,w) \psi^p(w) dV(w) \\
& \qquad \le C \frac{1}{(e^{\b-\pi/2}  -|z_2|^2)^{ap}}  
\int_{-\infty} ^{+\infty} 
e^{-h(|\Re z_1-t|)} e^{t-\Re z_1}
\bigl(1+|\Re z_1-t|\bigr) e^{-bpt}\, 
dt \\
& \qquad  = C \frac{1}{(e^{\b-\pi/2}  -|z_2|^2)^{ap}}  
\int_{-\infty} ^{+\infty} 
e^{-h(|s|)} e^{s}
\bigl(1+|s|\bigr) e^{-bp(s+\Re z_1)}\, 
dt \\
& \qquad  = C \psi^p (z)
\end{align*}
as long as $0<1-bp<h$, that is $\frac1p -h<b<\frac1p$.

The same condition must be valid for the conjugate exponent $p'$,
that
is 
$\frac{1}{p'} -h<b<\frac{1}{p'}$, which translates into
$1-\frac1p -h<b<1-\frac1p$.  
Assume that $p>2$ for the moment.  
Then the intervals
$(\frac1p -h,\frac1p)$ and $(1-\frac1p -h,1-\frac1p)$ have non-empty
intersection if and only if 
$$
\frac1p > 1-\frac1p -h\, ;
$$
that is if and only if $h>1-\frac2p$.  

We recall here that $h$ is assumed to satisfy $\nu_b < h < \min(1,2\nu_b)$
(see Theorems \ref{THM1}, \ref{THM2}), and that $p$ is momentarily
assumed to satisfy 
$2<p< 2/[1-\nu_b]$. Then $1-\frac2p<\nu_b$ and we are done.

The argument for $p<2$ is similar, and we skip the details.
\medskip  \\
	       
Next, arguing again as in the proof of Proposition \ref{prop-B1},
we see that
\begin{align*}
&\int_{E_3} A_1(z,w) \psi^p(w) dV(w) \\
& \qquad \le C \frac{1}{(e^{\b-\pi/2}  -|z_2|^2)^{ap}}  
\int_{-\infty} ^{+\infty} \int_{\b-\pi}^{\b} 
\frac{
e^{-h(|\Re z_1-t|)} e^{t-\Re z_1}
\bigl(1+|\Re z_1-t|\bigr) e^{-bpt}}
{\bigl[ (\Re z_1-t)^2+(2\b-\Im z_1 -u)^2 \bigr](\b-u)^{ap}} 
\, dudt  \\
& \qquad  = C 
\frac{e^{-bp\Re z_1} }{(e^{\b-\pi/2}  -|z_2|^2)^{ap}}  
\int_{-\infty} ^{+\infty} \int_{\b-\pi}^{\b} 
\frac{
e^{-h(|s|)} e^{s(1-bp)}
\bigl(1+|s|\bigr) e^{-bp(s+\Re z_1)}}{
\bigl[ s^2+(2\b-\Im z_1 -u)^2 \bigr](\b-u)^{ap}} \, 
duds \\
& \qquad \le C \frac{e^{-bp\Re z_1} }{(e^{\b-\pi/2}-|z_2|^2) ^{ap} }
\int_{\b-\pi}^{\b} 
\frac{1}{2\b-(\Im z_1 +u) } 
\frac{1}{(\b-u) ^{ap} }\, du \\
& \qquad  
\le C \frac{e^{-bp\Re z_1} }{(e^{\b-\pi/2}-|z_2|^2) ^{ap} }
\frac{1}{(\b-\Im z_1)^{ap}} 
= C \psi^p (z)\ .
\end{align*}
This concludes the proof of the proposition.
\endpf

The fact that the
operators with integral kernels $A_2$ and $A_3$ are 
bounded on $L^p(D_\b)$ for $2/[1+\nu_\b]<p<2/[1-\nu_\b]$ 
is proved following the same pattern as in the proof of 
\ref{prop-B1}-\ref{prop-B3} and \ref{prop-A1}.
To avoid  repetitions we leave the details to the reader. 

\proof[End of the Proof of Theorem \ref{THM3}]
In order to finish the proof of Theorem \ref{THM3}, we show that
the integral operator with kernel $A_7(z,w)$ is bounded on
$L^p(D_\b)$ for $2/[1+\nu_\b]<p<2/[1-\nu_\b]$.
Here
$$
A_7(z,w)
\equiv \frac{e^{\nu_\b|\Re z_1 -\Re w_1|}(1+|\Re z_1-\Re w_1|) 
e^{\Re z_1+\Re w_1}}{
\big|e^{[-i(z_1-\ovw_1)+\pi]/2}  -z_2\ovw_2\big|^2}\, ,
$$
and we define
$$
\psi(w) = e^{-b\Re w_1} (e^{[\Im w_1 +\pi/2]}-|w_2|^2)^{-a} \, , 
$$
for some $a,b>0$ to be selected later.

Then, arguing as in the previous proof, we see that, for $0<a<1/p$,
\begin{align*}
& \int_{D_\b} A_7(z,w)\psi^p(w)\, dV(w)\\
& \qquad \le C \frac{1}{(e^{[\Im z_1 +\pi/2]}-|z_2|^2)^{-ap}}
\int_{-\infty}^{+\infty}
e^{-\nu_\b|\Re z_1 -t|}(1+|\Re z_1-t|) 
e^{\Re z_1+t} e^{-bpt}\, dt\\
 & \qquad 
\le C \frac{e^{-bp\Re z_1}}{(e^{[\Im z_1 +\pi/2]}-|z_2|^2)^{-ap}}
\int_{-\infty}^{+\infty}
e^{-\nu_\b|s|}(1+|s|) e^{(1-bp)s}\, ds \\
& \qquad C \psi^p(z) 
\end{align*}
as long as $0<1-bp<\nu_b$. 

The analogous condition with $p'$ in place
of $p$ (and the same $b$) must hold.
Again, we argue as in the proof of Proposition \ref{prop-A1} to see
that
$$
b \in \biggl(1-\frac1p,\frac1p-\nu_b\biggr)\quad\text{if}\ p<2
\qquad\text{and}\qquad
b\in \biggl(1-\frac1p-\nu_b,\frac1p\biggr) \quad\text{if}\ p>2
$$
satisfies the required conditions.

This concludes the proof of Theorem \ref{THM3}.
\endpf

\bigskip

\section{Decomposition of the Bergman Space}
\label{decomposition-B-space}
\vspace*{.12in}

In this section we begin our detailed analysis of the Bergman kernels
on $D_\beta$ and $D_{\beta'}$, leading to the asymptotic expansions in
Theorem \ref{THM1}
 and Theorem \ref{THM2}.  
We follow the calculations in \cite{KIS} and
\cite{BAR2} in order to obtain a decomposition of the Bergman space
on the domains $D_\b$ and $D_\b'$. We in fact concentrate our
attention on the latter domain. An analogous argument holds on
$D_\b$.  Note that in the present context the usual transformation
rule for the Bergman kernel will serve us well.  Calculations on
$D'_\beta$ will transfer to $D_\beta$ automatically.

Let the Bergman space $A^2(D'_\b)$ be the collection of
holomorphic functions that are square integrable with
respect to Lebesgue volume measure $dV$ on $D'_\b$.
Following Kiselman \cite{KIS} and Barrett \cite{BAR2}, 
we decompose the Bergman space as follows.
Using the rotational invariance in $z_2$ and 
elementary Fourier series, each $f \in A^2(D'_\b)$  can be  written as 
$$
f = \sum_{j=-\infty}^\infty f_j \, ,
$$
where each $f_j$ is holomorphic and satisfies 
$f_j(z_1, e^{i\theta} z_2) = e^{ij\theta} f(z_1, z_2)$ 
for $\theta$ real.  In fact, such an $f_j$ is given by
$$
f_j(z_1,z_2) = \frac1{2\pi} \int_0^{2\pi} f(z_1,e^{i\theta}z_2)
e^{-ij\theta} \,d\theta \, . 
$$
Therefore
$$
{\cal H} = \bigoplus_{j\in\ZZ} {\cal H}^j\, ,
$$
where
$$
{\cal H}^j = \bigl \{f \in L^2: f \ \hbox{is holomorphic and} \ 
          f(w_1, e^{i\theta} w_2) = e^{ij\theta} f(w_1, w_2) \bigr \} \, .
$$

Notice that the function 
$$
g_j(z_1,z_2)\equiv
\biggl[\frac1{2\pi} \int_0^{2\pi} f(z_1,e^{i\theta}z_2)
e^{-ij\theta} \,d\theta \biggr] z_2^{-j} 
$$
is holomorphic on $D'_\b$ and depends only on $|z_2|$.  Therefore it
must be locally constant in $z_2$.  But, for all $z_1$, the set
$\{ z_2:\, (z_1,z_2)\in D'_\b\}$ is connected, so that
$g_j(z_1,z_2)\equiv h_j(z_1)$.  Since $f_j(z_1,z_2)=h_j(z_1)z_2^j$ is
holomorphic on $D'_\b$, it is easy to see that $h_j$ must be
holomorphic on the strip $\{ z_1:\, |y|<\b\}$, where
$z_1=x+iy$.

The following result is  contained in
\cite{BAR2} Section 3.
(Further details can also be found in
\cite{KP2}.) 
\begin{lemma}\label{local-constantcy}	\sl
Let $\b>\pi/2$ and 
$f_j\in {\cal H}^j$.  Then there exists a function $h_j$
of one complex variable,
holomorphic in the strip $S_\b=\{ z_1=x+iy,:\, |y|<\b\}$, such that
$f_j(z_1,z_2)=h_j(z_1)z_2^j$.  Moreover, $h_j$ is 
square-integrable on $S_\b$
with respect to the weight
$$
\l_j(y) = \bigl( \chi_{\pi/2} * \bigl[
   e^{(j+1)( \, \cdot \, )} \chi_{\b - \pi/2} \bigr]
   \bigr)(y) 
\, . 
$$
Here $\chi_s$ denotes the characteristic function of the interval
$[-s,s]$ for $s \geq 0$. 
\end{lemma}
\vspace*{.12in}

If $K$ is the Bergman kernel for $A^2(D'_\b)$  and $K_j$ the
Bergman kernel for ${\cal H}^j$, then we may write
$$
K = \sum_{j=-\infty}^\infty K_j \, .
$$
Notice that, by rotational the invariance property of ${\cal H}^j$,
with $z = 
(z_1, z_2)$ and $w = (w_1, w_2)$, we have that 
$$
K_j(z,w) = H_j(z_1,
w_1) z_2^j \overline{w}_2^j \,.  
$$
Our job, then, is to calculate each
$H_j$, and thereby each $K_j$.  The first step of this calculation is 
already done in \cite{BAR2}, Section 2:
\begin{proposition}\label{H-j}  \sl
Let $\b>\pi/2$.  Then 
\begin{equation}\label{H-j-equation}
H_j(z_1,w_1) = \frac{1}{2\pi} \int_{-\infty}^\infty
   \frac{ e^{i(z_1 - \wbar_1)\xi} 
       \xi \left (\xi - \frac{j+1}{2} \right ) }{\sinh(\pi \xi)
        \sinh \left [(2\b - \pi) \left (\xi - \frac{j+1}{2}
         \right ) \right ] } 
        \, d\xi \, .
\end{equation}
\end{proposition}

The papers \cite{KIS} and \cite{BAR2} calculate and analyze only the
Bergman 
kernel for ${\cal H}^{-1}$ (i.e., the Hilbert subspace with index $j = -1$).  
This is attractive to do because
certain ``resonances'' cause cancellations that make the
calculations tractable when $j = -1$.  One of the main thrusts
of the present work is to perform the more difficult calculations for
all $j$. \medskip  \\

Therefore
it follows that the Bergman kernel
$K'$ for $D_\b'$ can be written as
\begin{align}
K'(z,w) 
& = \sum_{j \in \ZZ} H_j(z_1,w_1)\bigl(z_2\wbar_2\bigr)^j\notag\\
& =  \sum_{j \in \ZZ} \biggl(
\frac{1}{2\pi} \int_{-\infty}^\infty
   \frac{ e^{i(z_1 - \wbar_1)\xi} 
       \xi  \left (\xi - \frac{j+1}{2} \right ) }{\sinh(\pi \xi)
       \sinh \left ((2\b - \pi) \left (\xi - \frac{j+1}{2}
         \right ) \right ) } 
        \, d\xi \biggr) \bigl(z_2\wbar_2\bigr)^j
\, .
\label{K-prime-equation}
\end{align}

For $(z_1,z_2)$, $(w_1,w_2)\in D_\b'$ we set
$$
\tau=z_1-\wbar_1\, , \qquad \l =z_2\wbar_2\, .
$$
Notice that, if $z_1,w_1$ vary in $S_\b$, then $\tau$ varies in $S_{2\b}$. 
Moreover, we set
$$
g_j(\xi) = \frac{1}{2\pi}  \frac{ e^{i\tau\xi} 
      \xi  \bigl(\xi - \frac{j+1}{2} \bigr) }{\sinh(\pi \xi)
       \sinh \bigl((2\b - \pi)(\xi -\frac{j+1}{2} ) \bigr) } 
$$
and
\begin{equation}\label{I-j}
I_j(\tau)= 
\int_{-\infty}^{+\infty} g_j(\xi)\, d\xi \, .
\end{equation}

\begin{remark}\label{remark-1}{\rm 
Thus, in order to determine the expression of the Bergman kernel
for $D_\b'$, we  have reduced ourselves to calculate
\begin{equation}\label{our-job}
\sum_{j\in\ZZ} I_j(\tau)\lj \ ,
\end{equation}
for $\tau\in S_{2\b}$. 
The first task is to calculate $I_j(\tau)$ for each $j$.
We shall distinguish two cases according to whether
$|\Re\tau| > c_0$ or $|\Re\tau|\le c_0$, for some fixed (small)
constant $c_0$.
  
When $|\Re\tau|>c_0$ we shall
use the method of contour integrals, thus splitting
$I_j$ as sum of a residue $R_j$ and of a term $J_j$ coming from the
contour integral.  In this case, the expression of the Bergman kernel
will be given by
$$
\sum_{j \in \ZZ} R_j(\tau) \l^j +
     \sum_{j \in \ZZ} J_j(\tau) \l^j \, .
$$
The term $\sum_{j \in \ZZ} R_j(\tau) \l^j$ coming from the sum of the
residues contains the main singularity as $\Re\tau\ra\pm\infty$.

When $|\Re\tau|$ remains bounded we shall compute the expression of the
Bergman kernel by computing the sum $\sum_{j\in\ZZ} I_j(\tau)\lj$
directly. \medskip  
}
\bigskip
\end{remark}

We begin with the analysis of the case $|\Re(\tau)|>c_0$, for some
fixed constant $c_0>0$.
The following result is elementary.

\begin{lemma}\label{residues}   \sl
The function $g_j$ is holomorphic in the plane except at the points
$$
\xi=ki \ , \ \ k\in\ZZ\sm\{0\} \ ,\qquad 
\xi=\j2+ik\nu_\b \ , \ \  k\in\ZZ\sm\{0\} \, ,
$$
where $\nu_\b$ is defined in (\ref{nu-beta}).  Moreover, we have the
residue 
$$
\Res{\xi=\frac{j+1}{2} \pm i\nu_\b} g_j = 
 \frac{1}{2\pi i} \cdot
  \frac{ \nu_\b^2}{\pi} \bigl(\nu_b \mp i\frac{j+1}{2} \bigr) 
\frac{ e^{\mp \nu_\b \tau + i\tau(j+1)/2}}{
\sinh(\pi \j2 \pm i \nu_\b \pi)} \, .
$$
\end{lemma}

\proof
Since
$$
\Res{\xi=\frac{j+1}{2} \pm i\nu_\b} \frac{1}{\sinh
  \bigl((\2bp)(\xi-\j2)\bigr)} = -\frac{1}{\2bp} \, ,
$$
the desired residues are
$$
\pm \frac{1}{2\pi i}\cdot \frac{-i\nu_\b}{\2bp}\cdot
\frac{e^{i\tau\xi}\xi}{\sinh(\pi\xi)}\bigg|_{\xi=\j2\pm i\nu_\b}\, .
$$
The assertion follows at once.
\endpf

\begin{proposition}\label{R-j-plus-J-j}  \sl
Let $\b>\pi$ and fix $h$ such that
$$
\nu_\b < h< \min(1,2\nu_\b)\, .
$$
For $\Re\tau \ge 0$ and $\Re\tau<0$, 
define respectively
\begin{equation}\label{R-j-plus-J-j-equation}
R_j(\tau) = 2\pi i
\Res{\xi=\frac{j+1}{2} \pm i\nu_\b} g_j \ ,
\qquad\qquad
J_j(\tau) =
\int_{-\infty}^{+\infty} g_j(\xi\pm ih)\, d\xi\, ,
\end{equation}
and recall that $I_j$ is defined by (\ref{I-j}).
Then, for all $j\in\ZZ$,
$$
I_j(\tau) = R_j(\tau)+J_j(\tau)\, .
$$
\end{proposition}

\proof
We assume that $\nu_\b<1$, i.e. $\b>\pi$, let $\nu_\beta$ be as in
(\ref{nu-beta}), and fix $h$ such  that
$\nu_\b < h< \min(1,2\nu_\b)$ (that is, as in (\ref{(*)})).

According to whether $\Re\tau\ge0$ or $\Re\tau<0$, 
we choose
the contour of integration $\gamma^\pm_N$ to be a rectangular box,
centered on the imaginary axis, with corners at $N + i0$, $-N + i0$,
$N \pm ih$, and $-N \pm ih$.  Then
\begin{align*}
2\pi i
\Res{\xi=\frac{j+1}{2} \pm i\nu_\b} g_j 
& = \oint_{\gamma^\pm_N} g_j(\xi) \, d\xi  \\
   & =  \int_{-N}^N g_j(\xi) \, d\xi +
   \oint_{I_N}g_j d\gamma + \oint_{I_{-N}}g_j d\gamma 
- \int_{-N}^N g_j(\xi \pm ih) \, d\xi \, .
\end{align*}
Here $I_N$ and $I_{-N}$ are the left and right edges of the box.
Of course the first and
last integrals on the right are taken over the (long) top and (long) bottom
of the rectangles.  The other two (i.e., second and third) integrals 
on the right are over the (short, vertical) ends
of the rectangles.  We next verify that the 
latter two short integrals tend to zero as $N \ra +\infty$.  We now restrict
ourselves 
to the case $\Re\tau\ge0$.  The argument in the  case $\Re\tau< 0$ is
completely analogous.  

In this case, the integral over the right vertical segment equals 
\begin{equation}\label{double-star}
i \int_0^h g(N + it) \, dt = i \int_0^h 
\frac{\pi e^{i (z - \ovw)(N + it)} \bigl( 
\frac{j+1}{2} - (N + it) 
     \bigr) }{\sinh \bigl((2\b - \pi)(N + it)\bigr) 
\sinh \bigl(\pi(N + it) 
          - \frac{j+1}{2} \bigr) }  \,dt\,  .
\end{equation}
Now we write $\tau = u + iv$.  Then
$$
\Re [i\tau(N + it)] =
     \Re (i(u + iv)(N + it)) = - ut - vN \, .
$$
Hence the size of  the numerator of the integrand in
(\ref{double-star}) 
is $Ne^{-vN}$ (as a function of $N$). 
As for the denominator, notice that (for $a$ and $b$ real)
$$
\sinh(a + ib) = \sinh a \cos b + i \cosh a \sin b \, ,
$$
so that
$$
|\sinh(a + ib)|^2 
 =  \sinh^2 a + \sin^2 b \, .
$$
Thus we may estimate the size of the denominator in
(\ref{double-star})  by
$C \cdot e^{2\b N}$.

But $|v| < 2\b$, so the denominator swamps the numerator
as $N \rightarrow \infty$.  In conclusion, the integrals
over $I_N$, (and, analogously, the corresponding integral over $I_{-N}$) 
disappear for $N\rightarrow+\infty$.  The result follows.
\endpf

\section{The Sum of the $R_j$}\label{sec6}
\vspace*{.12in}

Recall that, for $\tau\in S_{2\b}$, the integral defining $I_j$ converges
absolutely. Recall also that 
we decompose $I_j$ as $I_j=R_j+J_j$. 

We now define 
 ${\cal D}$ to be the domain in $\CC^2$ given by 
\begin{equation}\label{cal-D}
{\cal D} =
\bigl\{ (\tau,\l)\in\CC^2:\, \big| \Im\tau -\log|\l|^2 \big|< \pi \, ,\ 
e^{-\bmp} < |\l| < e^{\bmp}\bigr\}\, .
\end{equation}
Then we have
\begin{proposition}\label{sum-of-Rj}  \sl 
Let $c_0>0$ be fixed.  There exists a function
$E(\tau,\l) \in C^\infty$ in a neighborhood 
of $\overline{\cal D}\cap\{(\tau,\l):\, |\Re\tau|\ge c_0\}$ in $\CC^2$
such that $D^\alpha E = \mathcal O (|\Re\tau|^\alpha)$
as
$|\Re\tau|\ra+\infty$ and 
such that, for $(\tau,\l)\in {\cal D}$ with $|\Re\tau|\geq c_0$, 
\begin{multline}
\sum_{j \in \ZZ} R_j(\tau) \l^j  \\
= e^{-\nb \tau 
  \sgn(\Re \tau)} \cdot \left ( 
\frac{\varphi_1 (\tau)}{(1 - e^{[i\tau - \pi]/2}\l)^2}
   + \frac{\varphi_2 (\tau,\l)}{(\l - e^{-[\pi + i\tau]/2})^2} +
    E(\tau,\l) \right ) \, ,  \label{sum-of-Rj-equation}
\end{multline}
where
\begin{align*}
\varphi_1 (\tau) 
& = \frac{2\nu_\b^2}{\pi}  e^{-i\pi\nb} e^{[i\tau - \pi]/2}
\biggl(  \frac{i}{2}+\nb\bigl(1-e^{[i\tau-\pi]/2}\bigr) \biggr) 
\, , \\
\varphi_2 (\tau,\l)
& = \frac{2\nu_\b^2}{\pi}  
e^{i\pi\nb} \l e^{-[i\tau + \p]/2}
\biggl( \nb e^{-[i\tau+\pi]/2}-\bigl(\nb -\frac{i}{2}\bigr)\l\biggr)
\, .
\end{align*}
The convergence of the series
is uniform on compact subsets of ${\cal D}$.
\end{proposition}
\proof
We obtain the formula (\ref{sum-of-Rj-equation}) by approximating the  
$\sinh$ function by an exponential, summing the resulting geometric
series, and then showing that the performed approximation gives rise
to error terms that are more regular than the explicit term
obtained by replacing the $\sinh$ functions by exponentials.

Notice that, for $a,b\in\RR$, $a \neq 0$, $|b|<\pi$, 
\begin{equation}\label{sinh}
\frac{e^{|a|}}{\sinh (a + ib)}
= 2\, \sgn (a) e^{-i\sgn(a)b} \biggl(1 +
\frac{e^{-2\sgn(a)(a+ib)}}{1- e^{-2\sgn(a)(a+ib)}}\biggr) \, .
\end{equation}

For simplicity of notation we momentarily
restrict ourselves to the case
$\Re\tau\ge0$.  The argument in  the case $\Re\tau<0$ is analogous.
Using  Lemma \ref{residues} and Proposition \ref{R-j-plus-J-j}, 
we write
\begin{eqnarray*}
R_j(\tau) 
& =  &\frac{\nu_\b^2}{\pi} 
            \frac{e^{-\nb\tau}\left (\nb + i \j2 \right )
             e^{i \j2 \tau}}{\sinh (\pi \j2 + i\nb \pi)} \\
& = & \frac{\nu_\b^2}{\pi}  e^{-\nb\tau} \left [ \nb 
                     \frac{e^{i\j2\tau}}{\sinh
       (\pi \j2 + i\nb\pi)} + i \j2 \frac{e^{i\j2\tau}}{\sinh
       (\pi \j2 + i \nb \pi)} \right ] \, .
\end{eqnarray*}

We now use (\ref{sinh}). 
The first sum that we must consider is
\begin{align}
& \sum_{j \in \ZZ} \frac{e^{i\j2\tau}}{\sinh(\pi \j2 + i \nb \pi)}\lj \notag\\
& \qquad = 2 
\biggl(
\sum_{j \in \ZZ,\, j\neq-1} \sigma(j) e^{-i\sigma(j)\pi\nu_\b} 
e^{i\j2 \tau - \pi |\j2|} \lj \notag \\ 
& \qquad\qquad   +  \sum_{j \in \ZZ} \sigma(j) e^{-i\sigma(j)\pi\nu_\b} e^{i\j2 \tau} 
\frac{e^{- 3\pi |\j2| -2i\pi\nb\sigma(j)}}{1-
  e^{-2\pi|\j2|-2i\pi\nb\sigma(j)}} \lj 
+\frac{1}{\sinh
  (i\pi\nu_b)}  
\biggr) 
\notag  \\ 
& \qquad \equiv 2
(F_1 +E_1) \, ; \label{first-sum}
\end{align}
here $\sigma(j)=\sgn(j+1)$.
On the other hand, the second sum is
\begin{align}
& \sum_{j \in \ZZ} \j2  
\frac{e^{i\j2\tau}}{\sinh(\pi \j2 + i \pi\nb)}\lj \notag\\
& \qquad 
= 2 
\biggl(
\sum_{j \in \ZZ} \j2 \sigma(j)e^{-i\sigma(j)\pi\nu_\b} 
e^{i\j2 \tau - \pi |\j2|} \lj \notag\\
& \qquad\qquad\qquad
+   \sum_{j \in \ZZ} \j2 \sigma(j)e^{-i\sigma(j)\pi\nu_\b}  e^{i\j2 \tau} 
\frac{e^{- 3\pi |\j2| -2i\pi\nb\sigma(j)}}{1- e^{-\pi|\j2|-2i\pi\nb\sigma(j)}} \lj \biggr)
\notag  \\
& \qquad \equiv 2  (F_2 +E_2) \, . \label{second-sum}
\end{align}
Notice that then
\begin{equation}\label{dagger}
\sum_{j \in \ZZ} R_j(\tau) \lj = 
\frac{2\nu_\b^2}{\pi}  
e^{-\nb\tau} \left [ \nb (F_1 + E_1)
                      + i (F_2+E_2) \right ] \, .
\end{equation}

Now
\begin{eqnarray*}
F_1 & = & \sum_{j \in \ZZ,\, \neq-1} 
\sigma(j)e^{-i\sigma(j)\pi\nu_\b} e^{i(\j2)\tau - \pi|\j2|} \l^j \\
   & = & e^{-i\pi\nb}
\sum_{j \ge 0} e^{(i\tau - \pi)(\j2)} \lj 
      e^{i\pi\nb} \sum_{j < -1} e^{(i\tau + \pi)\j2} \lj \\
   & = & e^{-i\pi\nb+[i\tau - \pi]/2} \sum_{j \ge 0}
             \left ( e^{[i\tau - \pi]/2} \l \right )^j - 
         e^{-i\pi\nb+[i\tau + \pi]/2} \sum_{j \geq 1}
             \left ( e^{-[i\tau + \pi]/2} \l^{-1} \right )^j  + \frac{e^{i\pi\nb}}{\l}\\
   & = & 
\frac{e^{-i\pi\nb+[i\tau - \pi]/2} }{1 - e^{[i\tau-\pi]/2}\l}
  - \frac{e^{i\pi\nb-[i\tau + \pi]/2}}{\l - e^{-[i\tau + \pi]/2}} + \frac{e^{i\pi\nb}}{\l}
\, .
\end{eqnarray*}
Here, of course, the first sum converges for 
$|e^{[i\tau - \pi]/2}\l| < 1$ and the second
sum converges for $|e^{-[i\tau + \pi]/2}\l^{-1}| < 1$.
So the full series, summed over $j \in \ZZ$, converges
on the annulus $e^{[\Im \tau - \pi]/2} < |\l|
< e^{[\Im \tau + \pi]/2}$.  

For $F_2$, notice that 
the sum is precisely the same as the one we just computed
except for a factor of $\j2$ in front. Thus 
we can formally differentiate in $\tau$ and obtain
\begin{align*}
F_2 
& =  \sum_{j \in \ZZ} \j2  e^{i(\j2)\tau - \pi|\j2|} \l^j 
= \frac{1}{i}\frac{d}{d\tau} F_1 \\
& = \frac{1}{2}  \frac{e^{-i\pi\nb+[i\tau - \pi]/2}}{(1 -
                  e^{[i\tau - \p]/2}\l)^2} +
      \frac{1}{2}  \frac{\l e^{i\pi\nb-[i\tau + \p]/2}}{(\l
        - e^{-[i\tau + \p]/2})^2} \, .
\end{align*}
Therefore
\begin{align}
& \nb F_1+iF_2  \notag \\ 
& =  \nb \biggl( 
\frac{e^{-i\pi\nb+[i\tau - \pi]/2} }{1 - e^{[i\tau-\pi]/2}\l}
  -  \frac{e^{i\pi\nb-[i\tau + \pi]/2}}{\l - e^{-[i\tau + \pi]/2}} \biggr)
+  \frac{i}{2}\biggl(  
\frac{e^{-i\pi\nb+[i\tau - \pi]/2}}{(1 - e^{[i\tau - \p]/2}\l)^2} 
      \frac{\l e^{i\pi\nb-[i\tau + \p]/2}}{(\l- e^{-[i\tau + \p]/2})^2} 
\biggr)  \notag \\
& = \frac{e^{-i\pi\nb+[i\tau - \pi]/2}}{(1 - e^{[i\tau - \p]/2}\l)^2}
\biggl(  \frac{i}{2}+\nb\bigl(1-e^{[i\tau-\pi]/2}\bigr) \biggr) 
\notag \\
& \qquad
+  \frac{e^{i\pi\nb-[i\tau + \p]/2}}{(\l- e^{-[i\tau + \p]/2})^2} 
\biggl( \nb e^{-[i\tau+\pi]/2} -\bigl(\nb -\frac{i}{2}\bigr)\l\biggr)
\, .   \label{m1+m2}
\end{align}
Hence, by (\ref{dagger}) and (\ref{m1+m2}),
and recalling that $\Re\tau\ge c_0>0$, we have
\begin{multline}
\sum_{j \in \ZZ} R_j(\tau) \lj =  \\
e^{-\nb\tau} 
\left [ 
\frac{\varphi_1 (\tau,\l)}{(1 - e^{[i\tau - \p]/2}\l)^2}
+  \frac{\varphi_2 (\tau,\l)}{(\l- e^{-[i\tau + \p]/2})^2} \right]
+ 
\frac{\nu_\b^2}{2\pi}  e^{-\nb\tau+i\pi\nb}
\left [ \nb E_1 + i E_2
 \right ] \, , \label{ddagger}
\end{multline}
where
\begin{align*}
\varphi_1 (\tau,\l)
& = \frac{2\nu_\b^2}{\pi}  e^{-i\pi\nb} e^{[i\tau - \pi]/2}
\biggl(  \frac{i}{2}+\nb\bigl(1-e^{[i\tau-\pi]/2}\bigr) \biggr) 
\, ,
\\
\varphi_2 (\tau,\l)
& = \frac{2\nu_\b^2}{\pi}  
e^{i\pi\nb} \l e^{-[i\tau + \p]/2}
\biggl( \nb e^{-[i\tau+\pi]/2} -\bigl(\nb -\frac{i}{2}\bigr)\l\biggr) \,
,
\end{align*}
as we claimed.
Notice in particular $\varphi_1 ,\varphi_2 $ are entire functions, 
and are bounded as
$|\Re\tau|\rightarrow+\infty$, together with all their derivatives.
\medskip   \\

We now study the error term
$$
E \equiv \frac{2\nu_\b^2}{\pi}  e^{-\nb\tau-i\pi\nb}
\left [ \nb E_1 + i E_2
 \right ] \, . 
$$

Notice that there exists a positive constant $c_0 > 0$
such that, for all $j$, 
$$
|1- e^{-2\pi|\j2|-2i\pi\nb}|\ge c_0 \, .
$$
Hence the series
$$
\sum_{j \in \ZZ} e^{i\j2 \tau} 
\frac{e^{- 3\pi |\j2| -2i\pi\nb\sigma(j)}}{1- e^{-2\pi|\j2|-2i\pi\nb\sigma(j)}} \lj
\qquad\text{and}\quad
\sum_{j \in \ZZ} \j2  e^{i\j2 \tau} 
\frac{e^{- 3\pi |\j2| -2i\pi\nb\sigma(j)}}{1- e^{-2\pi|\j2|-2i\pi\nb\sigma(j)}} \lj
$$
converge  when
$$
\left | e^{-[i\tau + 3\p]/2} \right | < |\l| <
          \left | e^{-[i\tau - 3\p]/2} \right | \, , 
$$
which is a strictly
larger annulus than 
$e^{[\Im \tau - \pi]/2} < |\l|
< e^{[\Im \tau + \pi]/2}$.  Thus
the sums of the two 
series are smooth and bounded,
with all derivatives smooth and bounded, on a neighborhood
of the closure $\overline{\cal D}$ of ${\cal D}$ (where this domain
is defined at the beginning of Section \ref{sec6}).

For the case $\Re\tau<-c_0$, notice that
\begin{align*}
R_j(\tau) 
& = \frac{\nu_\b^2}{\pi} 
            \frac{e^{\nb\tau}\left (\nb - i \j2 \right )
             e^{i \j2 \tau}}{\sinh (\pi \j2 - i\nb \pi)} \\
& = \frac{\nu_\b^2}{\pi}  e^{\nb\tau} \left [ \nb 
                     \frac{e^{i\j2\tau}}{\sinh
       (\pi \j2 - i\nb\pi)} - i \j2 \frac{e^{i\j2\tau}}{\sinh
       (\pi \j2 - i \nb \pi)} \right ] \, .
\end{align*}
Then we consider the sums
\begin{align*}
& 2 
\biggl(
\sum_{j \in \ZZ} \sigma(j)e^{-i\sigma(j)\pi\nu_\b} e^{i\j2 \tau - \pi |\j2|} \lj 
   +  \sum_{j \in \ZZ} \sigma(j)e^{-i\sigma(j)\pi\nu_\b} e^{i\j2 \tau} 
\frac{e^{- 3\pi |\j2| -2i\sigma(j)\pi\nb}}{1-
  e^{-2\pi|\j2|-2i\sigma(j)\pi\nb}} \lj \biggr) \, ,\\ 
& 2
\biggl(
\sum_{j \in \ZZ} \j2 \sigma(j)e^{-i\sigma(j)\pi\nu_\b} 
e^{i\j2 \tau - \pi |\j2|} \lj
+   \sum_{j \in \ZZ} \j2  \sigma(j)e^{-i\sigma(j)\pi\nu_\b} e^{i\j2 \tau} 
\frac{e^{- 3\pi |\j2| -2i\sigma(j)\pi\nb}}{1-
  e^{-2\pi|\j2|-2i\sigma(j)\pi\nb}} \lj \biggr) \ .
\end{align*} 

Therefore, when $\Re\tau<c_0$,
$$
\sum_{j\in\ZZ} R_j(\tau)\lj 
= \frac{2\nu_\b^2}{\pi} e^{\nb\tau} \bigl[
\nb (F_1 +E_1') -i(F_2+E_2')\bigr]\, ,
$$
where $F_1,F_2$ are as in the previous case, while $E_1',E_2'$ are
error terms that are dealt with as in the case $\Re\tau\ge c_0$.
\bigskip
\endpf

\section{The Sum of the $J_j$}\label{decomposition-of-Jj} 
\vspace*{.12in}

We wish to study the sum
$$
\sum_{j \in \ZZ} J_j(\tau) \lj \, ,
$$
where
$$
J_j(\tau) = \frac{1}{2\pi} \int_{-\infty}^{+\infty} 
    \frac{e^{i\tau(\xi + ih)}(\xi + ih)(\xi + ih - \j2) }{\sinh(
        \pi(\xi + ih))\sinh((\2bp)(\xi + ih - \j2))}\, d\xi \, .
$$

\begin{remark}\label{h}{\rm 
We make an important observation.  We arrived at the definition of
$J_j(\tau)$ by computing a contour integral along the rectangular path
$\gamma_N^\pm$, whose upper (lower) side is the segment $[-N\pm
ih,N\pm ih]$, for $\Re\tau>0$  ($\Re\tau<0$) respectively.  However, the
definition of $J_j(\tau)$ makes sense for all real values of the
parameter $h$, in particular for $h=0$.  Notice that, in this case,
$J_j(\tau)=I_j(\tau)$.  We shall make use of this fact when we compute  
the expression of the Bergman kernel for $|\Re(z_1-\wbar_1)|\le
c_0$. 
}
\end{remark}

We use (\ref{sinh}) to write
\begin{align*}
& \frac{1}{\sinh(\pi(\xi + ih))\sinh((\2bp)(\xi +ih-\j2))} \\
& \qquad = 
\frac{2\sgn(\xi)e^{-i\sgn(\xi)\pi h}}{e^{\pi|\xi|}}
\biggl( 1+ 
\frac{e^{-2\sgn(\xi)\pi(\xi+ih)}}{1-e^{-2\sgn(\xi)\pi(\xi+ih)}}
\biggr) \\ 
& \qquad \qquad \times 
\frac{2\sgn(\xi-\j2)e^{-i\sgn(\xi-\j2)(\2bp)h}}{e^{(\2bp)|\xi-\j2|}}
\biggl( 1+ 
\frac{e^{-2\sgn(\xi-\j2)(\2bp)(\xi-\j2+ih)}}{1
  -e^{-2\sgn(\xi-\j2)(\2bp)(\xi-\j2+ih)}}  
\biggr) \\ 
& \qquad =
\frac{4\sgn(\xi)\sgn(\xi-\j2)e^{
    -ih\bigl(\sgn(\xi)\pi +\sgn(\xi-\j2)(2\b-\pi)\bigr) }}{
              e^{\pi|\xi|}e^{(\2bp)|\xi-\j2|}}  \\
& \qquad \qquad \times \biggl( 1 + 
\frac{e^{-2\sgn(\xi)\pi(\xi+ih)}}{1-e^{-2\sgn(\xi)(\xi+ih)}} 
 + \frac{e^{-2\sgn(\xi-\j2)(\2bp)(\xi-\j2+ih)}}{1
  -e^{-2\sgn(\xi-\j2)(\2bp)(\xi-\j2+ih)}} \\
& \qquad\qquad \qquad +
\frac{e^{-2\sgn(\xi)\pi(\xi+ih)}}{1-e^{-2\sgn(\xi)\pi(\xi+ih)}} \cdot
\frac{e^{-2\sgn(\xi-\j2)(\2bp)(\xi-\j2+ih)}}{1
  -e^{-2\sgn(\xi-\j2)(\2bp)(\xi-\j2+ih)}}  
\biggr) \, .
\end{align*}
We set
\begin{equation}\label{sigma}
\sigma(\xi)=
\sgn(\xi)\sgn(\xi-\j2)
e^{-ih\bigl(\sgn(\xi)\pi +\sgn(\xi-\j2)(2\b-\pi)\bigr)}
\end{equation}
and decompose $J_j$ accordingly as 
\begin{equation}\label{decompose-J-j}
J_j(\tau) =  M_j(\tau) + E_j^{(1)}(\tau) 
+E_j^{(2)}(\tau) +E_j^{(3)}(\tau) \, ,
\end{equation}
where
\begin{align}
M_j(\tau)
& = 
\frac{2}{\pi}\int_{-\infty}^{+\infty} 
    \sigma(\xi) e^{i\tau(\xi + ih)}
\frac{
(\xi +ih)(\xi +ih -\j2)}{
e^{\pi |\xi|}e^{(\2bp)|\xi-\j2|}} \, d\xi  \notag \\
&  =  \frac{2}{\pi}e^{-h\tau}
\int_{-\infty}^{+\infty} \sigma(\xi)
e^{i\tau\xi}
\frac{(\xi +ih)(\xi +ih -\j2) 
}{e^{\pi|\xi|}e^{(\2bp)|\xi-\j2|}} \, d\xi  \, , \label{M-j}
\end{align}
and, analogously,
\begin{align} 
E_j^{(1)}(\tau) 
& = \frac{2}{\pi}e^{-h\tau} 
\int_{-\infty}^{+\infty} 
    \sigma(\xi) e^{i\tau\xi}
\frac{(\xi +ih)(\xi +ih -\j2)}{e^{\pi|\xi|}e^{(\2bp)|\xi-\j2|}} 
\cdot
\frac{e^{-2\sgn(\xi)\pi(\xi+ih)}}{1-e^{-2\sgn(\xi)(\xi+ih)}} 
\, d\xi\, , \notag \\
E_j^{(2)}(\tau) 
& = \frac{2}{\pi}e^{-h\tau} 
\int_{-\infty}^{+\infty}
    \sigma(\xi) e^{i\tau\xi}
\frac{(\xi +ih)(\xi +ih -\j2)}{e^{\pi|\xi|}e^{(\2bp)|\xi-\j2|}} \, ,
\notag\\ 
& \qquad\qquad\qquad\qquad \times
\frac{e^{-2\sgn(\xi-\j2)(\2bp)(\xi-\j2+ih)}}{1
  -e^{-2\sgn(\xi-\j2)(\2bp)(\xi-\j2+ih)}} 
\, d\xi\,  \label{E3-j} \\
E_j^{(3)}(\tau) 
& = \frac{2}{\pi}e^{-h\tau} 
\int_{-\infty}^{+\infty}
    \sigma(\xi) e^{i\tau\xi}
\frac{(\xi +ih)(\xi +ih -\j2)}{e^{\pi|\xi|}e^{(\2bp)|\xi-\j2|}} \notag\\
& \qquad\qquad\qquad\qquad \times
\frac{e^{-2\sgn(\xi)\pi(\xi+ih)}}{1-e^{-2\sgn(\xi)\pi(\xi+ih)}} 
\cdot
\frac{e^{-2\sgn(\xi-\j2)(\2bp)(\xi-\j2+ih)}}{1
  -e^{-2\sgn(\xi-\j2)(\2bp)(\xi-\j2+ih)}}  
\, d\xi\, . \notag
\end{align}

Thus we have reduced the situation to computing the sums
$$
\sum_{j\in\ZZ} M_j(\tau)\lj +
\sum_{j\in\ZZ} E_j^{(1)}(\tau)\lj 
+ \sum_{j\in\ZZ} E_j^{(2)}(\tau)\lj
+ \sum_{j\in\ZZ} E_j^{(3)}(\tau)\lj  
\, .
$$
We shall prove the following results.
\begin{proposition}\label{sum-of-Mj-proposition}  \sl
There exist entire functions $\psi_j^\pm$,
$j=1,\dots,4$,
which are uniformly (together with all their derivatives) of
size $\mathcal{O}(|\Re\tau|)$ 
as $|\Re\tau|\ra+\infty$ in any set  
$\{|\Im\tau|\le C,\ |\l|\le C'\}$,  
such that  
\begin{align*}
\sum_{j\in\ZZ} M_j(\tau)\lj 
& = \frac{1}{\pi} e^{-\sgn(\Re\tau)h\tau}  \biggl\{ 
  \frac{e^{i2\b h}}{(i\tau+2\b)^2} 
\biggl[
\frac{e^{\bmp}}{2(e^{\bmp} - \l)^2}  
   \bigl(1+\psi_1^+(\tau,\l) \bigr) \\    
&\qquad\quad
+ \frac{e^{-[i\tau +\pi]/2}}{2(\l - e^{-[i\tau +\pi]/2})^2}
   \bigl( 1 + \psi_2^-(\tau,\l) \bigr)  
+ \frac{e^{-i2\b h-[i\tau+\pi]/2}\psi_4^+(\tau)}{
2(e^{\bmp} - \l)(e^{-[i\tau + \pi]/2}-\l)} \biggr] 
\\
& \qquad +  
\frac{e^{-i2\b h}}{(i\tau-2\b)^2} 
\biggl[
\frac{e^{[\pi- i\tau]/2}}{2(e^{[\pi - i\tau]/2} - \l)^2}
        \bigl(1 + \psi^+_2(\tau,\l) \bigr) \\
&\qquad\quad
+ \frac{  e^{-\bmp}}{2(\l - e^{-\bmp})^2}  
        \bigl(1 +\psi_1^-(\tau,\l) \bigr)  
+ \frac{e^{-i2\b h+[\pi-i\tau]/2}\psi_4^-(\tau)}
{2(e^{-\bmp} - \l)(e^{[\pi-i\tau]/2}-\l)}  \biggr] \\
&\qquad
+
\frac{e^{ih(2\b-2\pi)+ \bmp}}{2(e^{[\pi-i\tau]/2}-\l)^2(e^{\bmp}-\l)^2}
\psi_3^+(\tau,\l) \\
&\qquad
+
\frac{e^{-ih(2\b-2\pi)- \bmp}}{2(e^{-[i\tau+\pi]/2}-\l)^2(e^{-\bmp}-\l)^2} 
\psi_3^-(\tau,\l) \biggr\}
\ .
\end{align*}
Here $h$ is as in (\ref{(*)}).

The convergence is uniform on compact subsets of
$D_\b'$.  Moreover, the functions $\psi_j^\pm$
can be computed explicitly (see  
(\ref{psi1+psi2+}), (\ref{psi1-psi2-}), (\ref{new-psi3+}), 
(\ref{new-psi3-}),
and 
(\ref{psi4+-}) respectively).
\end{proposition}
   
\begin{proposition}\label{sum-of-error-terms-proposition}  \sl
There exist functions $\Psi_1^{(k)},\,  
\Psi_2^{(k)}$ and $ \Psi_3^{(k)}$, holomorphic in a
neighborhood of $\overline{{\cal D}}$, bounded together with all their
derivatives 
as $|\Re\tau|\ra+\infty$, such that, for $k=1, 2, 3$,
$$
\sum_{j\in\ZZ} E_j^{(k)}(\tau)\lj = 
e^{-\sgn(\Re\tau)h\tau}\biggl\{ \frac{\Psi_1^{(k)}(\tau,\l)}{(\l - e^{-\bmp})^2} 
+\frac{\Psi_2^{(k)}(\tau,\l)}{(e^{\bmp}-\l)^2} +\Psi_3^{(k)}(\tau,\l)\biggr\} \ .
$$
Here $h$ is as in (\ref{(*)}).
\end{proposition}
(Recall that the $E_j^{(k)}$ are defined in (\ref{E3-j}).)

We shall compute these sums, and prove these two propositions, in Sections 
\ref{new-section}   and \ref{new-section-bis}.  
\bigskip	    

\section{The Bergman Kernels for $D_\b$  and $D'_\b$ --   
Proof of Theorems \ref{THM1} and \ref{THM2} }  
\label{K-K'}
\vspace*{.12in}

In this section, assuming the validity of
Propositions \ref{sum-of-Mj-proposition}  
and \ref{sum-of-error-terms-proposition},
we complete the proof  of the asymptotic
expansions 
for the Bergman kernels for the worm domains $D_\b$ and $D'_\b$.  

\proof[Proof of Theorem \ref{THM1}]
Before we move ahead, let us review what we have accomplished, and how
all the parts fit together.

{From} Remark \ref{remark-1} we know that the Bergman kernel
$K_{D'_\b}$ for 
$D'_\b$ can be written as
$$
K_{D'_\b}(z_1,z_2,w_1,w_2)
= \sum_{j\in\ZZ} I_j(\tau)\lj\, ,
$$
where we set $z_1-\wbar_1=\tau$ and $z_2\wbar_2=\l$.
 
When $|\Re\tau|>c_0$ we use the decomposition
$$
\sum_j I_j(\tau) \l^j = \sum_j R_j(\tau) \l^j +  
	   \sum_j M_j(\tau)\l^j + \sum_{k=1}^3 \sum_j E_j^{(k)}
           (\tau) \l^j \, . 
$$

We put together
(\ref{our-job}), Proposition \ref{R-j-plus-J-j}, 
Proposition \ref{sum-of-Rj}, 
(\ref{decompose-J-j}), Proposition \ref{sum-of-Mj-proposition},
and Proposition \ref{sum-of-error-terms-proposition}. 
For $|\Re (z_1- \ovw_1)|>c_0$, we have that
\begin{align*}
&K_{D'_\b}(z_1,z_2,w_1,w_2) \notag \\ 
& =  e^{-\nb (z_1- \ovw_1)   \sgn(\Re (z_1- \ovw_1))} \biggl\{ 
 \frac{\varphi_1 (z_1- \ovw_1)}{(1 - e^{[i(z_1- \ovw_1) - \pi]/2}z_2 \ovw_2)^2} 
   + \frac{\varphi_2 (z_1- \ovw_1,z_2 \ovw_2)}{
(z_2 \ovw_2 - e^{-[\pi + i(z_1- \ovw_1)]/2})^2} \biggr\} \notag \\
& \quad +
\frac{1}{\pi} e^{-\sgn(\Re(z_1- \ovw_1))h(z_1- \ovw_1)}  \biggl\{ 
  \frac{e^{i2\b h}}{(i(z_1- \ovw_1)+2\b)^2} 
\biggl[
\frac{e^{\bmp}}{2(e^{\bmp} - z_2 \ovw_2)^2}  
   \bigl(1+\psi_1^+(z_1- \ovw_1,z_2 \ovw_2) \bigr) \notag \\    
&\qquad\quad
+ \frac{e^{-[i(z_1- \ovw_1) +\pi]/2}}{2(z_2 \ovw_2 - e^{-[i(z_1- \ovw_1) +\pi]/2})^2}
   \bigl( 1 + \psi_2^-(z_1- \ovw_1,z_2 \ovw_2) \bigr)\notag \\  
& \quad
+ \frac{e^{-i2\b h-[i(z_1- \ovw_1)+\pi]/2}\psi_4^+(z_1- \ovw_1)}{
2(e^{\bmp} - z_2 \ovw_2)(e^{-[i(z_1-\ovw_1)+\pi]/2}-z_2 \ovw_2)} \biggr] 
\notag \\
\end{align*}
\begin{align}
& \qquad +  
\frac{e^{-i2\b h}}{(i(z_1- \ovw_1)-2\b)^2} 
\biggl[
\frac{e^{[\pi- i(z_1- \ovw_1)]/2}}{2(e^{[\pi - i(z_1- \ovw_1)]/2} - z_2 
\ovw_2)^2}
        \bigl(1 + \psi^+_2(z_1- \ovw_1,z_2 \ovw_2) \bigr) \notag \\
&\qquad\quad
+ \frac{  e^{-\bmp}}{2(z_2 \ovw_2 - e^{-\bmp})^2}  
        \bigl(1 +\psi_1^-(z_1- \ovw_1,z_2 \ovw_2) \bigr)  \notag \\
&\quad 
+ \frac{e^{-i2\b h+[\pi-i(z_1- \ovw_1)]/2}\psi_4^-(z_1- \ovw_1)}
{2(e^{-\bmp}-z_2 \ovw_2)(e^{[\pi-i(z_1- \ovw_1)]/2}-z_2 \ovw_2)}  
\biggr] \notag \\
&\quad
+
\frac{e^{ih(2\b-2\pi)+ \bmp}}{2(e^{[\pi-i(z_1- \ovw_1)]/2}-z_2 \ovw_2)^2
(e^{\bmp}-z_2 \ovw_2)^2}
\psi_3^+(z_1- \ovw_1,z_2 \ovw_2) \notag \\
&\quad
+
\frac{e^{-ih(2\b-2\pi)- \bmp}}{2(e^{-[i(z_1- \ovw_1)+\pi]/2}-z_2 \ovw_2)^2
(e^{-\bmp}-z_2 \ovw_2)^2} 
\psi_3^-(z_1- \ovw_1,z_2 \ovw_2)  \notag \\
&\quad
+ \frac{\Psi_1(z_1- \ovw_1,z_2 \ovw_2)}{
(z_2 \ovw_2 - e^{-\bmp})^2} 
+\frac{\Psi_2(z_1- \ovw_1,z_2 \ovw_2)}{
(e^{\bmp}-z_2 \ovw_2)^2} +\Psi_3(z_1- \ovw_1,z_2 \ovw_2)
\biggr\} \ , \label{exp-K'-1}
\end{align}
where $\Psi_j \equiv \sum_{k=1}^3 \Psi_j^{(k)}$.
We recall that $h$ is a fixed parameter, $\nu_\b<h<\min(1,2\nu_\b)$, 
and that the functions
$\psi_j^\pm$s, $\varphi$s and $\Psi$s depend (smoothly) on
$h$.  Moreover, we recall that 
these functions all sasisfy the condition 
${\cal O}(|\Re z_1-\Re w_1|)$, together with all their
derivatives, for $z,w \in \overline{D'_\b}$.  
\medskip       \\

For $|\Re(z_1-w_1)|\le c_0$ we use the result from Section
\ref{decomposition-of-Jj}, with $h=0$, for which see Remark \ref{h}. 
By Propositions \ref{sum-of-Mj-proposition} and
\ref{sum-of-error-terms-proposition}, with $h=0$, we have that
when $|\Re(z_1-w_1)|\le c_0$,
\begin{align}
&K_{D'_\b}(z_1,z_2,w_1,w_2) \notag \\  
& =
  \frac{e^{\bmp}\bigl(1+\psi_1^+(z_1-\ovw_1,z_2\ovw_2)\bigr) 
}{2(i(z_1- \ovw_1)+2\b)^2 (e^{\bmp} - z_2 \ovw_2)^2}  
\notag \\    
&\qquad\quad
+ \frac{e^{-[i(z_1- \ovw_1) +\pi]/2}\bigl( 1 + \psi_2^-(z_1-
  \ovw_1,z_2 \ovw_2)   \bigr)}{2(i(z_1- \ovw_1)+2\b)^2
(z_2 \ovw_2 -e^{-[i(z_1-\ovw_1)+\pi]/2})}
\notag \\  
& \quad
+ \frac{e^{-[i(z_1- \ovw_1)+\pi]/2}\psi_4^+(z_1- \ovw_1)}{
4(i(z_1- \ovw_1)+2\b)^2(e^{\bmp} - z_2 \ovw_2)
(e^{-[i(z_1-\ovw_1)+\pi]/2}
-z_2 \ovw_2)
} 
\notag \\
& \qquad +  
\frac{e^{[\pi- i(z_1- \ovw_1)]/2} \bigl(1 + \psi^+_2(z_1-
  \ovw_1,z_2 \ovw_2) \bigr) }{2(i(z_1- \ovw_1)-2\b)^2 
(e^{[\pi - i(z_1- \ovw_1)]/2} - z_2 \ovw_2)^2}
\notag \\
&\qquad\quad
+ \frac{  e^{-\bmp}\bigl(1 +\psi_1^-(z_1- \ovw_1,z_2 \ovw_2)
  \bigr)}{2(i(z_1- \ovw_1)-2\b)^2  
(z_2 \ovw_2 - e^{-\bmp})^2}  
\notag \\
\end{align}

\begin{align}
&\quad 
+ \frac{e^{[\pi-i(z_1- \ovw_1)]/2}\psi_4^-(z_1- \ovw_1)}
{4(i(z_1- \ovw_1)-2\b)^2  
(e^{-\bmp}-z_2 \ovw_2)(e^{[\pi-i(z_1- \ovw_1)]/2}-z_2 \ovw_2)}  
 \notag \\
&\quad
+
\frac{e^{\bmp}\psi_3^+(z_1- \ovw_1,z_2 \ovw_2) 
}{2(e^{[\pi-i(z_1- \ovw_1)]/2}-z_2 \ovw_2)^2
(e^{\bmp}-z_2 \ovw_2)^2}
\notag \\
&\quad
+
\frac{e^{-\bmp}\psi_3^-(z_1- \ovw_1,z_2 \ovw_2)  
}{2(e^{-[i(z_1- \ovw_1)+\pi]/2}-z_2 \ovw_2)^2
(e^{-\bmp}-z_2 \ovw_2)^2} 
+E(z_1- \ovw_1,z_2 \ovw_2)\ .
\label{exp-K'-2}
\end{align}
Here $E$ denotes a function which is holomorphic 
$(z_1,z_2)$ and anti-holomorphic in $(w_1,w_2)$, for 
$(z_1,z_2)$ and $(w_1,w_2)$ varying 
in a neighborhood
of $D'_\b$.

We should point out that the functions $\psi^{\pm}_j$ in the formula above
differ from
the ones appearing in (\ref{exp-K'-1}), since they
are obtained
from Propositions \ref{sum-of-Mj-proposition} and
\ref{sum-of-error-terms-proposition}
for different values of the parameter $h$.
\medskip    \\

Now, from (\ref{exp-K'-1}) and (\ref{exp-K'-2}),  
the proof of Theorem \ref{THM1} follows at once.
\endpf

\bigskip
\proof[Proof of Theorem \ref{THM2}]
We now derive the explicit expression for the asymptotic
expansion  
of the Bergman kernel for the worm domain $D_\b$.  We use the formula
from Theorem \ref{THM1}, together with the transformation formula
under biholomorphic mappings (see \cite{KRA1}).  
Specifically, recall that
the domains $D'_\b$ and $D_\b$ are biholomorphic via the mapping
\begin{eqnarray*}
\Phi: D'_\b & \rightarrow & D_\b \\
(z_1, z_2) & \mapsto & (e^{z_1}, z_2) \equiv (\z_1, \z_2) \, .
\end{eqnarray*}

Hence
$$
\begin{matrix}
\Phi^{-1}:& D_\beta & \longrightarrow & D'_\b  \\
         & (\om_1, \om_2) & \mapsto & 
          (\log\om_1,\om_2) = (w_1, w_2) \, . \end{matrix}
$$
We see that
\begin{align*}
K_{D_\beta} \bigl((\z_1, \z_2), (\om_1, \om_2)\bigr) 
& = K_{D'_\b}\bigl(\Phi^{-1}(\z), \Phi^{-1}(\om)\bigr) 
\cdot \det \Jac \Phi^{-1} (\z)
	          \overline{\det \Jac_{\Phi^{-1}} (\om)} \, .
\end{align*}

Notice that, on $D_\b$, the function $\log\z_1$ is well defined once
$\z_2$ is fixed, since $\z_1$ lies in the half-plane
$\Re(\z_1 e^{i\log|\z_2|^2})>0$.  Hence, for $\z=(\z_1,\z_2)\in D_\b$,
the number $\log\z_1$ is the unique value obtained by analytic continuation
starting from the principal branch of the $\log$ when
$\z_2$ lies on the positive real line.  Also, 
observe that, via the transformation $\Phi$, 
$$
(z_1 - \wbar_1)  \mapsto 
\log (\z_1/\overline\om_1)
\quad\text{and}\quad
 e^{[i(z_1 - \wbar_1)+ \pi]/2} \mapsto 
\biggl(\frac{\z_1}{\overline{\om}_1}\biggr)^{i/2} e^{\pi/2}
\ , 
$$
where, by $\log (\z_1/\overline\om_1)$ we mean 
$\log\z_1 - \log\overline\om_1$.

Thus, since $\Jac \Phi^{-1}(\om)=1/\om_1$, 
writing $t=|\z_1|-|\om_1|$,
we have that
\begin{align}
 & K_{D_\b}(\z,\om) =    
\frac{1}{\z_1 \overline{\om}_1} 
K_{D'_\b}\bigl( (\log \z_1, \z_2),
         (\log \om_1, \om_2)\bigr )  \notag \\
  & = \frac{\chi_1(t)}{\z_1 \overline{\om}_1} 
 K_b\bigl((\log\z_1,\z_2),(\log\om_1,\om_2)\bigr)
\notag\\  
& \quad + 
\frac{\chi_2(t)}{\z_1 \overline{\om}_1} 
\biggl\{ e^{-h| \log(|\z_1|/|\overline{\om}_1|)|}
e^{-h\sgn t \cdot (\arg\z_1+\arg\om_1)}   
K_{\tilde b}\bigl((\log\z_1,\z_2),(\log\om_1,\om_2)\bigr)  \notag\\
& \quad
+ e^{-\nu_b| \log(|\z_1|/|\overline{\om}_1|)|}
e^{-\nu_b\sgn t \cdot (\arg\z_1+\arg\om_1)}   
\frac{\phi_1 (\log\z_1,\om_1)}{\bigl(
( \z_1/\overline{\om}_1)^{-i/2}
e^{\pi/2}- z_2 \ovw_2 \bigr)^2}\notag\\
& \quad  
+e^{-\nu_b| \log(|\z_1|/|\overline{\om}_1|)|}
e^{-\nu_b\sgn t \cdot (\arg\z_1+\arg\om_1)}   
\frac{ \phi_2 \bigl( (\log\z_1,\z_2),(\om_1,\om_2)\bigr)}
{\bigl(
( \z_1/\overline{\om}_1)^{-i/2}
e^{-\pi/2}- z_2 \ovw_2 \bigr)^2} \biggr\}
\notag\\  
  & = \chi_1(t)
\frac{K_b\bigl((\log\z_1,\z_2),(\log\om_1,\om_2)\bigr)
}{\z_1 \overline{\om}_1} \notag\\  
& \quad + 
\chi_2(t) 
\biggl\{  \biggl(\frac{|\z_1|}{|\om_1|}\biggr)^{-h\sgn t} 
e^{-h\sgn t \cdot (\arg\z_1+\arg\om_1)}   
\frac{K_{\tilde b}\bigl((\log\z_1,\z_2),(\log\om_1,\om_2)\bigr)}{ \z_1
  \overline{\om}_1}  
  \notag\\
& \quad
+ \biggl(\frac{|\z_1|}{|\om_1|}\biggr)^{-\nu_\b\sgn t} 
e^{-\nu_b\sgn t \cdot (\arg\z_1+\arg\om_1)}\biggl(    
\frac{\phi_1 (\log\z_1,\om_1)}{ \z_1\overline{\om}_1}\cdot
\frac{1}{\bigl(
( \z_1/\overline{\om_1})^{-i/2}
e^{\pi/2}- z_2 \ovw_2 \bigr)^2} \notag\\
& \quad\qquad\qquad   
+\frac{ \phi_2 \bigl( (\log\z_1,\z_2),(\om_1,\om_2)\bigr)}
{\z_1\overline{\om}_1} \cdot \frac{1}{\bigl( 
( \z_1/\overline{\om}_1)^{-i/2}
e^{-\pi/2}- z_2 \ovw_2 \bigr)^2} \biggr)  \biggr\}
\ .
\end{align}

This is the promised expression for the Bergman kernel of $D_\b$.  We
only 
need to check that the resulting functions $g_1,g_2$, 
$G_1,G_2,\dots$, $\tilde G_1,\tilde G_2,\dots$ have the announced
properties. Denoting by $G$ any of these functions, then
$$
G(\z,\om)= 
\psi \bigl( \log(\z_1 /\overline{\om}_1),\z_2\overline{\om}_2\bigr)\ ,
$$
where $\psi$ denotes any of the functions $\psi$s in Proposition
\ref{technical}, or $\Psi$s in Proposition
\ref{sum-of-error-terms-proposition}.
Then it is easy to see that $G\in{\cal
  C}^\infty(\overline{D_\b}\setminus\{0\}\times 
\overline{D_\b}\setminus\{0\})$ and that
$$
\partial_\z^\alpha \partial_{\overline{\om}}^\gamma G(\z,\om)=
{\cal O} \bigl( |\z_1|^{-|\alpha|} |\om_1|^{-|\gamma|} \bigr)
\qquad\qquad\text{as}\quad |\z_1|,|\om_1|\ra 0\ .  \eqno \BoxOpTwo
$$
\bigskip

\section{Foundational Steps in the Proofs of Propositions 7.2 and
  7.3}\label{new-section} 
\vspace*{.12in}

The core of our calculation for the expansion of the Bergman kernel is
contained in the following result. We begin the proof of this
proposition now, 
but postpone the completion of its proof until the end of the
paper---see Section \ref{appendix}.

\begin{proposition}\label{technical}  	  \sl
Let $R,S>0$ and let ${\cal D}_{R,S}$ be the domain in $\CC^2$ defined
by
$$
{\cal D}_{R,S} =
\bigl\{ (\tau,\l)\in\CC^2:\, \big|\Im\tau-\log|\l|^2\big|<S,\, 
e^{-R/2} < |\l| < e^{R/2}\bigr\}\, .
$$
Let 
$$
{\cal I}_j(\tau)=
\int_{-\infty}^\infty \sigma_{R,S}(\xi)
\bigl[
\xi^2 + b\xi+ c \bigr] 
         e^{i\tau\xi} e^{-R|\xi -\j2|} e^{-S|\xi|} \, d\xi \, ,
$$
where 
\begin{equation}\label{sigmaRS}
\sigma_{R,S}(\xi)=
\sgn(\xi)\sgn(\xi-\j2)
e^{-ih\bigl(\sgn(\xi-\j2)R+\sgn(\xi)S\bigr)}\, , 
\end{equation}
$$
b = 2ih - \j2, \quad\text{and} \quad c = ih \bigl( ih - \j2 \bigr) \, .
$$
Then there exist entire functions $\psi_j^\pm$,   $j=1,\dots,4$,
uniformly $\mathcal{O}(|\Re\tau|)$ as 
$|\Re\tau|\ra+\infty$ in any set  $\{|\Im\tau|\le C,\ |\l|\le C'\}$,  
together with all their derivatives, such
that  
\begin{align*}
\sum_{j\in\ZZ} {\cal I}_j(\tau)\lj  
& =  \frac{e^{ih(R+S)}}{(i\tau+R+S)^2} 
\biggl[
\frac{e^{R/2}}{2(e^{R/2} - \l)^2}  
   \bigl(1+\psi_1^+(\tau,\l) \bigr) \\    
&\qquad\quad
+ \frac{e^{-[i\tau +S]/2}}{2(\l - e^{-[i\tau +S]/2})^2}
   \bigl( 1 + \psi_2^-(\tau,\l) \bigr)  
+ \frac{e^{-[i\tau+S]/2}\psi_4^+(\tau)}{
2(e^{R/2} - \l)(e^{-[i\tau + S]/2}-\l)} \biggr] 
\\
& \qquad +  
\frac{e^{-ih(R+S)}}{(i\tau-R-S)^2} 
\biggl[
\frac{e^{[S- i\tau]/2}}{2(e^{[S - i\tau]/2} - \l)^2}
        \bigl(1 + \psi^+_2(\tau,\l) \bigr) \\
&\qquad\quad
+ \frac{  e^{-R/2}}{2(\l - e^{-R/2})^2}  
        \bigl(1 +\psi_1^-(\tau,\l) \bigr)  
+ \frac{e^{[S-i\tau]/2}\psi_4^-(\tau)}
{2(e^{-R/2} - \l)(e^{[S-i\tau]/2}-\l)}  \biggr] \\
&\qquad
+
\frac{e^{ih(R-S)+ R/2}}{2(e^{[S-i\tau]/2}-\l)^2(e^{R/2}-\l)^2}
\psi_3^+(\tau,\l) \\
&\qquad
+
\frac{e^{-ih(R-S)- R/2}}{2(e^{-[i\tau+S]/2}-\l)^2(e^{-R/2}-\l)^2} 
\psi_3^-(\tau,\l) 
\ .
\end{align*}

The convergence is uniform on compact subsets of
${\cal D}_{R,S}$.  Moreover, the functions $\psi_j^\pm$
can be computed explicitly (see
(\ref{psi1+psi2+}), (\ref{psi1-psi2-}), (\ref{new-psi3+}), 
(\ref{new-psi3-}),
and 
(\ref{psi4+-}) respectively).
\end{proposition}

\proof[Basic Steps in the Proof of Proposition \ref{technical}]
The proof of Proposition 
\ref{technical} is broken up into Lemmas
\ref{NEW-LEMMA-IANDIII}--\ref{NEW-LEMMA-II*}. 
We enunciate those lemmas here and complete the first part of the
proof.  The proofs 
of those technical lemmas are deferred until Section \ref{appendix}.
Also the concluding parts  
of the proof of Proposition \ref{technical} will be presented at that time.
\medskip

We wish to compute the sum
 $\sum_{j\in\ZZ}{\cal I}_j(\tau)\lj$, 
for $(\tau,\l)\in{\cal D}_{R,S}$.

We begin with the terms having index
$j\ge0$.  Writing   $\sigma_{R,S}$ explicitly,
the integral 
${\cal I}_j$ becomes 
\begin{align} 
& \int_{-\infty}^\infty  
\sgn(\xi)\sgn(\xi-\j2)
e^{-ih\bigl(\sgn(\xi-\j2)R+\sgn(\xi)S\bigr)}
\bigl[
\xi^2 + b\xi+ c \bigr] 
         e^{i\tau\xi} e^{-R|\xi -\j2|} e^{-S|\xi|} \, d\xi \notag\\
& \qquad  
= e^{ih(R+S)}  e^{-R\j2}  \biggl(
\int_{-\infty}^{-\d} [\xi^2 + b\xi + c] 
e^{(i\tau+R+S)\xi} 
    \, d\xi \notag\\ 
& \qquad\qquad\qquad\qquad\qquad\qquad\qquad\qquad
+ \int_{-\d}^\d [\xi^2 + b\xi + c] e^{(i\tau+R+S)\xi} 
    \, d\xi \biggr)
\notag\\
&  \qquad\qquad  - e^{ih(R-S)}
e^{-R\j2} \int_\d^{\j2-\d} [\xi^2 + b\xi + c] e^{(i\tau+R-S)\xi} 
    \, d\xi \notag\\
& \qquad\qquad + e^{-ih(R+S)}
e^{R\j2} \biggl(
\int_{\j2-\d}^{\j2+\d} [\xi^2 + b\xi + c] e^{(i\tau-R-S)\xi}  \, d\xi   
\notag \\ 
& \qquad\qquad\qquad\qquad\qquad\qquad\qquad\qquad
+ \int_{\j2+\d}^{+\infty} [\xi^2 + b\xi + c] e^{(i\tau-R-S)\xi}  \, d\xi
\biggr) 
    \notag \\
& \qquad
\equiv  I+{\cal E}_1 - I\!I +{\cal E}_2 + I\!I\!I 
\ . 
\label{I-II+III}
\end{align}
Here $\d$ denotes a non-negative parameter.  While in the course of the
proof of Proposition \ref{technical} the parameter $\d$ is unnecessary
(and it will be taken to be 0), we need this further refinement in the 
proof of Proposition \ref{sum-of-error-terms-proposition} in order to
control certain error terms. 
At this stage we concentrate our attention on the terms $I,I\!I$ 
and $I\!I\!I$ and we will deal with the remaining terms at a later
time. 
\medskip  \\

We shall use the following
trivial fact from calculus:
\begin{equation}\label{calculus}
\int (\xi^2 +b\xi +c) e^{\a \xi} \, d\xi =
  \left [ \frac{\xi^2}{\a} 
    + \left(\frac{b}{\a} - \frac{2}{\a^2} \right)\xi 
 + \frac{c}{\a} 
 - \frac{b}{\a^2} 
 +  \frac{2}{\a^3}\right ]  e^{\a \xi} 
\, .
\end{equation}

Notice that $(\tau,\l)\in{\cal D}_{R,S}$ implies that
$|\Im \tau | < R+S$.  Then 
the evaluation of our integrals at $\pm \infty$ will
always be 0.  

We break up the proof into a series of lemmas.

\begin{lemma}\label{NEW-LEMMA-IANDIII}	 \sl
With the notation above, there exist entire functions
$\psi_1^+$ and $\psi_2^+$, 
uniformly $\mathcal{O}(|\Re\tau|)$ as 
$|\Re\tau|\ra+\infty$ in any fixed compact set,
together with all their derivatives, 
such that
\begin{align*}  
\sum_{j \geq 0} (I +I\!I\!I) \l^j  
& =  \frac{e^{ih(R+S)+R/2}}{2(e^{R/2} - \l)^2(i\tau+R+S)^2}
e^{-\d[i\tau+R+S]} \bigl(1+\psi_1^+(\tau,\l) \bigr)     \\
& \quad + 
\frac{e^{-ih(R+S)+[S - i\tau]/2}}{2(e^{[S - i\tau]/2} - \l)^2(i\tau-R-S)^2}
        e^{-\d[R+S-i\tau]} \bigl( 1 + \psi_2^+(\tau,\l) \bigr)  \\
& \quad +
\biggl( \frac{2e^{ih(R+S)-\d[i\tau+R+S]}}{(e^{R/2} - \l)(i\tau+R+S)^3}
     - \frac{2e^{-ih(R+S)-\d[R+S-i\tau]}}{(e^{[S - i\tau]/2} -
       \l)(i\tau-R-S)^3} \biggr)  
 \, .
\end{align*}
Moreover,
for any $M>0$, 
there exist
constants $C_M$ such that, for all $(\tau,\l)$
such that
$|\Im\tau|,\, |\l|\le M$, as $|\Re\tau|$,
$R,S$ tend to $+\infty$, we have
$$
|\psi_j^+(\tau,\l)| \le C_M (|\Re\tau|+R+S) 
\qquad\qquad \text{for}\quad j=1,2\ .
$$
\end{lemma}

\begin{lemma}\label{NEW-LEMMA-II}    \sl
There exists an entire function $\psi_3^+$,
uniformly $\mathcal{O}(1)$ as 
$|\Re\tau|\ra+\infty$ in any fixed compact set,
together with all its derivatives, such
that  
$$
\sum_{j\ge0} I\!I\lj =
\frac{e^{ih(R-S)-\d[S-R-i\tau]+
    R/2}}{2(e^{[S-i\tau]/2}-\l)^2(e^{R/2}-\l)^2} 
\psi_3^+(\tau,\l)
\ .
$$
Moreover,
$|\psi_3^+ (\tau,\l)|
 \le C_M (e^{(R+S)/2} +e^S)$, for $(\tau,\l)$ in any  fixed compact set,
 as
$R,S \ra +\infty$.
\end{lemma}

Now we turn to the case $j\le0$.  In this situation the integral
${\cal I}_j(\tau)$ equals
\begin{align} 
& \int_{-\infty}^\infty 
\sgn(\xi)\sgn(\xi-\j2)
e^{-ih\bigl(\sgn(\xi-\j2)R+\sgn(\xi)S\bigr)}
\bigl[
\xi^2 + b\xi+ c \bigr]  
         e^{i\tau\xi -R|\xi -\j2| -S|\xi|} \, d\xi \notag \\
& \qquad = e^{ih(R+S)}
e^{-R(\j2)}  
\biggl(
\int_{-\infty}^{\j2-\d}  
    [\xi^2 + b\xi + c] e^{(i\tau + R+S)\xi} \, d\xi \notag \\ 
& \qquad\qquad\qquad\qquad\qquad\qquad\qquad\qquad
\int_{\j2-\d}^{\j2+\d}
    [\xi^2 + b\xi + c] e^{(i\tau - R-S)\xi} \, d\xi \biggr) \notag \\
& \qquad\qquad - e^{ih(S-R)} e^{R(\j2)}  \int_{\j2+\d}^{-\d} 
    [\xi^2 + b\xi + c] e^{(i\tau - R+S)\xi} \, d\xi \notag \\
& \qquad\qquad+ e^{-ih(R+S)}e^{R(\j2)}  
\biggl(
\int_{-\d}^{\d}
    [\xi^2 + b\xi + c] e^{(i\tau - R-S)\xi} \, d\xi \notag \\
& \qquad\qquad\qquad\qquad\qquad\qquad\qquad\qquad
\int_{\d}^{+\infty} [\xi^2 + b\xi + c] e^{(i\tau - R-S)\xi} \, d\xi 
\biggr) \notag \\
& \qquad \equiv  I^* +{\cal E}_1^*
- I\!I^* +{\cal E}_2^*+ I\!I\!I^* 
\, . \label{I-II+III*}
\end{align}
Here, as in the case $j\ge0$, we fix (the same) non-negative parameter
$\d$. For the proof of Proposition \ref{technical} we in fact take
$\d=0$.  Notice that this decomposition makes sense when $j<-1$ and
$\d$ is a small positive parameter.  The case $j=-1$ is somewhat
special, and can be treated along the same lines.  In particular, one
can simply take $\d=0$ in this case as well.

\begin{lemma}\label{NEW-LEMMA-I*ANDIII*}  \sl
There exist entire functions
$\psi_1^-$ and $\psi_2^-$, 
uniformly $\mathcal{O}(|\Re\tau|)$ as 
$|\Re\tau|\ra+\infty$ in any fixed compact set,
together with all their derivatives, 
such that
\begin{align*}  
\sum_{j <0} (I^* +I\!I\!I^*) \l^j  
& = 
\frac{  e^{-ih(R+S)-R/2}}{2(\l - e^{-R/2})^2(i\tau-R-S)^2}  
e^{-\d[R+S-i\tau]} \bigl(1 +
\psi_1^-(\tau,\l) \bigr)     \\
& \quad + \frac{e^{ih(R+S)}
e^{-[i\tau +S]/2}}{2(\l - e^{-[i\tau + S]/2})^2(i\tau+R+S)^2}
   e^{-\d[i\tau+R+S]} \bigl( 1 + \psi_2^-(\tau,\l) \bigr)  \\
& \quad+  \biggl( \frac{2e^{-ih(R+S)-\d[R+S-i\tau]}}{
(e^{-R/2} -  \l)(i\tau-R-S)^3} 
 - \frac{2e^{ih(R+S)-\d[i\tau+R+S]}}{(e^{-[i\tau + S]/2} - \l)(i\tau+R+S)^3}
 \biggr) \ .
\end{align*}
Moreover,
for any $M>0$, 
there exist
constants $C_M$ such that, for all $(\tau,\l)$
such that
$|\Im\tau|,\, |\l|\le M$, as $|\Re\tau|$,
$R,S$ tend to $+\infty$, we have
$$
|\psi_j^-(\tau,\l)| \le C_M (|\Re\tau|+R+S) 
\qquad\qquad \text{for}\quad j=1,2\ .
$$
\end{lemma}

\begin{lemma}\label{NEW-LEMMA-II*}   \sl
There exists an entire function $\psi_3^-$,
uniformly $\mathcal{O}(1)$ as 
$|\Re\tau|\ra+\infty$ in any fixed compact set,
together with all its derivatives, such
that  
$$
\sum_{j<0} I\!I^*\lj =
\frac{e^{-ih(R-S)+\d[i\tau-R+S]-R/2}}{2(e^{-[i\tau+S]/2}-\l)^2(e^{-R/2}-\l)^2} 
\psi_3^-(\tau,\l) \ .
$$
Moreover, for any $M>0$, there exists $C_M>0$ such that
for $(\tau,\l)$ satisfying $|\tau|\le M$, $|\l|\le e^M$, we have
$$
|\psi_3^- (\tau,\l)| 
 \le C_M \bigl( e^{R+S/2}   +e^{S(2\d+1/2)} 
+ e^{2R} \bigr)\ ,
$$
 as
$R,S \ra +\infty$.
\end{lemma}

These four lemmas constitute the first part of the proof of
Proposition \ref{technical}. 
We postpone the proof of Lemmas
\ref{NEW-LEMMA-IANDIII}--\ref{NEW-LEMMA-II*} to Section
\ref{appendix}, and complete 
the proof of Proposition \ref{technical} at that time.
\bigskip
	      
\section{Proof of Propositions 
\ref{sum-of-Mj-proposition} and
\ref{sum-of-error-terms-proposition} }\label{new-section-bis}   
\vspace*{.12in}

Assuming Proposition \ref{technical} for the moment, we are now in a
position to compute the sum of the main terms---Proposition
\ref{sum-of-Mj-proposition}---and of the error terms---Proposition
\ref{sum-of-error-terms-proposition}.

\proof[Proof of Proposition \ref{sum-of-Mj-proposition}]
It suffices to notice that 
$$
M_j(\tau) = \frac{2}{\pi}
e^{-\sgn(\Re\tau)h\tau}{\cal I}_j(\tau)\, ,
$$
where ${\cal I}_j(\tau)$ is as in Proposition \ref{technical} with
$R=\2bp$ and $S=\pi$. 
Moreover, the functions $\psi_j^\pm$, $j=1,\dots,4$ can be obtained by
(\ref{psi1+psi2+}), (\ref{psi1-psi2-}), 
(\ref{new-psi3+}), (\ref{new-psi3-})
and 
(\ref{psi4+-})
by taking the same values for $R$ and $S$ as above. 
\bigskip
\endpf

\proof[Proof of Proposition \ref{sum-of-error-terms-proposition}]
We now compute the sum of the error terms, that is we prove 
Proposition \ref{sum-of-error-terms-proposition}.  Recall that these
error 
terms arose from the decomposition of $J_j(\tau)$ in Section
\ref{decomposition-of-Jj}.  

We wish to evaluate the three sums $\sum_{j\in\ZZ} E_j^{(k)}(\tau)\lj$,
for $k = 1,2,3$. 

\subsection{Sum of the \boldmath $E_j^{(1)}$}

We recall that, when $\Re\tau\ge0$,  
$$
E^{(1)}_j(\tau) 
= 
\frac{2}{\pi}e^{-h\tau} 
\int_{-\infty}^{+\infty} 
    \sigma(\xi) e^{i\tau\xi}
\frac{(\xi +ih)(\xi +ih -\j2)}{e^{\pi|\xi|}e^{(\2bp)|\xi-\j2|}} 
\cdot
\frac{e^{-2\pi\sgn(\xi)(\xi+ih)}}{1-e^{-2\pi\sgn(\xi)(\xi+ih)}} 
\, d\xi\, ,
$$
where $\sigma$ is defined in
(\ref{sigma}).  At this time, for simplicity of notation, we
concentrate on the case $\Re\tau\ge0$.   
The case $\Re\tau<0$ will follow by a completely analogous argument.
 
We decompose the integral defining
$E^{(1)}_j$ as in (\ref{I-II+III}) and (\ref{I-II+III*}), according to
whether $j\ge0$ or $j<0$. Then, for a fixed $\d>0$,
\begin{equation}\label{E1+}
E^{(1)}_j(\tau) 
= 
\frac{2}{\pi}e^{-h\tau} \bigl( I-I\!I+I\!I\!I+{\cal E}_1+{\cal
  E}_2\bigr) 
\end{equation}
when $j\ge0$, and
\begin{equation}\label{E1-}
E^{(1)}_j(\tau) 
= 
\frac{2}{\pi}e^{-h\tau} \bigl( I^*-I\!I^*+I\!I\!I^*+{\cal E}_1^* +{\cal
  E}_2^*\bigr) 
\end{equation}
when $j<0$.

\subsection{The Cases of \boldmath $\sum_{j\ge0} {\cal E}_1\lj$  and
\boldmath $\sum_{j\ge0} {\cal E}_2 \lj$}

We begin by considering the sum over positive indices $\sum_{j\ge0} {\cal E}_1\lj$ (recall
that ${\cal E}_1$ depends on $j$ and $\tau$). Writing
$m_1(\xi)=e^{-2\pi\sgn(\xi)(\xi+ih)}/(1-e^{-2\pi\sgn(\xi)(\xi+ih)})$,
we have that
\begin{align*}
&\sum_{j\ge0} 
{\cal E}_1\lj \\
& = \sum_{j\ge0} \biggl(
\int_{-\d}^{\d} \sigma(\xi)
(\xi+ih)(\xi+ih-\j2) e^{i\tau\xi} e^{-\pi|\xi|}
e^{-(\2bp)|\xi-\j2|}  m_1(\xi)\, d\xi \biggr) \lj \\
& = \sum_{j\ge0} 
e^{-(\2bp)\j2} \biggl( e^{2\b ih}
\int_{-\d}^{0} (\xi+ih)^2 e^{(i\tau+2\b)\xi} m_1(\xi)\, d\xi
- 
\int_0^{\d} (\xi+ih)^2 e^{(i\tau+2\b-2\pi)\xi} m_1(\xi)\, d\xi
\\
& \qquad
-\j2 \int_{-\d}^{0} (\xi+ih)e^{(i\tau+2\b)\xi} m_1(\xi)\, d\xi 
+\j2 \int_0^{\d} (\xi+ih)e^{(i\tau+2\b-2\pi)\xi} m_1(\xi)\, d\xi 
\biggr)\lj \, .
\end{align*}
Of course we must sum in $j$.  We know that, when
$\big|e^{-\gamma/2}\l\big|<1$, 
\begin{equation}\label{summation-formulas}
\sum_{j=0}^\infty e^{-\gamma \j2} \lj = \frac{1}{e^{\gamma/2} - \l}
\qquad\text{and}
\qquad
\sum_{j=0}^\infty e^{-\gamma\j2}  \j2  \lj =
   \frac{e^{\gamma/2}}{2(e^{\gamma/2} - \l)^2} \, .
\end{equation}
Therefore
$$
\sum_{j\ge0} 
{\cal E}_1\lj
= \frac{h_1(\tau)}{e^{(\b-\pi/2)}-\l}
+\frac{h_2(\tau)}{(e^{(\b-\pi/2)}-\l)^2} 
= \frac{h_1(\tau)(e^{(\b-\pi/2)}-\l) +h_2(\tau)}{(e^{(\b-\pi/2)}-\l)^2}\ ,
$$
where $h_1$, $h_2$
are entire functions of {\it exponential type}, that is, $h_1$ and
$h_2$ decay exponentially, together with all their derivatives, in
every closed horizontal strip.

Now we turn to the sum $\sum_{j\ge0} {\cal E}_2\lj$. 
Now we notice that
$m_1(\xi)=\sum_{k=1}^{+\infty} e^{-2k\pi(\xi+ih)}$, for $\xi>0$, where
$m_1$ is as before. Notice that the series converges uniformly on
compact sets, with 
bounds  uniform in 
$j\ge0$, so we can interchange the order of 
integration and summation.
Then 
\begin{align*}
\sum_{j\ge0} 
{\cal E}_2\lj
& = \sum_{j\ge0} \biggl(
e^{-2\b ih}\int_{\j2-\d}^{\j2+\d} (\xi+ih)(\xi+ih-\j2) e^{i\tau\xi-\pi\xi
-(\2bp)(\xi-\j2)}  m_1(\xi)\, d\xi \biggr)\lj \\
& = e^{-(2\b -\pi)ih}\sum_{j\ge0}\sum_{k\ge1}  
\biggl(
\int_{\j2-\d}^{\j2+\d} (\xi+ih)(\xi+ih-\j2) e^{i\tau\xi-(\2bp)(\xi-\j2)}\\
& \qquad\qquad\qquad\qquad\qquad\qquad\qquad
\times e^{-(2k+1)\pi(\xi+ih)} \, d\xi \biggr)\lj \\
& = e^{-(2\b -\pi)ih}\sum_{j\ge0}\sum_{k\ge1}  
\biggl( e^{i\tau\j2-(2k+1)\pi\j2}
\int_{-\d}^\d (\xi+ih)(\xi+ih+\j2)
e^{i\tau\xi-(\2bp)\xi}\\
& \qquad\qquad\qquad\qquad\qquad\qquad\qquad
\times  e^{-(2k+1)\pi(\xi+ih)} \, d\xi \biggr)\lj \\
& = e^{-(2\b -\pi)ih} \sum_{k\ge1}\biggl(  
\frac{1}{e^{[(2k+1)\pi-i\tau]/2}-\l}
\int_{-\d}^\d (\xi+ih)^2 e^{i\tau\xi-(\2bp)\xi}
e^{-(2k+1)\pi(\xi+ih)} \, d\xi \\
& \qquad\qquad
+\frac{e^{[(2k+1)\pi-i\tau]/2}}{(e^{[(2k+1)\pi-i\tau]/2}-\l)^2} 
\int_{-\d}^\d (\xi+ih) e^{i\tau\xi-(\2bp)\xi}
e^{-(2k+1)\pi(\xi+ih)} \, d\xi \biggr) \\
&  = e^{-(2\b -\pi)ih} \sum_{k\ge1}\biggl(  
\frac{h_1^{(k)}(\tau)}{e^{[(2k+1)\pi-i\tau]/2}-\l} 
+\frac{h_2^{(k)}(\tau) 
e^{[(2k+1)\pi-i\tau]/2}}{(e^{[(2k+1)\pi-i\tau]/2}-\l)^2} 
\biggr)\ ,
\end{align*}
where $h_1^{(k)},h_2^{(k)}$ 
are entire functions of exponential type, with bounds uniform in $k$.
Notice that,
when summing in $j$, the series  converges for $|\l|<e^{\Im\tau+(2k+1)\pi}$.

\subsubsection*{\bf Claim 1}
We now claim that the above sum converges to a 
function $E_1(\tau,\l)$, 
holomorphic on  the domain
$$
{\cal D}_{\infty,2\pi}
=\bigl\{ (\tau,\l):\ 
\big| \Im\tau-\log|\l|^2\big| < 2\pi,\ |\l|>0 \bigr\}\, ,
$$
that is of exponential type in $\tau$
uniformly in $\l$ when $\l$ varies in a compact set.  
\medskip  \\

For, by the observation above, the functions to be summed are all
holomorphic in ${\cal D}_{\infty,2\pi}$.
Next, for fixed
$M>0$ we can select  $k_0$ large enough so that for all $k\ge k_0$,
when $(\tau,\l)\in{\cal D}_{\infty,2\pi}$, with $|\Im\tau|\le M$ and 
 $|\l|\le e^M$, we have that
$$
\big|e^{[(2k+1)\pi-i\tau]/2}-\l\big| \ge 
 e^{[(2k+1)\pi-M]/2} -e^M\ge  \frac12 e^{\pi k}\ .
$$
Therefore
$$
\bigg| \frac{h_1^{(k)}(\tau)}{e^{[(2k+1)\pi-i\tau]/2}-\l} \bigg|
\le c e^{-\pi k} \, ,
$$
and similarly
$$
\bigg|
\frac{h_2^{(k)}(\tau) 
e^{[(2k+1)\pi-i\tau]/2}}{(e^{[(2k+1)\pi-i\tau]/2}-\l)^2} \bigg|
\le c e^{-\pi k} 
\frac{e^{[(2k+1)\pi+M]/2}}{e^{[(2k+1)\pi-M]/2}-e^M}
\le c e^{-\pi k} \ ,
$$
so that the two series above converge uniformly in the
fixed compact set.  This proves the claim. 
\medskip   \\

\subsection{The Cases of \boldmath $\sum_{j<0} {\cal E}_1^*\lj$  and
\boldmath $\sum_{j<0} {\cal E}_2^* \lj$}
Now we turn to the sums over negative indices $j$.  Clearly it
suffices to consider the sum for $j\le j_0<-1$.  On the relevant
interval of integration, 
we can then write $m_1(\xi)=\sum_{k=1}^{+\infty} e^{2k\pi(\xi+ih)}$, and
the series converges uniformly there.  
We have
\begin{align} 
\sum_{j\le j_0} 
{\cal E}_1^*\lj
& = \sum_{j\le j_0} 
e^{2\b ih} \biggl( 
\int_{\j2-\d}^{\j2+\d} (\xi+ih)(\xi+ih-\j2) e^{i\tau\xi+\pi\xi+(\2bp)(\xi-\j2)}  
m_1(\xi)\, d\xi \biggr) \lj \notag \\
& = \sum_{j\le j_0}  \sum_{k\ge1}
e^{2\b ih} \biggl(
\int_{\j2-\d}^{\j2+\d} (\xi+ih)(\xi+ih-\j2)
e^{i\tau\xi+\pi\xi+(\2bp)(\xi-\j2)}\notag \\  
& \qquad\qquad\qquad
\times e^{2k\pi(\xi+ih)} \, d\xi \biggr) \lj \notag \\
& = 
e^{2\b ih} \sum_{k\ge1} \sum_{j\le j_0}
e^{[i\tau+(2k+1)\pi]\j2} 
\biggl( 
\int_{-\d}^{\d} (\xi+ih)^2 e^{i\tau\xi+(2k+1)\pi\xi+(\2bp)\xi}
e^{2kih} \, d\xi  \notag \\
& \qquad\qquad\qquad
+\j2\int_{-\d}^{\d} (\xi+ih)
e^{i\tau\xi+(2k+1)\pi\xi+(\2bp)\xi}e^{2kih} \, d\xi  \biggr)\lj \
. \label{calE1*-first-display} 
\end{align}
Again we must sum in $j$.  Since
\begin{equation}\label{neg-sum-1}
\sum_{j\le j_0} e^{\gamma \j2} \lj = \frac{(e^{\gamma/2}\l)^{j_0+1}}{\l-e^{-\gamma/2}}
\end{equation}
and
\begin{equation}\label{neg-sum-2}
\sum_{j\le j_0} \j2 e^{\gamma\j2}  \lj = - 
   \frac{(e^{\gamma/2}\l)^{j_0+1} 
\bigl( (j_0+1)(e^{-\gamma/2}-\l)+e^{-\gamma/2}\bigr)
}{2(\l- e^{-\gamma/2})^2} \, ,
\end{equation}
we have, when $\big|e^{\gamma/2}\l\big|>1$, i.e. when $|\l|>e^{-\gamma/2}$, that 
\begin{align}
\sum_{j\le j_0} 
{\cal E}_1^*\lj
& = 
e^{2\b ih} \sum_{k\ge1} 
\biggl( e^{2kih} 
\frac{(\l-e^{[i\tau+(2k+1)\pi]/2})^{j_0+1}}{\l-e^{-[i\tau+(2k+1)\pi]}}  
\int_{-\d}^{0} (\xi+ih)^2 e^{i\tau\xi+(2k+1)\pi\xi+(\2bp)\xi}
 \, d\xi  \notag \\
& \qquad
+e^{2kih}  \frac{(e^{[i\tau+(2k+1)\pi]/2}\l)^{j_0+1} \bigl(
  (j_0+1)(e^{-[i\tau+(2k+1)\pi]/2}-\l)
+e^{-[i\tau+(2k+1)\pi]/2}\bigr)
}{2(\l- e^{-[i\tau+(2k+1)\pi]/2})^2} \notag \\
& \qquad\qquad
\times \int_{-\d}^{0} (\xi+ih)
e^{i\tau\xi+(2k+1)\pi\xi+(\2bp)\xi}
\, d\xi  \biggr) \notag \\
& =\sum_{k\ge1} 
\biggl( \frac{h^{(k)}_3(\tau,\l)}{\l-e^{-[i\tau+(2k+1)\pi]/2}}
+\frac{h^{(k)}_4(\tau,\l)}{2(\l- e^{-[i\tau+(2k+1)\pi]/2})^2} 
\biggr)  \ , \label{calE1*-second-display} 
\end{align}
where $h_3^{(k)},h_4^{(k)}$ are entire functions that are of
exponential type in $\tau$, uniformly in $k$ and $\l$, for $\l$ varying
in any compact set.  Notice that the above series in $j$ converge when
$|\l|> e^{[\Im\tau-(2k+1)]/2}$ hence, in particular, when $(\tau,\l)\in
{\cal D}_{\infty,2\pi}$.

Arguing as in Claim 1 (Subsection 10.2), we see that the above sum
converges to a
function $E_2(\tau,\l)$, 
holomorphic on  ${\cal D}_{\infty,2\pi}$, 
which is of exponential type in $\tau$, uniformly in
$\l$ for $\l$ varying in any compact sets.  
For, having fixed an
$M>0$, we can select $k_0$ such that if $k\ge k_0$, if $|\Im\tau|\le M$ and
$|\l|\le e^M$, 
$$
\bigg|
\frac{1}{(e^{[(2k+1)\pi-i\tau]/2}-\l)^2} \bigg|
\le c e^{-\pi k} 
\qquad\text{and}\quad |h_3^{(k)}(\tau,\l)|,\,|h_4^{(k)}(\tau,\l)| \le C
\ .
$$

Now we turn to the sum 
\begin{align*}
\sum_{j\le-1}{\cal E}_2^*\lj
& =    \sum_{j\le-1} e^{2\b ih} 
\int_{-\d}^\d (\xi+ih)(\xi+ih-\j2)
e^{i\tau\xi -\pi\xi-(\2bp)(\xi-\j2)} m_1(\xi)\, d\xi \lj \\
& = e^{2\b ih} \sum_{j\le-1} 
e^{(\2bp)\j2} \biggl( \int_{-\d}^\d (\xi+ih)^2 
e^{(i\tau-\2bp)\xi} m_1(\xi)\, d\xi\\
& \qquad\qquad\qquad\qquad
-\j2 \int_{-\d}^\d (\xi+ih)
e^{(i\tau-\2bp)\xi} m_1(\xi)\, d\xi \biggr)\lj \ .
\end{align*}
Using formulas (\ref{neg-sum-1}) and
(\ref{neg-sum-2}) 
with $j_0=-1$ we see that, for $|\l|>e^{-(\b-\pi/2)}$,
$$
\sum_{j<0} 
{\cal E}_2^*\lj
= \frac{g_1(\tau)}{\l-e^{-(\b-\pi/2)}}
+\frac{g_2(\tau)}{(\l-e^{-(\b-\pi/2)})^2} 
= \frac{g_1(\tau)(\l-e^{-(\b-\pi/2)}) +g_2(\tau)}{(\l-e^{-(\b-\pi/2)})^2}\ ,
$$
where $g_1$, $g_2$
are entire functions of exponential type, that is, they decay
exponentially, together with all their derivatives, in every closed
horizontal strip.
\medskip   \\

According to the decompositions (\ref{E1+}) and (\ref{E1-}), 
in order to conclude the analysis of the term
$\sum_{j\in\ZZ}E_j^{(1)}\lj$, it remains to consider the sums of the
terms involving $I$ through $I\!I\!I$ and $I^*$ through $I\!I\!I^*$,
that is, 
$$
\sum_{j\ge0} \bigl(I-I\!I+I\!I\!I\bigr)\lj +
\sum_{j<0} \bigl(I^*-I\!I^*+I\!I\!I^*\bigr)\lj \ .
$$

\subsection{The Case of 
\boldmath $\sum_{j\ge0} \bigl(I+I\!I\!I\bigr)\lj + 
\sum_{j<0} \bigl(I^*+I\!I\!I^*\bigr)\lj $}  
Notice that, when $\Re\tau\ge0$, since $h$ is fixed with $0<h<1/2$, 
the
series $\sum_{k\ge1} e^{-2k\pi\sgn(\xi)(\xi+ih)}$ converging  
  to $m_1(\xi)$ has partial sums that are uniformly bounded.
  Thus, by the Lebesgue dominated convergence theorem,
we can interchange the order of summation and integration. 
Recalling  the definitions of $\sigma$ and
$\sigma_{R,S}$, (\ref{sigma}) and
  (\ref{sigmaRS}), we then write 
\begin{align}
& \sum_{j\ge0} \bigl(I+I\!I\!I\bigr)\lj \notag \\ 
& =  
\sum_{j\ge0} \lj \sum_{k\ge1}   
\biggl( \int_{-\infty}^{-\d} +\int_{\j2+\d}^{+\infty} 
\biggr)  
\sigma(\xi)(\xi +ih)(\xi +ih -\j2)\notag \\ 
& \qquad\qquad\qquad\times 
e^{i\xi\tau} e^{-\pi|\xi|}e^{-(\2bp)|\xi-\j2|}  
e^{-2\pi k\sgn(\xi)(\xi+ih)} 
\,d\xi \notag \\
& = \sum_{k=1}^{+\infty}  
\sum_{j\ge0} \lj 
\biggl( \int_{-\infty}^{-\d} +\int_{\j2+\d}^{+\infty} 
\biggr)  
\sigma_{(\2bp),(2k+1)\pi}(\xi)(\xi +ih)(\xi +ih -\j2)\notag \\ 
& \qquad\qquad\qquad\times 
e^{i\xi\tau}
e^{-(2k+1)\pi|\xi|}
e^{-(\2bp)|\xi-\j2|} \, d\xi \notag \\
& = \sum_{k=1}^{+\infty}  \sum_{j\ge0} \lj {\cal I}_j^{(k)}(\tau)
\, .  \label{I+III-first-error-estimate}
\end{align} 

Analogous arguments apply when $j<0$ to give
\begin{align} 
& \sum_{j<0} \bigl(I^* +I\!I\!I^*\bigr)\lj \notag \\ 
& =  
\sum_{j<0} \lj \sum_{k\ge1}   
\biggl( \int_{-\infty}^{\j2-\d}  +\int_{\d}^{+\infty} 
\biggr)  
\sigma(\xi)(\xi +ih)(\xi +ih -\j2)\notag \\ 
& \qquad\qquad\qquad\times 
e^{i\xi\tau} e^{-\pi|\xi|}e^{-(\2bp)|\xi-\j2|}  
e^{-2\pi k\sgn(\xi)(\xi+ih)} 
\,d\xi \notag \\
& = \sum_{k=1}^{+\infty}  
\sum_{j<0} \lj 
\biggl( \int_{-\infty}^{\j2-\d}  +\int_{\d}^{+\infty} 
\biggr)  
\sigma_{(\2bp),(2k+1)\pi}(\xi)(\xi +ih)(\xi +ih -\j2)\notag \\ 
& \qquad\qquad\qquad\times 
e^{i\xi\tau}
e^{-(2k+1)\pi|\xi|}
e^{-(\2bp)|\xi-\j2|} \, d\xi \notag \\
& = \sum_{k=1}^{+\infty}  \sum_{j<0} \lj {\cal I}_j^{(k)}(\tau)
\, . \label{I+III-second-error-estimate}
\end{align}

In order to calculate the sum $\sum_{j\in\ZZ} {\cal I}_j^{(k)}(\tau)\lj$, 
we apply 
Lemma
\ref{NEW-LEMMA-I*ANDIII*} and the identity
(\ref{I-II+III*}) 
with 
$R=\2bp$, $S=(2k+1)\pi$ and $\d>0$ a fixed parameter.
The double sum above equals
\begin{align}
& \sum_{k\ge1}
\frac{e^{ih(\2bp+S)-\d[i\tau+\2bp+S]}}{(i\tau+\2bp+S)^2} 
\biggl[
\frac{e^{(\2bp)/2}}{2(e^{(\2bp)/2} - \l)^2}  
   \bigl(1+\psi_1^+(\tau,\l) \bigr)  \notag \\    
&\qquad\qquad
+ \frac{e^{-[i\tau +S]/2}}{2(\l - e^{-[i\tau +S]/2})^2}
   \bigl( 1 + \psi_2^-(\tau,\l) \bigr)  
+ \frac{e^{-[i\tau+S]/2}\psi_4^+(\tau)}{
2(e^{(\2bp)/2} - \l)(e^{-[i\tau + S]/2}-\l)} \biggr] 
\notag \\
& \qquad +  
\frac{e^{-ih(\2bp+S)-\d[\2bp+S-i\tau]}}{\bigl(i\tau-(\2bp)-S\bigr)^2} 
\biggl[
\frac{e^{[S- i\tau]/2}}{2(e^{[S - i\tau]/2} - \l)^2}
        \bigl(1 + \psi^+_2(\tau,\l) \bigr) \notag \\
&\qquad\qquad
+ \frac{  e^{-(\2bp)/2}}{2(\l - e^{-(\2bp)/2})^2}  
        \bigl(1 +\psi_1^-(\tau,\l) \bigr)  
+ \frac{e^{[S-i\tau]/2}\psi_4^-(\tau)}
{2(e^{-(\2bp)/2} - \l)(e^{[S-i\tau]/2}-\l)}  \biggr] 
\ , \label{estimate-double-sum-error-term}
\end{align}
where $S=(2k+1)\pi$ and
$\psi_j^\pm\equiv\psi_{(\2bp,(2k+1)\pi),j}^\pm$ depend also on the
parameter $\d$ and
are defined in 
(\ref{psi1+psi2+}), (\ref{psi1-psi2-}), 
and 
(\ref{psi4+-}).  Here we clearly need to stress the dependence on
$S$. 

We are going to show that the functions depending on $k$, that is
on $S$, can be summed, and their sums are functions of $(\tau,\l)$,
holomorphic in  a neighborhood of $\overline{{\cal D}}$, bounded 
(together with their derivatives) as
$|\Re\tau| \rightarrow+\infty$.
\medskip   \\

We let $(\tau,\l)$ vary in the closure of the domain
$$ 
{\cal D}_{2\b,3\pi/2}=
\bigl\{ \big|\Im\tau -\log|\l|^2 \big|< 3\pi/2,\
e^{-\b}<|\l|<e^{\b} \bigr\}\, , 
$$
a domain that contains the closure of the domain
${\mathcal D}$. 
Using 
Lemma \ref{NEW-LEMMA-IANDIII}, we now obtain that
\begin{align} 
\bigg| \frac{e^{ih(\2bp+S)-\d[i\tau+\2bp+S]}}{(i\tau+\2bp+S)^2} 
\bigg|\cdot
\big| e^{(\2bp)/2}  \bigl(1+\psi_1^+(\tau,\l) \bigr) \big|
& \le  C \bigl( |\Re\tau|+S) \frac{e^{-\d[i\tau+R+S]}}{|i\tau+R+S|^2} 
\notag \\    
& \le C e^{-\d S}\  \label{est-E1-psi1+} 
\end{align}
as $S\rightarrow+\infty$, uniformly in $(\tau,\l)\in{\cal D}_{2\b,3\pi/2}$.

Next notice that the functions
$
\frac{e^{-[i\tau +S]/2}}{2(\l - e^{-[i\tau +S]/2})^2}
$
are bounded for  $(\tau,\l)\in{\cal D}_{2\b,3\pi/2}$, uniformly in
$S=(2k+1)\pi$. 
Then 
\begin{align}
& \bigg|
\frac{e^{ih(\2bp+S)-\d[i\tau+\2bp+S]}}{(i\tau+\2bp+S)^2} 
\cdot
\frac{e^{-[i\tau +S]/2}}{2(\l - e^{-[i\tau +S]/2})^2}
   \bigl( 1 + \psi_2^-(\tau,\l) \bigr)  \bigg| \notag \\
& \qquad\qquad\qquad
\le C e^{-\d S} \frac{|\Re\tau|+S}{|i\tau+(\2bp)+S|^2} \notag \\
& \qquad\qquad\qquad
\le C e^{-\d S}  \label{est-E1-psi2-}
\end{align}
as $S\rightarrow+\infty$, uniformly in $(\tau,\l)\in{\cal D}_{2\b,3\pi/2}$.

Arguing as above, we see that 
\begin{align}
\bigg|
\frac{e^{ih(\2bp+S)-\d[i\tau+\2bp+S]}}{(i\tau+\2bp+S)^2} 
\cdot  \frac{e^{-[i\tau+S]/2}\psi_4^+(\tau)}{
2(e^{-[i\tau + S]/2}-\l)} \biggr| 
& \le C  e^{-\d S} \frac{1}{|i\tau+(\2bp)+S|^2} \notag \\
& \le C e^{-\d S}\ , \label{est-E1-psi4+}
\end{align}
again, 
as $S\rightarrow+\infty$, uniformly in $(\tau,\l)\in{\cal D}_{2\b,3\pi/2}$.

Next we notice that the functions
$
\frac{e^{[S-i\tau]}}{2(e^{[S-i\tau]/2}-\l)^2}
$
are bounded for  $(\tau,\l)\in{\cal D}_{2\b,3\pi/2}$ uniformly in
$S=(2k+1)\pi$. For it suffices to notice that
$$
\big| \l e^{-[S-i\tau]/2}\big|^2 \le e^{-3\pi/2}
$$
for all $(\tau,\l)\in{\cal D}_{2\b,3\pi/2}$ and all $k\ge 1$.

Then, using Lemma \ref{NEW-LEMMA-IANDIII} again, we have
\begin{align}
& \bigg| 
\frac{e^{-ih(\2bp+S)-\d[\2bp+S-i\tau]}}{\bigl(i\tau-(\2bp)-S\bigr)^2} 
\cdot 
\frac{e^{[S- i\tau]/2}}{2(e^{[S - i\tau]/2} - \l)^2}
        \bigl(1 + \psi^+_2(\tau,\l) \bigr) \bigg| 
\notag \\
& \qquad\qquad\qquad
 \le C e^{-\d S} \frac{(|\Re\tau| +S)e^{-S/2}}{|i\tau-(\2bp)-S|^2} 
\notag \\
&  \qquad\qquad\qquad
\le C e^{-(1+\d) S}\ , \label{est-E1-psi2+} 
\end{align}
\begin{align}
 \bigg| 
\frac{e^{-ih(\2bp+S)-\d[\2bp+S-i\tau]}}{\bigl(i\tau-(\2bp)-S\bigr)^2} 
(1+\psi_1^-(\tau,\l)
\bigg| 
& \le  C e^{-\d S} \frac{|\Re\tau|+S}{|i\tau-(\2bp)-S|^2} 
\notag\\
& \le C e^{-\d S}\ , \label{est-E1-psi1-} 
\end{align}
and
\begin{align}
\bigg| 
\frac{e^{-ih(\2bp+S)-\d[\2bp+S-i\tau]}}{\bigl(i\tau-(\2bp)-S\bigr)^2} 
\cdot 
\frac{e^{[S-i\tau]/2}\psi_4^-(\tau)}{
2(e^{[S-i\tau]/2}-\l)} \biggr| 
&
\le C 
\frac{e^{-\d S} }{|i\tau-(\2bp)-S|^2} \notag \\
& 
\le C e^{-\d S} \ .\label{est-E1-psi4-} 
\end{align}

Next we have
\subsection{The Case of 
\boldmath $\sum_{j\ge0} I\!I\lj + 
\sum_{j<0} I\!I^* \lj $}\label{summation-by-parts}  

In this case we will
use Lemmas \ref{NEW-LEMMA-II} and \ref{NEW-LEMMA-II*},
and  summation by parts.  
We begin with the case $j\ge0$.
Notice that, for $\xi\ge\d$, for any positive integer $N$,
$
m_1(\xi) 
= \sum_{k\ge1}^N e^{-2k\pi (\xi+ih)} 
+ e^{-2N\pi(\xi+ih)}m_1(\xi) 
$.
Therefore
\begin{align*}
& \sum_{j\ge0} I\!I \lj \\
& = e^{ih(2\b-2\pi)} \sum_{j\ge0} e^{-(2\b-\pi)\j2}\lj
\int_{\d}^{\j2-\d} \bigl[ \xi^2 +b\xi +c\bigr] 
e^{(i\tau+(\2bp)-\pi)\xi} m_1(\xi)\, d\xi\\
& =  e^{ih(2\b-2\pi)} \sum_{j\ge0} e^{-(2\b-\pi)\j2}
\lj \biggl( \sum_{k=1}^N
e^{-2\pi i kh} 
\int_{\d}^{\j2-\d} \bigl[ \xi^2 +b\xi +c\bigr] 
e^{(i\tau+(\2bp)-(2k+1)\pi)\xi} \, d\xi\\
& \qquad\qquad + \int_{\d}^{\j2-\d} \bigl[ \xi^2 +b\xi +c\bigr] 
e^{(i\tau+(\2bp)-(2N+1)\pi)\xi} m_1(\xi)\, d\xi \biggr) \\
& \equiv A+B\ .
\end{align*}

Now we apply Lemma \ref{NEW-LEMMA-II} with
$R=\2bp$ and $S=(2k+1)\pi$, $k=1,\dots,N$.  Then we  see that
there exists a function $\Psi_1(\tau,\l)$ holomorphic on the closure of 
${\cal D}_{2\b,3\pi/2}$, bounded 
(together with all their derivatives)
as $|\Re\tau|\rightarrow+\infty$,
such that 
$$
A= \frac{\Psi_1(\tau,\l)}{(e^{\bmp} -\l)^2} \ . 
$$

In order to evaluate $B$ we use the summation by parts formula
\begin{equation}\label{summation-by-parts-formula}
\sum_{j=0}^N a_j b_j = \sum_{k=0}^{N-1} s_k (b_{k+1} -b_k) +s_N b_N-
s_0b_0\ ,
\end{equation}
where $s_k=\sum_{j=0}^k a_j$, $a_j=e^{-(\2bp)\j2} \lj$ and
$$
b_j= \int_{\d}^{\j2-\d} \bigl[ \xi^2 +b\xi +c\bigr] 
e^{(i\tau+(\2bp)-(2N+1)\pi)\xi} m_1(\xi)\, d\xi\ .
$$

Notice that
$$
s_k =\sum_{j=0}^k e^{-(\2bp)\j2} \lj = 
\frac{1-\l^k e^{-k\bmp}}{e^{\bmp} -\l}\ ,
$$
and that
\begin{align*}
b_{k+1}-b_k 
& = \biggl(\int_{\d}^{\frac{k+3}{2} -\d} -\int_{\d}^{\frac{k+1}{2}-\d} 
\biggr) 
\bigl[ \xi^2 +b\xi +c\bigr] 
e^{(i\tau+(\2bp)-(2N+1)\pi)\xi} m_1(\xi)\, d\xi  \\
& =   \int_{1-\d}^{2-\d} 
\bigl[ \bigl(\xi+\frac{k-1}{2}\bigr)^2 
+b\bigl(\xi+\frac{k-1}{2}\bigr) +c\bigr] 
e^{(i\tau+(\2bp)-(2N+1)\pi)(\xi+\frac{k-1}{2})} \\
& \qquad\qquad\qquad\qquad\times
 m_1\bigl(\xi+\frac{k-1}{2}\bigr)\, 
d\xi  \\
& = e^{(i\tau+(\2bp)-(2N+1)\pi)(\frac{k-1}{2})} 
\int_{1-\d}^{2-\d} 
\bigl[ \bigl(\xi+\frac{k-1}{2}\bigr)^2 
+b\bigl(\xi+\frac{k-1}{2}\bigr) +c\bigr] \\
& \qquad\qquad\qquad\qquad\times
e^{(i\tau+(\2bp)-(2N+1)\pi)(\xi+\frac{k-1}{2})} 
m_1\bigl(\xi+\frac{k-1}{2}\bigr)\, d\xi \ .
\end{align*}
Therefore
\begin{multline}
\sum_{k=0}^{N-1} 
s_k (b_{k+1}-b_k )  \\
 = \sum_{k=0}^{N-1} \frac{1-\l^k e^{-k\bmp}}{e^{\bmp} -\l}
e^{(i\tau+(\2bp)-(2N+1)\pi)(\frac{k-1}{2})} \bigl( 
c_2 k^2 +c_1 k +c_0 \bigr) \ , \label{Psi2}
\end{multline}
where $c_j=c_j(\tau,k)$ are entire functions of $\tau$, of exponential type,
with bounds uniform in $k$.  

Now it is not hard to see that, as we let $N\ra+\infty$, the term
on the righthand side of equation (\ref{Psi2})
above converges to a 
function $\Psi_2(\tau,\l)$, holomorphic in a neighborhood of
the
closure of ${\cal D}_{2\b,3\pi/2}$, times $1/(e^{\bmp}-\l)$.
\medskip    \\

The sum for $j<0$ is treated similarly, and we do not repeat the
argument here.  Hence (relabeling the functions) 
we obtain that there exist functions $\Psi_1,\Psi_2$, holomorphic on the
closure of ${\cal D}_{2\b,3\pi/2}$, bounded 
(together with all their derivatives)
as $|\Re\tau|\rightarrow+\infty$,
such that the 
sum $ \sum_{j\ge0} \bigl(I\!I\bigr)\lj + 
\sum_{j<0} \bigl(I\!I^*\bigr)\lj$, 
relative to the error term $E_1$,
equals
$$
\frac{\Psi_1(\tau,\l)}{(e^{R/2}-\l)^2}
+ \frac{\Psi_2(\tau,\l)}{(e^{-R/2}-\l)^2} \ .
$$

This proves the statement for the error term $E_1$.
\medskip    \\

\subsection{The Sum of the \boldmath $E_j^{(2)}$}

The analysis of the sum $\sum_{j\in\ZZ}E_j^{(2)}\lj$ is quite similar
to the previous case of $E_j^{(1)}$.  In fact, the change of variables
$\xi'=\xi+\j2$ and rescaling allow one to apply the previous proof
to this case.

\subsection{The Sum of the \boldmath $E_j^{(3)}$}
Finally, we consider the sum 
$\sum_{j\in\ZZ} E_j^{(3)}\lj$.  Recall that 
\begin{multline*}
E_j^{(3)}(\tau) 
= \frac{2}{\pi}e^{-h\sgn(\Re\tau)\tau} 
\int_{-\infty}^{+\infty}
    \sigma(\xi) e^{i\tau\xi}
\frac{(\xi +ih)(\xi +ih -\j2)}{e^{\pi|\xi|}e^{(\2bp)|\xi-\j2|}} \\
\times
\frac{e^{-2\sgn(\xi)\pi(\xi+ih)}}{1-e^{-2\sgn(\xi)\pi(\xi+ih)}} 
\cdot
\frac{e^{-2\sgn(\xi-\j2)(\2bp)(\xi-\j2+ih)}}{1
  -e^{-2\sgn(\xi-\j2)(\2bp)(\xi-\j2+ih)}}  
\, d\xi\, . 
\end{multline*}

Again, we begin with the case $\Re\tau\ge0$.  
\medskip    \\

We divide the integral above as in 
(\ref{I-II+III}).  We write
\begin{align*}
& \int_{-\infty}^{+\infty}
    \sigma(\xi) e^{i\tau\xi} m_1(\xi)\tilde m(\xi)
\frac{(\xi +ih)(\xi +ih -\j2)}{e^{\pi|\xi|}e^{(\2bp)|\xi-\j2|}} 
\, d\xi\, \\
& = \biggl(
\int_{-\infty}^{-\d} + \int_{-\d}^\d +\int_\d^{\j2-\d} +
\int_{\j2-\d}^{\j2+\d}+\int_{\j2+\d}^{+\infty} \biggr)\\
& \qquad \qquad \qquad \qquad \times
\sigma(\xi) m_1(\xi)\tilde m(\xi)
\frac{(\xi +ih)(\xi +ih -\j2)e^{i\tau\xi} }{e^{\pi|\xi|}e^{(\2bp)|\xi-\j2|}} 
\, d\xi\, \\
& \equiv  I+{\cal E}_1 - I\!I +{\cal E}_2 + I\!I\!I 
\ ,
\end{align*}
where, as before,
$$
m_1(\xi)
=
\frac{e^{-2\sgn(\xi)\pi(\xi+ih)}}{1-e^{-2\sgn(\xi)\pi(\xi+ih)}} 
\qquad\text{and}\quad
\tilde m(\xi)
\frac{e^{-2\sgn(\xi-\j2)(\2bp)(\xi-\j2+ih)}}{1
  -e^{-2\sgn(\xi-\j2)(\2bp)(\xi-\j2+ih)}}  
\ .
$$

We first consider the sum involving $\mathcal{E}_1$:
\begin{align*}
&\sum_{j\ge0} {\cal E}_1\lj \\
& = e^{2\b ih} \sum_{j\ge0} \biggl(
\int_{-\d}^\d m_1(\xi)
(xi+ih)\bigl(\xi+ih-\j2\bigr) e^{i\tau\xi+\pi\xi+(\2bp)(\xi-\j2)} \tilde m(\xi)
\, d\xi \biggl) \lj \\
& = e^{2\b ih} \sum_{j\ge0} 
\biggl(
\sum_{k\ge1} e^{2k(\2bp)ih} 
\int_{-\d}^\d m_1(\xi)
(\xi+ih)\bigl(\xi+ih-\j2\bigr) \\
& \qquad\qquad\qquad \times
e^{i\tau\xi+\pi\xi+(\2bp)(2k+1)(\xi-\j2)} 
\, d\xi \biggl) \lj \\
\end{align*}

\begin{align*}
& = e^{2\b ih} 
\sum_{k\ge1} e^{2k(\2bp)ih} 
\sum_{j\ge0} e^{-(\2bp)(2k+1)\j2} \biggl(
\int_{-\d}^\d m_1(\xi) (\xi+h)^2 e^{i\tau\xi+\pi\xi+(\2bp)(2k+1)\xi} 
\, d\xi \\
& \qquad\qquad\qquad -\j2 \int_{-\d}^\d m_1(\xi) (\xi+h)
e^{i\tau\xi+\pi\xi+(\2bp)(2k+1)\xi}  \, d\xi \biggr) \lj\\
& = e^{2\b ih} 
\sum_{k\ge1} e^{2k(\2bp)ih} 
\biggl( \frac{1}{e^{(\bmp)(2k+1)}-\l} 
\int_{-\d}^\d m_1(\xi) (\xi+h)^2 
e^{i\tau\xi+\pi\xi+(\2bp)(2k+1)\xi}\, d\xi \\
& \qquad\qquad\qquad - 
\frac{e^{(\bmp)(2k+1)}}{2(e^{(\bmp)(2k+1)}-\l)^2}  
\int_{-\d}^\d m_1(\xi) (\xi+h)
e^{i\tau\xi+\pi\xi+(\2bp)(2k+1)\xi}  \, d\xi \biggr) .
\end{align*} 
The last series converge when $|\l|<e^{\bmp(2k+1)}$.

Now we fix $0<\d<1$. Then the series
$$
\sum_{k\ge1} e^{2k(\2bp)ih} 
\frac{e^{(\2bp)(2k+1)\xi}}{e^{(\bmp)(2k+1)}-\l}\, ,
\qquad\text{and}\qquad
\sum_{k\ge1} e^{2k(\2bp)ih} 
\frac{e^{(\2bp)(2k+1)(\xi+1)}}{2(e^{(\bmp)(2k+1)}-\l)^2}
$$
converge uniformly for $\xi\in[-\d,\d]$ to functions of $\xi$ and $\l$,
$C^\infty$ in $\xi$ and holomorphic in $\l$ for $|\l|<e^{3\bmp}$
Therefore we obtain that 
$\sum_{j\ge0} {\cal E}_1\lj$ converges to a function $E(\tau,\l)$,
holomorphic for $|\l|<e^{3\bmp}$,
entire in $\tau$ and of exponential decay in any closed horizontal
strip, uniformly in $\l$ for $|\l|\le C<e^{3\bmp}$.
\medskip    \\

We now consider the sum involving $\mathcal{E}_2$.  We have
\begin{align*}
&\sum_{j\ge0} {\cal E}_2\lj \\
& = e^{-2\b ih} \sum_{j\ge0} \int_{\j2-\d}^{\j2+\d} m_1(\xi) \tilde m(\xi)
(\xi+ih)\bigl(\xi+ih-\j2\bigr) e^{i\tau\xi-\pi\xi-(\2bp)(\xi-\j2)} 
\, d\xi \lj \\
& = 
e^{-2\b ih} \sum_{j\ge0} 
\sum_{k\ge1} e^{-2k\pi ih} 
\int_{\j2-\d}^{\j2+\d} \tilde m(\xi)
(\xi+ih)\bigl(\xi+ih-\j2\bigr) e^{i\tau\xi-(2k+1)\pi\xi-(\2bp)(\xi-\j2)} 
\, d\xi \lj \\
& = e^{-(2\b-\pi) ih} \sum_{j\ge0} 
\sum_{k\ge1} e^{-(2k+1)\pi ih} 
e^{(i\tau-(2k+1)\pi)\j2}\biggl(
\int_{-\d}^\d m_2(\xi) (\xi+ih)^2  e^{i\tau\xi-(2k+1)\pi\xi-(\2bp)\xi} \, d\xi \\
& \qquad\qquad +\j2 
\int_{-\d}^\d m_2(\xi) (\xi+ih)  e^{i\tau\xi-(2k+1)\pi\xi-(\2bp)\xi} \, d\xi 
\biggr) \lj \\
& = e^{-(2\b-\pi) ih} \sum_{k\ge1} e^{-(2k+1)\pi ih} 
\biggl( \frac{1}{ e^{[(2k+1)\pi-i\tau]/2} -\l}
\int_{-\d}^\d m_2(\xi) (\xi+ih)^2  e^{i\tau\xi-(2k+1)\pi\xi-(\2bp)\xi} \, d\xi \\
& \qquad\qquad + \biggl(
\frac{e^{[(2k+1)\pi-i\tau]/2}}{2(e^{[(2k+1)\pi-i\tau]/2} -\l)^2} 
\int_{-\d}^\d m_2(\xi) (\xi+ih)  e^{i\tau\xi-(2k+1)\pi\xi-(\2bp)\xi} \, d\xi 
\biggr) \ ,
\end{align*}
where we have set
$\tilde m(\xi+\j2)\equiv m_2(\xi)$, since it does not depend on $j$.  
When summing in $j$, the series  converges for $|\l|<e^{\Im\tau+(2k+1)\pi}$.

The term on the
righthand side above equals 
$$
e^{-(2\b -\pi)ih} \sum_{k\ge1}\biggl(  
\frac{h_1^{(k)}(\tau)}{e^{[(2k+1)\pi-i\tau]/2}-\l} 
+\frac{h_2^{(k)}(\tau) 
e^{[(2k+1)\pi-i\tau]/2}}{(e^{[(2k+1)\pi-i\tau]/2}-\l)^2} 
\biggr)\ ,
$$
where $h_1^{(k)},h_2^{(k)}$ 
are entire functions of exponential type, with bounds uniform in $k$.
Arguing as in Claim 1 (Subsection 10.2), we obtain that
the above sum converges to a 
function $E_2(\tau,\l)$ 
holomorphic in the domain
${\cal D}_{\infty,2\pi}$
that is of exponential type in $\tau$,
uniformly in $\l$, when $\l$ varies in a compact set.  
\medskip  \\

Now we consider the sums for the other two error terms, ${\cal E}_1^*$
and ${\cal E}_2^*$.  These sums involve $j$ when $j<0$.

We have
\begin{multline*}
\sum_{j<0} 
{\cal E}_1^*\lj
= \sum_{j<0} 
e^{2\b ih} \biggl( 
\int_{\j2-\d}^{\j2+\d} m_1(\xi) \tilde m(\xi)
(\xi+ih)(\xi+ih-\j2) \\
  \times
e^{i\tau\xi+\pi\xi+(\2bp)(\xi-\j2)}  
\, d\xi \biggr) \lj \ ,
\end{multline*}
and $m_1(\xi)\equiv\tilde m(\xi+\j2)$ does not depend on $j$. 
We can
argue as in (\ref{calE1*-first-display}) and
 (\ref{calE1*-second-display}) to obtain that
the above sum
converges to a
function $E_3(\tau,\l)$, 
holomorphic on the closure of ${\cal D}_{\infty,2\pi}$, 
which is of exponential type in $\tau$, uniformly in
$\l$ for $\l$ varying in any compact set.  

Next, expanding $\widetilde m(\xi)$ in Taylor series, we see that
\begin{align*}  
& \sum_{j\le-1}{\cal E}_2^*\lj \\
& =   e^{2\b ih}  \sum_{j\le-1} 
\int_{-\d}^\d (\xi+ih)(\xi+ih-\j2)m_1(\xi)\tilde m(\xi)
e^{i\tau\xi -\pi\xi-(\2bp)(\xi-\j2)} \, d\xi \ \lj \\
& = \sum_{k\ge1}  e^{-(2k-1)(\2bp)ih}  \sum_{j\le-1} 
e^{ (\2bp)(2k+1)\j2} \biggl( \int_{-\d}^\d (\xi+ih)^2 m_1(\xi)
e^{i\tau\xi -\pi\xi-(\2bp)(2k+1)\xi} \, d\xi \\
& \qquad\qquad\qquad\qquad -\j2
\int_{-\d}^\d (\xi+ih) m_1(\xi)
e^{i\tau\xi -\pi\xi-(\2bp)(2k+1)\xi} \, d\xi \biggr) \lj \\
\end{align*}

\begin{align*}
& = \sum_{k\ge1}  e^{-(2k-1)(\2bp)ih}  
\frac{1}{\l-e^{-(\2bp)(2k+1)/2}} 
\biggl( \int_{-\d}^\d (\xi+ih)^2 m_1(\xi)
e^{i\tau\xi -\pi\xi-(\2bp)(2k+1)\xi} \, d\xi \\
& \qquad\qquad\qquad\qquad +
\frac{e^{-(\2bp)(2k+1)/2}}{2(\l-e^{-(\2bp)(2k+1)/2})^2}
\int_{-\d}^\d (\xi+ih) m_1(\xi)
e^{i\tau\xi -\pi\xi-(\2bp)(2k+1)\xi} \, d\xi \biggr)  \\
& = \sum_{k\ge1}  
\frac{g_1^{(k)}(\tau)}{\l-e^{-(\2bp)(2k+1)/2}} 
+
\frac{g_2^{(k)}(\tau) e^{-(\2bp)(2k+1)/2}}{(2\l-e^{-(\2bp)(2k+1)/2})^2} \, ,
\end{align*}
where $g_1^{(k)},g_2^{(k)}$ are entire functions of exponential type
with bounds uniform in $k$.

Thus the above series converges to a function $E_4(\tau,\l)$, entire in
$\tau$ and holomorphic when $|\l|>e^{-3\bmp}$.  
\medskip   \\

In order to conclude the proof, we finally need to analyze the sum
of the
terms involving $I$ through $I\!I\!I^*$;
that is
\begin{multline*}
\sum_{j\ge0} \bigl(I-I\!I+I\!I\!I\bigr)\lj +
\sum_{j<0} \bigl(I^*-I\!I^*+I\!I\!I^*\bigr)\lj \\
= 
\sum_{j\ge0} \bigl(I+I\!I\!I\bigr)\lj +
\sum_{j<0} \bigl(I^* +I\!I\!I^*\bigr)\lj 
-\biggl( \sum_{j\ge0} I\!I \lj
+\sum_{j<0} I\!I^*\lj \biggr)  
\ .
\end{multline*}
As in the case of $E^{(1)}$, we split the argument in two.

\subsection{The Case of 
\boldmath $\sum_{j\ge0} \bigl(I+ I\!I\!I\bigr)\lj + 
\sum_{j<0} \bigl(I^*+ I\!I\!I^*\bigr)\lj $}  

Arguing as in (\ref{I+III-first-error-estimate}) and 
(\ref{I+III-second-error-estimate}),
expanding $m$ and $\tilde m$ in Taylor series,
we see that the above sum equals
$$
\sum_{k,\ell=1}^{+\infty} \sum_{j\in\ZZ}
{\cal I}_j^{(k,\ell)} \lj\ ,
$$
where,  for $j\ge0$,
\begin{multline*}
{\cal I}_j^{(k,\ell)} 
= \biggl( \int_{-\infty}^{-\d} 
+\int_{\j2+\d}^{+\infty} 
\biggr)  
\sigma_{(\2bp)(2\ell+1),(2k+1)\pi}(\xi)\\
\times 
(\xi +ih)(\xi +ih -\j2)
e^{i\xi\tau}
e^{-(2k+1)\pi|\xi|}
e^{-(\2bp)|\xi-\j2|} \, d\xi \ ,
\end{multline*}
while for $j<0$
\begin{multline*}
{\cal I}_j^{(k,\ell)} 
= \biggl( \int_{-\infty}^{\j2-\d} 
+\int_{\d}^{+\infty} 
\biggr)  
\sigma_{(\2bp)(2\ell+1),(2k+1)\pi}(\xi)
\\
\times (\xi +ih)(\xi +ih -\j2) 
e^{i\xi\tau}
e^{-(2k+1)\pi|\xi|}
e^{-(\2bp)|\xi-\j2|} \, d\xi \ .
\end{multline*}

By applying Lemmas  \ref{NEW-LEMMA-IANDIII} and
\ref{NEW-LEMMA-I*ANDIII*}
we see that the series above equals
\begin{align*}
& \sum_{k,\ell=1}^{+\infty} 
\frac{e^{ih(R+S)-\d[i\tau+R+S]}}{(i\tau+R+S)^2} 
\biggl[
\frac{e^{R/2}}{2(e^{R/2} - \l)^2}  
   \bigl(1+\psi_1^+(\tau,\l) \bigr) \\    
&\qquad\qquad
+ \frac{e^{-[i\tau +S]/2}}{2(\l - e^{-[i\tau +S]/2})^2}
   \bigl( 1 + \psi_2^-(\tau,\l) \bigr)  
+ \frac{e^{-[i\tau+S]/2}\psi_4^+(\tau)}{
2(e^{R/2} - \l)(e^{-[i\tau + S]/2}-\l)} \biggr] 
\\
& \qquad +  
\frac{e^{-ih(R+S)-\d[R+S-i\tau]}}{(i\tau-R-S)^2} 
\biggl[
\frac{e^{[S- i\tau]/2}}{2(e^{[S - i\tau]/2} - \l)^2}
        \bigl(1 + \psi^+_2(\tau,\l) \bigr) \\
&\qquad\qquad
+ \frac{  e^{-R/2}}{2(\l - e^{-R/2})^2}  
        \bigl(1 +\psi_1^-(\tau,\l) \bigr)  
+ \frac{e^{[S-i\tau]/2}\psi_4^-(\tau)}
{2(e^{-R/2} - \l)(e^{[S-i\tau]/2}-\l)}  \biggr] 
\ ,
\end{align*}
where $R=(2\ell +1)(\2bp)$ and $S=(2k+1)\pi$;
see formula (\ref{final-sum}). 

We are going to show that the functions depending on $k$, that is
on $S$, can be summed, and their sums are functions of $(\tau,\l)$,
holomorphic in  a neighborhood of $\overline{{\cal D}}$, bounded 
(together with their derivatives) as
$|\Re\tau| \rightarrow+\infty$.
\medskip   \\

We let $(\tau,\l)$ vary in the closure of the domain
$$ 
{\cal D}_{4(\b-2\pi),2\pi}=
\bigl\{ \big|\Im\tau -\log|\l|^2 \big|< 2\pi,\
e^{-2(\b-2\pi)}<|\l|<e^{2(\b-2\pi)} \bigr\}\, , 
$$
a domain that contains the closure of the domain
${\mathcal D}$. 
Using 
Lemma \ref{NEW-LEMMA-IANDIII}, we now obtain that
\begin{align} 
\bigg| \frac{e^{ih(R+S)-\d[i\tau+R+S]}}{(i\tau+R+S)^2} 
\bigg|\cdot
\big| e^{R/2}  \bigl(1+\psi_1^+(\tau,\l) \bigr) \big|
& \le  C \bigl( |\Re\tau|+R+S) \frac{e^{-\d[i\tau+R+S]}}{|i\tau+R+S|^2} 
\notag \\    
& \le C e^{-\d (R+S)}\ , \label{est-E3-psi1+} 
\end{align}
as $R,S\rightarrow+\infty$, uniformly in $(\tau,\l)\in{\cal
  D}_{2(\b-2\pi),2\pi}$.

Next notice that the functions
$
\frac{e^{-[i\tau +S]/2}}{2(\l - e^{-[i\tau +S]/2})^2}
$
are bounded for  $(\tau,\l)\in{\cal D}_{4(\b-2\pi),2\pi}$, uniformly in
$S=(2k+1)\pi$ and $R=(2\ell+1)\pi$. 
Then 
\begin{multline}
\bigg|
\frac{e^{ih(R+S)-\d[i\tau+R+S]}}{(i\tau+R+S)^2} 
\cdot
\frac{e^{-[i\tau +S]/2}}{2(\l - e^{-[i\tau +S]/2})^2}
   \bigl( 1 + \psi_2^-(\tau,\l) \bigr)  \bigg| \\
\le C e^{-\d (R+S)} \frac{|\Re\tau|+R+S}{|i\tau+R+S|^2} 
\le C e^{-\d (R+S)} \label{est-E3-psi2-}
\end{multline}
as $R,S\rightarrow+\infty$, uniformly in $(\tau,\l)
\in{\cal D}_{4(\b-2\pi),2\pi}$.

The arguments for the remaining terms are analogous to the ones in
(\ref{est-E1-psi4+})--(\ref{est-E1-psi4-}), so we obtain that  they
are all bounded by a constant times $e^{-\d(R+S)}$.  This completes
the proof of this case.
\medskip    \\

Finally, we have
\subsection{The Case of 
\boldmath $\sum_{j\ge0} I\!I\lj + 
\sum_{j<0} I\!I^*\lj$}  
As in the previous Subsection \ref{summation-by-parts}
we use 
Lemmas \ref{NEW-LEMMA-II} and \ref{NEW-LEMMA-II*},
and  summation by parts.  We begin with the sum for $j\ge0$.
Since $\d\le\xi\le\j2-\d$, we can write
\begin{align*}
& \sum_{j\ge0} I\!I \lj \\
& = e^{ih(2\b-2\pi)} \sum_{j\ge0} e^{-(2\b-\pi)\j2}\lj
\int_{\d}^{\j2-\d} \bigl[ \xi^2 +b\xi +c\bigr] 
e^{(i\tau+(\2bp)-\pi)\xi} m_1(\xi)\tilde m(\xi) \, d\xi\\
& =  e^{ih(2\b-2\pi)} \sum_{j\ge0} e^{-(2\b-\pi)\j2}
\lj \biggl( \sum_{k=1}^N
e^{-2\pi i kh} 
\int_{\d}^{\j2-\d} \bigl[ \xi^2 +b\xi +c\bigr] 
e^{(i\tau+(\2bp)-(2k+1)\pi)\xi} \, d\xi\\
& \qquad\qquad + \int_{\d}^{\j2-\d} \bigl[ \xi^2 +b\xi +c\bigr] 
e^{(i\tau+(\2bp)-(2N+1)\pi)\xi} m_1(\xi)\, d\xi \biggr) \\
& \equiv A+B\ .
\end{align*}

Now we only sketch the details.  In order to evaluate $A$ we apply Lemma
\ref{NEW-LEMMA-II} with $R=2\b-\pi$ and $S=(2k+1)\pi$, $k=1,\dots,N$.
In order to evaluate $B$ we use summation by parts (see
\ref{summation-by-parts-formula}) and argue as in Subsection
\ref{summation-by-parts} 
This concludes the proof of  Proposition
\ref{sum-of-error-terms-proposition}.  \bigskip
\endpf

\section{Completion of the Proof of Proposition 
\ref{technical}.\label{appendix}}
\vspace*{.12in}

In this final section we complete
the proof Proposition \ref{technical}.

\proof[End of the proof of Proposition \ref{technical}]
We are in a position now to conclude the proof of Proposition 
\ref{technical}, modulo proving the lemmas from Section
\ref{new-section}. 
We shall also prove those lemmas at this time (after the
proof of the proposition).
Using Lemmas \ref{NEW-LEMMA-IANDIII} through
\ref{NEW-LEMMA-II*} we obtain that
\begin{align}
& \sum_{j \geq 0} (I - I\!I +I\!I\!I) \lj 
+ \sum_{j\le-1} (I^* - I\!I^* + I\!I\!I^*)\lj
\notag\\ 
& = 
\frac{e^{ih(R+S)+R/2}}{2(e^{R/2} - \l)^2(i\tau+R+S)^2} 
e^{-\d[i\tau+R+S]}
\bigl(1+\psi_1^+(\tau,\l) \bigr)     \notag\\
& \qquad + 
\frac{e^{-ih(R+S)+[S- i\tau]/2}}{2(e^{[S - i\tau]/2} - \l)^2(i\tau-R-S)^2}
        e^{-\d[R+S-i\tau]} 
\bigl( 1 + \psi_2^+(\tau,\l) \bigr)  \notag\\
& \qquad +
\frac{e^{ih(R-S)-\d[S-R-i\tau]+ R/2}}{2(e^{[S-i\tau]/2}-\l)^2(e^{R/2}-\l)^2}
\psi_3^+(\tau,\l)
\notag\\
\end{align}

\begin{align}
& \qquad + 
\biggl( \frac{2e^{ih(R+S)-\d[i\tau+R+S]}}{(e^{R/2} - \l)(i\tau+R+S)^3}
     - \frac{2e^{-ih(R+S)-\d[R+S-i\tau]}}{(e^{[S - i\tau]/2} -
       \l)(i\tau-R-S)^3} \biggr) 
\notag\\ 
& \qquad +
\frac{  e^{-ih(R+S)-R/2}}{2(\l - e^{-R/2})^2(i\tau-R-S)^2} 
e^{-\d[R+S-i\tau]}
 \bigl(1 +
\psi_1^-(\tau,\l) \bigr)     \notag\\
& \qquad + \frac{e^{ih(R+S)}
e^{-[i\tau +S]/2}}{2(\l - e^{-[i\tau +S]/2})^2(i\tau+R+S)^2}
   e^{-\d[i\tau+R+S]} \bigl( 1 + \psi_2^-(\tau,\l) \bigr)  \notag\\
& \qquad +
\frac{e^{-ih(R-S)+\d[i\tau-R+S]- R/2}}{2(e^{-[i\tau+S]/2}-\l)^2(e^{-R/2}-\l)^2} 
\psi_3^-(\tau,\l) \notag \\
& \qquad+  \biggl( \frac{2e^{-ih(R+S)-\d[R+S-i\tau]}}{(e^{-R/2} - \l)(i\tau-R-S)^3}
 - \frac{2e^{ih(R+S)-\d[i\tau+R+S]}}{(e^{-[i\tau + S]/2} -
   \l)(i\tau+R+S)^3} \biggr) 
       \, . \label{sum+and-}
\end{align}

Previously we have alluded to a certain cancellation of
the cubic terms.  We expect this because the worm domain
is locally like a product 
and the kernels for such domains have degree-two singularities (see
[APF]). 
Thus no terms with third-order singularity should be present.
Now we concentrate on the cubic terms and
calculate to confirm this expectation.

Notice that 
\begin{align}
&\frac{1}{(e^{-R/2} - \l)(i\tau-R-S)^3} -
       \frac{1}{(e^{[S - i\tau]/2} - \l)(i\tau - R-S)^3}\notag \\
& \qquad\qquad \qquad\qquad \qquad\qquad
=  \frac{1}{(i\tau-R-S)^3}  
\frac{e^{[S - i\tau]/2}-e^{-R/2}}{(e^{-R/2} - \l)(e^{[S-i\tau]/2}-\l)}
       \notag \\
& \qquad\qquad \qquad\qquad \qquad\qquad
=  \frac{e^{[S-i\tau]/2}\psi_4^-(\tau)}
{2(i\tau-R-S)^2(e^{-R/2} - \l)(e^{[S-i\tau]/2}-\l)} \,
       .\label{first-difference-cubes}
\end{align}  

Analogously, 
\begin{multline}
\frac{1}{(e^{R/2} - \l)(i\tau+R+S)^3}
- \frac{1}{(e^{-[i\tau + S]/2} - \l)(i\tau+R+S)^3} \\
= \frac{e^{-[i\tau+S]/2}\psi^+_4(\tau)}{
2(i\tau+R+S)^2(e^{R/2} - \l)(e^{-[i\tau + S]/2}-\l)} \, , 
\label{2nd-difference-cubes}
\end{multline}
where, in the last two displays, 
\begin{equation}\label{psi4+-}
 \psi_4^\pm  (\tau) 
 =  E(i\tau\pm(R+S))
\ ,\qquad\text{with}
\qquad
 E(x)
=  \frac{1-e^{x/2}}{x/2}  
\ \, . 
\end{equation}
\vspace*{.12in}

Now we plug formulas (\ref{first-difference-cubes}),
(\ref{2nd-difference-cubes}) 
into (\ref{sum+and-}) and factor common powers once again. 
In total, we finally obtain
that the sum on the lefthand side of (\ref{sum+and-}) equals
\begin{align}
& \frac{e^{ih(R+S)-\d[i\tau+R+S]}}{(i\tau+R+S)^2} 
\biggl[
\frac{e^{R/2}}{2(e^{R/2} - \l)^2}  
   \bigl(1+\psi_1^+(\tau,\l) \bigr) \notag \\    
&\qquad\qquad
+ \frac{e^{-[i\tau +S]/2}}{2(\l - e^{-[i\tau +S]/2})^2}
   \bigl( 1 + \psi_2^-(\tau,\l) \bigr)  
+ \frac{e^{-[i\tau+S]/2}\psi_4^+(\tau)}{
2(e^{R/2} - \l)(e^{-[i\tau + S]/2}-\l)} \biggr] 
\notag \\
& \qquad +  
\frac{e^{-ih(R+S)-\d[R+S-i\tau]}}{(i\tau-R-S)^2} 
\biggl[
\frac{e^{[S- i\tau]/2}}{2(e^{[S - i\tau]/2} - \l)^2}
        \bigl(1 + \psi^+_2(\tau,\l) \bigr) \notag \\
&\qquad\qquad
+ \frac{  e^{-R/2}}{2(\l - e^{-R/2})^2}  
        \bigl(1 +\psi_1^-(\tau,\l) \bigr)  
+ \frac{e^{[S-i\tau]/2}\psi_4^-(\tau)}
{2(e^{-R/2} - \l)(e^{[S-i\tau]/2}-\l)}  \biggr] \notag \\
&\qquad
+
\frac{e^{ih(R-S)-\d[S-R-i\tau]+R/2}}{2(e^{[S-i\tau]/2}-\l)^2(e^{R/2}-\l)^2}
\psi_3^+(\tau,\l) \notag \\
&\qquad 
+
\frac{e^{-ih(R-S)+\d[i\tau-R+S]- R/2}}{2(e^{-[i\tau+S]/2}-\l)^2(e^{-R/2}-\l)^2} 
\psi_3^-(\tau,\l) 
\ . \label{final-sum}
\end{align}

If we set $\d=0$ now, we obtain the expression of the sum as in the
statement of Proposition \ref{technical}.

Finally, notice that
from (\ref{psi4+-}) it follows at once that
\begin{equation}\label{estimate-for-psi+-}
|\psi_4^\pm(\tau,\l)|
\le C_M e^{\pm(R+S)/2} \, .
\end{equation}
This 
proves Proposition \ref{technical}.
\endpf	\\
\vspace*{.12in}

Finally, we need to prove Lemmas
\ref{NEW-LEMMA-IANDIII}--\ref{NEW-LEMMA-II*}. 
\smallskip \\
	       
\proof[Proof of Lemma \ref{NEW-LEMMA-IANDIII}]
We have
$$
I = e^{ih(R+S)} 
e^{-R\j2} e^{-\d[i\tau+R+S]}   
\biggl\{ 
 \frac{c-\d b+\d^2}{i\tau+R+S} - \frac{b-2\d}{(i\tau+R+S)^2}   
     + \frac{2}{(i\tau+R+S)^3} \biggr\} 
$$
and
\begin{align*}
I\!I\!I   
& =  - e^{-ih(R+S)-\d[R+S-i\tau]}e^{\j2[i\tau-S]}
\biggl\{ \biggl( \j2+\d \biggr)^2 \frac{1}{i\tau-R-S} \\
        & \qquad\qquad 
+ \biggl(\j2+\d\biggr)
\biggl( \frac{b}{i\tau-R-S} - \frac{2}{(i\tau-R-S)^2} \biggr) \\
& \qquad\qquad 
+   \biggl( \frac{c}{i\tau-R-S} 
   -\frac{b}{(i\tau-R-S)^2} + \frac{2}{(i\tau-R-S)^3} \biggr) 
\biggr\} \, .
\end{align*}
Therefore, recalling that $b=2ih-(j+1)/2$ and $c=-h^2-ih(j+1)/2$, we
see that 
\begin{align*}
& I  +I\!I\!I \\   
& =e^{ih(R+S)-\d[i\tau+R+S]}
e^{-R\j2} \biggl\{ \j2 \biggl( \frac{1}{(i\tau+R+S)^2} +
      \frac{\d-ih}{i\tau+R+S} \biggr)  \\
&\qquad\qquad  +  \biggl( \frac{2}{(i\tau+R+S)^3}
 +\frac{2(\d-ih)}{(i\tau+R+S)^2} +\frac{(\d-ih)^2}{i\tau+R+S} \biggr) 
\biggr\}  \\
& \quad -  
e^{\j2[i\tau -S]} \biggl\{ 
 e^{-ih(R+S)-\d[R+S-i\tau]} \j2 \biggl(
 \frac{\d+ih}{i\tau - R-S}  
      - \frac{1}{(i\tau -R-S)^2} \biggr)  \\
&\qquad\qquad
+ e^{-ih(R+S)-\d[R+S-i\tau]} \biggl( 
\frac{(\d+ih)^2}{i\tau - R-S} - \frac{2(\d+ih)}{(i\tau -R-S)^2}
+ \frac{2}{(i\tau -R-S)^3} \biggr)
            \biggr\}  
\ .
\end{align*}

Applying the summation formulas 
(\ref{summation-formulas})
to our last expression
we obtain for
$$
\big|e^{-R/2}\l\big|<1  \qquad\text{and}\quad
\big|e^{[i\tau-S]/2}\l\big|<1\ ,
$$
that is for $(\tau,\l)\in {\cal D}_{R,S}$, that
\begin{align*}
& \sum_{j \geq 0} (I +I\!I\!I) \l^j   \\
& =  \frac{e^{R/2}}{2(e^{R/2} - \l)^2}
     \biggl[ e^{ih(R+S)-\d[i\tau+R+S]}\biggl( 
\frac{1}{(i\tau+R+S)^2} + \frac{\d-ih}{i\tau+R+S} \biggr)
\biggr] \\
&\quad +\frac{1}{e^{R/2} - \l}  \biggl[  e^{ih(R+S)-\d[i\tau+R+S]}\biggl(
    \frac{2}{(i\tau+R+S)^3}  +\frac{2(\d-ih)}{(i\tau+R+S)^2} 
+\frac{(\d-ih)^2}{i\tau+R+S} \biggr) 
              \biggr]  \\
&\quad +\frac{e^{[S-i\tau]/2}}{2(e^{[S- i\tau]/2} - \l)^2}
    \biggl[ e^{-ih(R+S)-\d[R+S-i\tau]}\biggl(
\frac{1}{(i\tau-R-S)^2} -\frac{\d+ih}{i\tau-R-S} \biggr) 
 \biggr] \\
&\quad -\frac{1}{e^{[S - i\tau]/2} - \l} 
\biggl[ 
e^{-ih(R+S)-\d[R+S-i\tau]} \biggl( \frac{2}{(i\tau - R-S)^3}    
-\frac{2(\d+ih)}{(i\tau -R-S)^2} 
+\frac{(\d+ih)^2}{i\tau -R-S} \biggr)
           \biggr] \, . 
\end{align*}

We have
\begin{align*} 
& \sum_{j \geq 0} (I +I\!I\!I) \l^j \\
& \quad = \frac{e^{ih(R+S)-\d[i\tau+R+S]+ R/2}}{2(e^{R/2} - \l)^2} 
   \biggl[ \frac{1}{(i\tau+R+S)^2} +  \frac{\d-ih}{i\tau+R+S}
   \biggr]  \\  
&\qquad+ \frac{e^{ih(R+S)-\d[i\tau+R+S]}}{e^{R/2} - \l}  \biggl[
   \frac{2}{(i\tau+R+S)^3} + \frac{2(\d-ih)}{(i\tau+R+S)^2} 
+\frac{(\d-ih)^2}{i\tau+R+S} \biggr] \\
&\qquad+ \frac{e^{-ih(R+S)-\d[R+S-i\tau]+[S - i\tau]/2}}{2(e^{[S - i\tau]/2}
  - \l)^2}  
    \biggl[ \frac{1}{(i\tau-R-S)^2} -  \frac{\d+ih}{i\tau-R-S}
    \biggr]     \\
&\qquad
+ \frac{e^{-ih(R+S)-\d[R+S-i\tau]}}{e^{[S - i\tau]/2} - \l} \biggl[ 
   -\frac{2}{(i\tau -R-S)^3} + \frac{2(\d+ih)}{(i\tau -R-S)^2}
       - \frac{(\d+ih)^2}{i\tau -R-S}  \biggr]  
\, . 
\end{align*}
(Later on we shall see that the terms here that are of
cubic order will cancel with the cubic terms that
come from summing $j$ with $j < 0$.)

Factoring out common powers, we may simplify our expression
further and finally  write
\begin{align}  
& \sum_{j \geq 0} (I +I\!I\!I) \l^j \notag \\
& =  \frac{e^{ih(R+S)+R/2}}{2(e^{R/2} - \l)^2(i\tau+R+S)^2}
e^{-\d[i\tau+R+S]} \bigl(1+\psi_1^+(\tau,\l) \bigr)     \notag\\
& \quad + 
\frac{e^{-ih(R+S)+[S - i\tau]/2}}{2(e^{[S - i\tau]/2} - \l)^2(i\tau-R-S)^2}
        e^{-\d[R+S-i\tau]} \bigl( 1 + \psi_2^+(\tau,\l) \bigr)  \notag\\
& \quad +
\biggl( \frac{2e^{ih(R+S)-\d[i\tau+R+S]}}{(e^{R/2} - \l)(i\tau+R+S)^3}
     - \frac{2e^{-ih(R+S)-\d[R+S-i\tau]}}{(e^{[S - i\tau]/2} -
       \l)(i\tau-R-S)^3} \biggr)  
 \, ,\label{sum+}
\end{align}
where 
\smallskip
\begin{equation}\label{psi1+psi2+}
\begin{array}{ll}
 \psi_1^+ (\tau,\l)
 & =  {\displaystyle 
(\d-ih)(i\tau+R+S)  }\medskip \\
& \qquad  
{\displaystyle +2(1-e^{-R/2}\l) \bigl[ 
2(\d-ih) +(\d-ih)^2(i\tau+R+S) \bigr]}\ ,\medskip \\
 \psi_2^+ (\tau,\l)
& =  {\displaystyle 
-(\d+ih)(i\tau-R-S) }\medskip \\
& \qquad  
{\displaystyle 
+2 \bigl( 1-e^{-[S-i\tau]/2}\l\bigr)
\bigl[2(\d+ih)-(\d+ih)^2(i\tau-R-S)\bigr] } 
\ . 
\end{array} \medskip
\end{equation}
Finally, the estimate for $\psi_1^+$ and $\psi_2^+$   
 follows at once from the explicit expressions 
(\ref{psi1+psi2+}) \, .
\endpf 
\medskip \\

In order to calculate the sum for $I\!I$ (Lemma 9.3), we need 
 the following elementary but somewhat tedious
lemma, whose proof is, once again, postponed.
\begin{lemma}\label{lemma-Q}   \sl
Let $\a>0$ and $\d\ge0$ be fixed. We let
$f(x)=(\a e^{-x/2}-\l)^{-1}$,  $\Psi_\pm(x)=e^{\mp 2\d x}f(x)$, and define
\begin{multline}
Q_{\pm\d} (x)=
\frac{ihx-1}{x^2} \biggl[ 2\frac{\Psi_\pm(x)-\Psi_\pm(0)}{x}
-\bigl(\Psi_\pm'(x)+\Psi_\pm'(0)\bigr) 
\biggr]\\
 + \frac1x \biggl[
(h^2-\d^2) \bigl( \Psi_\pm(x)-\Psi_\pm(0) \bigr) \mp\d\bigl(e^{\mp2\d
  x}f'(x)-f'(0)\bigr) \biggr] \ . \label{----}
\end{multline}
Then there exist entire functions $G^\pm_1,\dots,G^\pm_6$, independent
of $\a$ and $\l$,  such that
\begin{multline}
Q_{\pm\d} (x)=
\frac{\a}{2(\a e^{-x/2}-\l)^2(\a-\l)^2} \biggl[
(\a e^{-x/2}-\l)(\a-\l)G^\pm_1(x) +(\a-\l)^2 G^\pm_2(x) \\
+\frac{(\a-\l)(\a e^{-x/2}-\l)^2}{\a} G^\pm_3 (x) +\a\l G^\pm_4(x)
+\l G^\pm_5(x) +\a^2 G^\pm_6(x) \biggr]
\ . \label{++++}
\end{multline}
\end{lemma}

\noindent The proof of Lemma \ref{lemma-Q} occurs at the very end of
this paper. 
\smallskip \\

\proof[Proof of Lemma \ref{NEW-LEMMA-II}]
Using (\ref{calculus}) we have that
\begin{align*}
I\!I & = e^{ih(R-S)}
e^{-R\j2} 
   \biggl\{ \biggl[ 
     \biggl( \j2 -\d\biggr)^2  \frac{1}{i\tau+R-S} \\
& \qquad 
      + \biggl( \j2 -\d\biggr)
\biggl( \frac{b}{i\tau+R-S} - \frac{2}{(i\tau+R-S)^2} \biggr)  \\
& \qquad  
+ \biggl( \frac{c}{i\tau+R-S} - \frac{b}{(i\tau+R-S)^2} 
+    \frac{2}{(i\tau+R-S)^3} 
            \biggr) 
  \biggr]  e^{(\j2-\d)[i\tau+R-S]} \\
& \qquad
 - \biggl[ \frac{c+\d b+\d^2}{i\tau+R-S}  - \frac{b}{(i\tau+R-S)^2} 
+ \frac{2}{(i\tau+R-S)^3} \biggr] e^{\d[i\tau+R-S]}  \biggr\} \, .
\end{align*}
Therefore, recalling that $b=2ih-(j+1)/2$ and $c=-h^2-ih(j+1)/2$, we
see that 
\begin{align*}
I\!I 
& = 
e^{\j2[i\tau -S]} \biggl\{ e^{ih(R-S)-\d[i\tau+R-S]} \j2 
\biggl( \frac{ih-\d}{i\tau+R-S} -\frac{1}{(i\tau+R-S)^2}\biggr) \\
&\qquad\qquad  
+ e^{ih(R-S)-\d[i\tau+R-S]} \biggl( 
 \frac{(\d-ih)^2}{i\tau+R-S}   + \frac{2(\d-ih)}{(i\tau+R-S)^2}  
+ \frac{2}{(i\tau+R-S)^3} \biggr) \\
&  \quad- e^{ih(R-S)+\d[i\tau+R-S]}
e^{-R\j2} \biggl\{ \j2 \biggl( \frac{1}{(i\tau+R-S)^2}
              - \frac{\d+ih}{i\tau+R-S}  \biggr)   \\
&\qquad\qquad 
+  \biggl( \frac{2}{(i\tau+R-S)^3} - \frac{2(\d+ih)}{(i\tau+R-S)^2}
           +\frac{(\d+ih)^2}{i\tau+R-S} \biggr) \biggr\} \, .
\end{align*}

Applying the summation formulas 
(\ref{summation-formulas})
we obtain that,
for
$$
\big|e^{-R/2}\l\big|<1  \qquad\text{and}\quad
\big|e^{[i\tau-S]/2}\l\big|<1\ ,
$$
that is for $(\tau,\l)\in {\cal D}_{R,S}$,
\begin{align}
& \sum_{j \geq 0} I\!I \l^j   \notag\\
& =  \frac{e^{R/2}}{2(e^{R/2} - \l)^2}
      e^{ih(R-S)+\d[i\tau+R-S]}
\biggl( -\frac{\d+ih}{i\tau+R-S} 
+ \frac{1}{(i\tau+R-S)^2}  \biggr)  \biggr] \notag\\
&\quad +\frac{1}{e^{R/2} - \l}  \biggl[  
e^{ih(R-S)+\d[i\tau+R-S]}\biggl(
\frac{2}{(i\tau+R-S)^3} - \frac{2(\d+ih)}{(i\tau+R-S)^2} 
+\frac{(\d+ih)^2}{i\tau+R-S} \biggr)
              \biggr]  \notag\\
&\quad +\frac{e^{[S-i\tau]/2}}{2(e^{[S- i\tau]/2} - \l)^2}
    \biggl[ e^{ih(R-S)-\d[i\tau+R-S]}\biggl(
\frac{1}{(i\tau+R-S)^2}  +\frac{(\d+ih)}{i\tau+R-S} 
  \biggr) \biggr] \notag\\
&\quad -\frac{1}{e^{[S - i\tau]/2} - \l} 
\biggl[ e^{ih(R-S)-\d[i\tau+R-S]}\biggl( \frac{2}{(i\tau+R-S)^3}
  +\frac{2(\d-ih)}{(i\tau+R-S)^2} 
+\frac{(\d-ih)^2}{i\tau+R-S} \biggr)  
\biggr] \, . \label{with-neg-powers-of-x}
\end{align}

Notice that there  is an apparent singularity as 
$i\tau \ra -(R-S)$.  We now analyze the terms 
containing negative powers of $(i\tau+R-S)$ and show that they
actually give rise to some cancellation. 
But we shall also see that there genuinely
{\it are} singularities as $i\tau \ra -(R+S)$ and $i\tau \ra (R + S)$.

Temporarily, we introduce the notation
$$
x = i\tau+R-S,  \qquad\text{and}\qquad
\a =  e^{R/2}\, 
$$
so that
$$
S-i\tau  = R-x
\quad\text{and}\quad
 e^{[S-i\tau]/2} = \a e^{-x/2} \, . 
$$
We now 
show
that the ostensible singularities in $x=0$ in fact all cancel out.
Then we may rewrite the 
terms containing negative powers of $(i\tau+R-S)$ (that is, of $x$)
 on the righthand side of
formula (\ref{with-neg-powers-of-x}) as
\begin{align}
& e^{ih(R-S)}\biggl\{ 
\frac{1}{x^3}  \biggl[ 
   \frac{2e^{\d x}}{\a - \l} - 
\frac{2e^{-\d x}}{\a e^{-x/2} - \l} \biggr]
   \notag \\
& \qquad+  \frac{1}{x^2}  \biggl[ 
\frac{\a e^{-x/2}e^{-\d x}}{2(\a e^{-x/2} - \l)^2}
    + \frac{2(ih-\d)e^{-\d x}}{\a e^{-x/2} - \l} 
+ \frac{\a e^{\d x}}{2(\a - \l)^2}
- \frac{2(\d+ih)e^{\d x}}{\a - \l}  \biggr]  
   \notag \\
& \qquad+ \frac{1}{x}  \biggl[ 
\frac{(\d-ih)\a e^{-x/2}e^{-\d x}}{2(\a 
     e^{-x/2} - \l)^2} 
-\frac{(\d-ih)^2 e^{-\d x}}{\a e^{-x/2} - \l}
- \frac{(\d+ih)\a e^{\d x}}{2(\a - \l)^2}
     +\frac{(\d+ih)^2 e^{\d x}}{\a - \l} 
    \biggr] \biggr\} \, . \label{RHS-of-neg-powers} 
\end{align} 

Now if  we set
$$
f(x) = \frac{1}{\a \ex2 - \l}
\ , \qquad\text{and}\qquad
F(x) =e^{-2\d x} f(x) \, ,
$$
then we see that we may rewrite the above expression
(\ref{RHS-of-neg-powers}) as $e^{ih(R-S)+\d x} Q_\d(x)$, where 
$Q_\d$ is given by (\ref{----}). 

Next we use Lemma \ref{lemma-Q} 
to see that $Q_\d$ is regular in $x=0$ and 
to write 
\begin{multline*}
Q_\d (x)=
\frac{\a}{2(\a e^{-x/2}-\l)^2(\a-\l)^2} \biggl[
(\a e^{-x/2}-\l)(\a-\l)G^+_1(x) +(\a-\l)^2 G^+_2(x) \\
+\frac{(\a-\l)(\a e^{-x/2}-\l)^2}{\a} G^+_3 (x) +\a\l G^+_4(x)
+\l G^+_5(x) +\a^2 G^+_6(x) \biggr]	 \, ,
\end{multline*}
where
$G^+_1,\dots, G^+_6$ are the functions appearing in Lemma \ref{lemma-Q} and
whose explicit expression is 
given in 
(\ref{G1G6}).

Therefore
the
above expression (\ref{RHS-of-neg-powers}) equals
\begin{equation*}
e^{ih(R-S)+\d[i\tau+R-S]} Q_\d \bigl(i\tau+R-S\bigr) 
 = \frac{e^{ih(R-S)+\d[i\tau+R-S] +
     R/2}}{2(e^{[S-i\tau]/2}-\l)^2(e^{R/2}-\l)^2} \psi_3^+(\tau,\l)
\ ,
\end{equation*}
where  we have set 
\begin{equation}\label{varphi1+-varphi2+}
 \varphi^+_j  (\tau) = G_j^+(i\tau+R-S)\ ,\qquad\text{for}\quad
j=1,\dots,6
\  ,
\end{equation}
and 
\begin{multline} \label{new-psi3+}
\psi_3^+(\tau,\l) =
(e^{[S-i\tau]/2}-\l)(e^{R/2}-\l)\varphi^+_1(\tau) +(e^{R/2}-\l)^2
\varphi^+_2(\tau) \\ 
+(1-e^{-R/2}\l)(e^{[S-i\tau]/2}-\l)^2 \varphi^+_3 (\tau) +e^{R/2}\l
\varphi^+_4(\tau) +\l \varphi^+_5(\tau) +e^R \varphi^+_6(\tau) 
\ .
\end{multline}

In view of the cancellation of negative powers of $x$---in
other words, any dependency on $x$ is in fact bounded and smooth---we 
may simplify formula (\ref{with-neg-powers-of-x})
to obtain
\begin{equation}\label{to-be-named}
\sum_{j \geq 0} I\!I \l^j =
\frac{e^{ih(R-S)-\d[S-R-i\tau]+ R/2}}{2(e^{[S-i\tau]/2}-\l)^2(e^{R/2}-\l)^2}
\psi_3^+(\tau,\l)\ .
\end{equation}

Recalling definition
(\ref{new-psi3+})), we write
$\psi_3^+$ explicitly as
\begin{align*}
\psi_3^+ (\tau,\l)
& =  (e^{[S-i\tau]/2}-\l)(e^{R/2}-\l)G^+_1(i\tau+R-S)  
+(e^{R/2}-\l)^2
G^+_2(i\tau+R-S) \\ 
& \qquad 
+(1-e^{-R/2}\l)^2(e^{[S-i\tau]/2}-\l) G^+_3 (i\tau+R-S) +e^{R/2}\l
G^+_4(i\tau+R-S) \\
& \qquad 
+\l(\a-\l) G^+_5(i\tau+R-S) +e^R G^+_6(i\tau+R-S) 
\, ,
\end{align*}
where $G^+_1,\dots,G_6^+$ are given in (\ref{G1G6}). Then, since
these  are entire functions, we have
\begin{align*}
|\psi_3^+ (\tau,\l)|
& \le C_M \biggl(
e^{(R+S)/2} |G^+_1(i\tau+R-S)| + e^{R/2} |G^+_2(i\tau+R-S)| \\
& \qquad 
+e^{S/2} |G^+_3(i\tau+R-S)| 
+e^{R/2} \bigl(|G^+_4(i\tau+R-S)| + |G^+_5(i\tau+R-S)| \bigr) 
\\
& \qquad 
+ e^R|G^+_6(i\tau+R-S)| \biggr) \\
& \le C_M
\bigl( e^{R+S/2}   +e^{S(2\d+1/2)} 
+ e^{2R} \bigr) \ .
\end{align*}
This proves Lemma \ref{NEW-LEMMA-II}.  
\endpf

\proof[Proof of Lemma \ref{NEW-LEMMA-I*ANDIII*}]
Recall formula (\ref{calculus})
and that $b=2ih-(j+1)/2$ and $c=-h^2-ih(j+1)/2$.  
We calculate:
\begin{align*}
I^* & =  e^{ih(R+S)}e^{-R\j2} \biggl\{\biggl[ \frac{(\j2-\d)^2}{i\tau+R+S} +
\biggl( \frac{2ih - \j2}{i\tau+R+S} - \frac{2}{(i\tau+R+S)^2} \biggr) 
\biggl(\j2-\d\biggr)  \\
& \qquad+  \biggl( \frac{2}{(i\tau+R+S)^3} - \frac{2ih -\j2}{(i\tau+R+S)^2}
        + \frac{ih(ih - \j2)}{i\tau+R+S} \biggr) \biggr] e^{(\j2-\d)(i\tau+R+S)} 
\biggr\} \\ 
& =  e^{ih(R+S)-\d[i\tau+R+S]}e^{[i\tau+S]\j2} \biggl\{ \frac{2}{(i\tau+R+S)^3}  
+\frac{2(\d-ih)}{(i\tau +R+S)^2} +\frac{(\d-ih)^2}{(i\tau +R+S)} \\
& \qquad 
- \j2 \biggl(  \frac{1}{(i\tau+R+S)^2} +\frac{\d-ih}{i\tau+R+S} \biggr)
\biggr\}
 \, .
\end{align*}

Next, for the term $I\!I\!I^*$, we have
\begin{align*}
I\!I\!I^* 
& = - e^{-ih(R+S)-\d[R+S-i\tau]} 
e^{R\j2} \biggl\{ \biggl(
  \frac{2}{(i\tau - R-S)^3} - \frac{2ih - \j2 +2\d}{(i\tau-R-S)^2}\\
& \qquad\qquad\qquad  
 +  \frac{ih(ih - \j2)+(2ih-\j2)\d+\d^2}{i\tau-R-S} \biggr) \biggr\}     \\
& =  - e^{-ih(R+S)-\d[R+S-i\tau]} e^{R\j2} \biggl\{
 \frac{2}{(i\tau -R-S)^3} -\frac{2(\d+ih)}{(i\tau-R-S)^2} 
+\frac{(\d+ih)^2}{i\tau-R-S}\\
& \qquad\qquad\qquad  
+ \j2 \biggl( \frac{1}{(i\tau-R-S)^2}
    - \frac{\d+ih}{i\tau-R-S} \biggr)  \biggr\} \, .
\end{align*}

Now we must sum in $j$.  Since
$$
\sum_{j\le-1} e^{\gamma \j2} \lj = \frac{1}{\l-e^{-\gamma/2}}
\qquad\text{and}
\qquad
\sum_{j\le-1} \j2 e^{\gamma\j2}  \lj = -
   \frac{e^{-\gamma/2}}{2(\l- e^{-\gamma/2})^2} \, ,
$$
when $\big|e^{\gamma/2}\l\big|>1$, i.e. when $|\l|>e^{-R/2}$, 
we have that 
\begin{align*}  
& \sum_{j\le-1} (I^* + I\!I\!I^*)\lj\notag \\
& = \frac{1}{\l - e^{-R/2}}
 \biggl[  e^{-ih(R+S)-\d[R+S-i\tau]} 
\biggl( -  \frac{2}{(i\tau-R-S)^3} 
+ \frac{2(\d+ih)}{(i\tau-R-S)^2}  - \frac{(\d+ih)^2}{i\tau-R-S} \biggr) \biggr]
      \notag \\ 
& \qquad 
+  \frac{e^{-R/2}}{2(\l - e^{-R/2})^2}
  \biggl[ 
e^{-ih(R+S)-\d[R+S-i\tau]}\biggl( \frac{1}{(i\tau-R-S)^2}
+ \frac{\d+ih}{i\tau-R-S} \biggr) \biggr] \notag\\
& \qquad +  \frac{1}{\l - e^{-[i\tau +S]/2}}  \biggl[
e^{ih(R+S)-\d[i\tau+R+S]}\biggl(     
\frac{2}{(i\tau+R+S)^3} 
+\frac{2(\d-ih)}{(i\tau+R+S)^2}
+\frac{(\d-ih)^2}{i\tau+R+S} \biggr) \biggr] \notag \\
&  \qquad +  \frac{e^{-[i\tau +S]/2}}{2(\l - e^{-[i\tau +S]/2})^2}
    \biggl[ 
e^{ih(R+S)-\d[i\tau+R+S]}\biggl( \frac{1}{(i\tau+R+S)^2} 
+\frac{\d-ih}{i\tau+R+S} \biggr)
\biggr] \biggr\}
\, . 
\end{align*}

Therefore,
factoring common denominators, we
we may simplify the above expression and
obtain that
\begin{align}
& \sum_{j\le-1} (I^*+ I\!I\!I^*)\lj\notag \\
& = e^{-ih(R+S)-\d[R+S-i\tau]} \biggl\{
\frac{1}{\l - e^{-R/2}}  \biggl[
  - \frac{2}{(i\tau-R-S)^3} + \frac{2(\d+ih)}{(i\tau-R-S)^2} 
- \frac{(\d+ih)^2}{i\tau-R-S}
     \biggr] \notag\\
& \quad +  \frac{e^{-R/2}}{2(\l - e^{-R/2})^2}
   \biggl[ \frac{1}{(i\tau-R-S)^2} - \frac{\d+ih}{i\tau-R-S} 
     \biggr]  \biggr\}\notag\\
&\quad +  e^{ih(R+S)-\d[i\tau+R+S]} \biggl\{
\frac{1}{\l - e^{-[i\tau + S]/2}} \biggl[
  \frac{2}{(i\tau+R+S)^3} +\frac{2(\d-ih)}{(i\tau+R+S)^2} 
+\frac{(\d-ih)^2}{i\tau+R+S}
      \biggr] \notag\\
&\quad  +  \frac{e^{-[i\tau + S]/2}}{2(\l - e^{-[i\tau +S]/2})^2}
    \biggl[ \frac{1}{(i\tau+R+S)^2} +\frac{\d-ih}{i\tau+R+S} 
        \biggr] \biggr\} \notag\\
& = 
\frac{  e^{-ih(R+S)-R/2}}{2(\l - e^{-R/2})^2(i\tau-R-S)^2}  
e^{-\d[R+S-i\tau]} \bigl(1 +
\psi_1^-(\tau,\l) \bigr)     \notag\\
& \quad + \frac{e^{ih(R+S)}
e^{-[i\tau +S]/2}}{2(\l - e^{-[i\tau + S]/2})^2(i\tau+R+S)^2}
   e^{-\d[i\tau+R+S]} \bigl( 1 + \psi_2^-(\tau,\l) \bigr)  \notag\\
& \quad+  \biggl( \frac{2e^{-ih(R+S)-\d[R+S-i\tau]}}{(e^{-R/2} - \l)(i\tau-R-S)^3}
 - \frac{2e^{ih(R+S)-\d[i\tau+R+S]}}{(e^{-[i\tau + S]/2} - \l)(i\tau+R+S)^3} \biggr)
       \label{sum-} \, , 
\end{align} 
where 
\begin{equation}\label{psi1-psi2-}
\begin{array}{ll}
 \psi_1^- (\tau,\l) 
& =  {\displaystyle 
-(\d+ih)(i\tau-R-S)  }\medskip \\
& \qquad  
{\displaystyle 
+2(e^{R/2}\l-1) \bigl[ 
2(\d+ih) -(\d+ih)^2(i\tau-R-S) \bigr]}\ ,\medskip \\
 \psi_2^- (\tau,\l) 
& =  {\displaystyle 
(\d-ih)(i\tau+R+S) }\medskip \\
& \qquad  
{\displaystyle 
+2 \bigl( \l e^{[i\tau+S]/2}-1\bigr)
\bigl[2(\d-ih)+(\d-ih)^2(i\tau+R+S)\bigr] } 
\ .  
\end{array} \smallskip
\end{equation}
Finally, the estimate for $\psi_1^-$ and $\psi_2^-$   
 follows at once from the explicit expressions above.
\endpf
\medskip \\

\proof[Proof of Lemma \ref{NEW-LEMMA-II*}]
Recalling formula (\ref{calculus}),
for the term $I\!I^*$ we have
\begin{align*}
I\!I^* 
& = e^{ih(S-R)} 
e^{R\j2} \biggl\{ \biggl( \frac{2}{(i\tau - R+S)^3}
        + \frac{2(\d-ih) + \j2}{(i\tau - R+S)^2} +
        \frac{(\d-ih)^2 +(\d-ih)\j2)}{(i\tau - R+S)} \biggr)
e^{-\d[i\tau-R+S]}  \\ 
& \qquad -  \biggl[ \frac{(\j2+\d)^2}{(i\tau - R+S)}
   + \biggl( \frac{2ih - \j2}{i\tau - R+S}
   - \frac{2}{(i\tau - R+S)^2} \biggr) \biggl( \j2 +\d\biggr)  \\
& \qquad +   \biggl( \frac{2}{(i\tau - R+S)^3} - 
     \frac{2ih - \j2}{(i\tau - R+S)^2} +
     \frac{ih(ih - \j2)}{i\tau - R+S} \biggr) \biggr]
     e^{(\j2+\d)[i\tau-R+S]} \biggr\}  \\ 
& =  e^{ih(S-R)-\d[i\tau-R+S]} e^{R\j2} \biggl\{ \frac{2}{(i\tau - R+S)^3} 
     +\frac{2(\d-ih)}{(i\tau-R+S)^2} +\frac{(\d-ih)^2}{i\tau-R+S} \\
 & \qquad\quad +  \j2  
\biggl( \frac{1}{(i\tau - R+S)^2} +\frac{\d-ih}{i\tau-R+S} \biggr)
\biggr\}  \\
& \quad - e^{ih(S-R)+\d[i\tau-R+S]}  e^{[i\tau +S]\j2} \biggl\{ \frac{2}{(i\tau-R+S)^3} 
-\frac{2(\d+ih)}{(i\tau-R+S)^2} +\frac{(\d+ih)^2}{i\tau-R+S}      \\
& \qquad\quad
+ \j2 \biggl( -\frac{1}{(i\tau-R+S)^2} +\frac{(\d+ih)}{i\tau-R+S} \biggr)  
\biggr\} \, .
\end{align*}

Now we must sum in $j$.  Since
$$
\sum_{j\le-1} e^{\gamma \j2} \lj = \frac{1}{\l-e^{-\gamma/2}}
\qquad\text{and}
\qquad
\sum_{j\le-1} \j2 e^{\gamma\j2}  \lj = -
   \frac{e^{-\gamma/2}}{2(\l- e^{-\gamma/2})^2} \, ,
$$
when $\big|e^{\gamma/2}\l\big|>1$, i.e. when $|\l|>e^{-R/2}$, 
we have that 
\begin{align}  
& \sum_{j\le-1} I\!I^*\lj\notag \\
& = \frac{1}{\l - e^{-R/2}}
 \biggl[ e^{-ih(R-S)-\d[i\tau-R+S] }\biggl( 
-\frac{2}{(i\tau-R+S)^3} 
-\frac{2(\d-ih)}{(i\tau-R+S)^2} -\frac{(\d-ih)^2}{i\tau-R+S} 
\biggr) \notag \\
& \qquad 
+  \frac{e^{-R/2}}{2(\l - e^{-R/2})^2}
  \biggl[ 
e^{-ih(R-S)-\d[i\tau-R+S] }\biggl( 
\frac{1}{(i\tau-R+S)^2} 
+\frac{\d-ih}{i\tau-R+S}  \biggr) \notag  \\  
& \qquad +  \frac{1}{\l - e^{-[i\tau +S]/2}}  \biggl[
e^{-ih(R-S)+\d[i\tau-R+S] }\biggl(     
\frac{2}{(i\tau-R+S)^3} -\frac{2(\d+ih)}{(i\tau-R+S)^2} 
+ \frac{(\d+ih)^2}{i\tau-R+S}  \biggr)  \notag \\
&  \qquad +  \frac{e^{-[i\tau +S]/2}}{2(\l - e^{-[i\tau +S]/2})^2}
    \biggl[ e^{-ih(R-S)+\d[i\tau-R+S] }\biggl( \frac{1}{(i\tau-R+S)^2}
-\frac{\d+ih}{i\tau-R+S} \biggr)    
\biggr] \biggr\}
\, . \label{starr}
\end{align}

Our purpose, as in the proof of Lemma
\ref{NEW-LEMMA-II}, is to examine the singularity
of the  kernel as $i\tau -R+S\ra0$.  
As in the analogous case of the sum for $j\ge0$,
some cancellation occurs.

Proceeding in a manner similar to our earlier calculation, we
temporarily introduce the notation 
$$
x = i\tau-R+S\, ,   \qquad
\text{and}\quad
\a =e^{-R/2}  \, ,
$$
so that
$$
i\tau+S
 = x+R \qquad
\text{and}\quad
e^{-[i\tau+S]/2} = \a e^{-x/2}  \, .
$$

As before we let
$f(x)  =  (\a\ex2-\l)^{-1}$ but now
$\Psi_-(x) = e^{2\d x}f(x)$.  Then we see that the sum of the terms
containing negative powers of $x$ equals $e^{-ih(R-S)+\d x} Q_{-\d}(x)$, 
where $Q_{-\d}$ is defined in (\ref{----}). 
Arguing as in the case
$j\ge0$, 
by Lemma \ref{lemma-Q} we obtain that
the sum of the terms containing the negative powers of
$i\tau-R+S$ on the righthand side of 
(\ref{starr}) equals
\begin{equation*}
e^{-ih(R-S)+\d x} Q_{-\d} \bigl(i\tau-R+S\bigr) 
= \frac{e^{-ih(R-S)+\d[i\tau-R+S] -
    R/2}}{2(e^{-[i\tau+S]/2}-\l)^2(e^{-R/2}-\l)^2} 
\psi_3^-(\tau,\l)
\ ,
\end{equation*} 
where we have set 
\begin{equation}\label{varphi1--varphi2-}
 \varphi^-_j  (\tau) = G^-_j(i\tau-R+S) \qquad\text{for}\quad
j=1,\dots,6
\  ,
\end{equation}
and 
\begin{align} 
\psi_3^-(\tau,\l)  & = 
(e^{-[i\tau+S]/2}-\l)(e^{-R/2}-\l)\varphi^-_1 (\tau)
\notag \\
& \quad
+(e^{-R/2}-\l)^2
\varphi^-_2(\tau)
+(1-e^{R/2}\l)(e^{[i\tau+S]/2}-\l)^2 \varphi^-_3 (\tau)\notag \\
& \quad 
+e^{-R/2}\l
\varphi^-_4(\tau)
+\l \varphi^-_5(\tau)+e^{-R} \varphi^-_6(\tau)
\ .\label{new-psi3-}
\end{align}

Therefore, using the argument above, 
we may simplify the expression (\ref{starr}) and 
obtain that
$$
\sum_{j\le-1} I\!I^* 
=
\frac{e^{-ih(R-S)+\d[i\tau-R+S]-R/2}}{2(e^{-[i\tau+S]/2}-\l)^2(e^{-R/2}-\l)^2} 
\psi_3^-(\tau,\l) \ .
$$

Next, recalling that $G_1^-,\dots,G_6^-$ are defined in 
(\ref{G1G6}), we have that
\begin{align*} 
|\psi_3^- (\tau,\l)| 
& \le C_M \biggl( |G^-_1 (i\tau-R+S)|
+ |G^-_2(i\tau-R+S)|
+e^{S+R/2} |G^-_3 (i\tau-R+S)| \\
& \qquad\qquad
+|G^-_4(i\tau-R+S)|
+|G^-_5(i\tau-R+S)|+e^{-R} |G^-_6(i\tau-R+S)|
\biggr) \\
& \le C_M \bigl( e^{S/2+2\d R} + e^{S+R(2\d +1/2)} +e^S + e^{2R} \bigr)\\
& \le C_M ( e^{S+R(2\d +1/2)} + e^{2R} )\ .
\end{align*}
This proves Lemma
\ref{NEW-LEMMA-II*}.
\endpf

Finally we prove Lemma \ref{lemma-Q}.

\proof[Proof of Lemma \ref{lemma-Q}]
For simplicity of notation, we provide the details for $Q_{+\d}$ only.
We write $F$ in place of $\Psi_\pm$.

We begin by noticing that
\begin{align*}
2\frac{F(x)-F(0)}{x} 
& = \frac{2}{x} \bigl( (e^{-2\d x} -1) f(x) +f(x)-f(0)\bigr) \\
& = \frac{e^{-2\d x} -1}{x/2} f(x) +2\frac{f(x)-f(0)}{x}\ ,
\end{align*} 
and that
\begin{align*}
F'(x)+F'(0) 
& = e^{-2\d x} \biggl( f'(x)-2\d f(x)\biggr) + +f'(0) -2\d f(0) \\
& = \bigl(e^{-2\d x} -1\bigr) f'(x) 
+\bigl( f'(x)+f'(0)\bigr) -2\d \bigl( e^{-2\d x} f(x)+f(0)\bigr) 
\ .
\end{align*}

Therefore
\begin{align*}
& 2\frac{F(x)-F(0)}{x} - \bigl(F'(x)+F'(0) \bigr) \\
& =  2\frac{f(x)-f(0)}{x} - \bigl( f'(x)-f'(0)\bigr) + \frac{e^{-2\d x}
  -1}{x/2} f(x) 
-\bigl( e^{-2\d x} -1\bigr) f'(x) 
\\
& \qquad\qquad 
+ 2\d \bigl( 
e^{-2\d x}f(x)+f(0)\bigr) \\
& = 2\frac{f(x)-f(0)}{x} - \bigl( f'(x)-f'(0)\bigr) 
+ \biggl( \frac{e^{-2\d x}   -1}{x/2} +2\d e^{-2\d x} +2\d \biggr)
f(x)\\
& \qquad\qquad 
-\bigl( e^{-2\d x} -1\bigr) f'(x) -2\d \bigl(
f(x)-f(0)\bigr) \\
& = 2\frac{f(x)-f(0)}{x} - \bigl( f'(x)-f'(0)\bigr) 
+ \biggl( \frac{e^{-2\d x}   -1}{x/2} +2\d e^{-2\d x} +2\d \biggr)
f(x)\\
& \qquad\qquad -\bigl( e^{-2\d x} -1+2\d x\bigr) f'(x) -2\d x\biggl(
\frac{f(x)-f(0)}{x} -f'(x) \biggr) \\
& \equiv A_1+A_2+A_3+A_4 \ .
\end{align*}

While $A_2$ and $A_3$ have simply expressions, we need to manipulate
$A_1$ and $A_4$ to obtain some simplifications. 
Setting $\Psi_1(x) =2\frac{1-e^{-x/2}}{x/2} - e^{-x/2} -1$,
we have that
\begin{align*}
A_1 
& = \frac{\a}{(\a e^{-x/2}-\l)(\a-\l)} \biggl(
\frac{1-e^{-x/2}}{x/2}\biggr)
- \biggl( \frac{\a e^{-x/2}}{2(\a e^{-x/2}-\l)^2} 
+\frac{\a}{2(\a-\l)^2} \biggr) \\
& = \frac{\a}{2(\a e^{-x/2}-\l)(\a-\l)} \biggl[
2\frac{1-e^{-x/2}}{x/2} -
\frac{e^{-x/2}(\a-\l)}{\a e^{-x/2}-\l} 
+\frac{\a e^{-x/2}-\l} {\a-\l} \biggr] \\
& = \frac{\a}{2(\a e^{-x/2}-\l)(\a-\l)} \biggl[
\Psi_1(x) -e^{-x/2} \biggl(  
\frac{\a-\l}{\a e^{-x/2}-\l} -1\biggr)
-\biggl( \frac{\a e^{-x/2}-\l} {\a-\l} -1\biggr)
\biggr] \\
& = \frac{\a}{2(\a e^{-x/2}-\l)^2(\a-\l)^2} \biggl[
\Psi_1(x)(\a e^{-x/2}-\l)(\a-\l) 
-\a\l (1-e^{-x/2})^2 \biggr] \ .
\end{align*}
Next,
\begin{align*}
A_4 
& = -2\d x\biggl[
 \frac{\a}{2(\a e^{-x/2}-\l)(\a-\l)} \biggl(\frac{1-e^{-x/2}}{x/2}
 \biggr) -\frac{\a e^{-x/2}}{2(\a e^{-x/2} -\l)^2} \biggr] \\
&  = -\frac{\a\d x}{(\a e^{-x/2}-\l)} \biggl[
\biggl(\frac{1-e^{-x/2}}{x/2} -1\biggr) \frac{1}{\a-\l} 
+\frac{\a(e^{-x/2} -1)}{(\a-\l)(\a e^{-x/2}-\l)} -\frac{e^{-x/2}-1}{\a
  e^{-x/2}-\l} \biggr] \\
& = -\frac{\a\d x}{(\a e^{-x/2}-\l)^2(\a-\l)} \biggl[
\biggl(\frac{1-e^{-x/2}}{x/2} -1\biggr) \bigl(\a e^{-x/2}-\l\bigr)  
+\bigl(e^{-x/2} -1\bigr)\l \biggr] \ .
\end{align*}

Now we turn to the second summand in (\ref{----}).  It holds that
\begin{align*}
F(x)-F(0)
& = 
\bigl( e^{-2\d x} -1\bigr) \frac{1}{\a e^{-x/2} -\l} 
+ \frac{\a}{(\a e^{-x/2} -\l)(\a-\l)}\bigl(1-e^{-x/2}\bigr) \\
& = \frac{1}{(\a e^{-x/2} -\l)(\a-\l)}\biggl[
\bigl(e^{-2\d x} -1\bigr)\bigl(\a  -\l\bigr) +\a \bigl(
1-e^{-x/2} \bigr) \biggr] 
\end{align*}
and
\begin{align*}
e^{-2\d x} f'(x) -f'(0) 
& = \bigl( e^{-2\d x} -1\bigr) 
\frac{\a e^{-x/2}}{2(\a e^{-x/2}  -\l)^2} +
\frac{\a(e^{-x/2}-1)}{2(\a e^{-x/2}-\l)^2} \\
& \qquad\qquad \qquad\qquad +
\frac{\a^2(e^{-x/2}-1)\bigl( \a(e^{-x/2}+1)+2\l\bigr)}{
2(\a e^{-x/2}-\l)^2 (\a-\l)^2} \\
& = \frac{\a}{2(\a e^{-x/2}-\l)^2 (\a-\l)^2} 
\biggl[ e^{-x/2} \bigl(e^{-2\d x}-1\bigr)(\a-\l)^2  \\
& \qquad
+ \bigl(e^{-x/2}-1\bigr)(\a-\l)^2  
+\a\bigl(1-e^{x/2}\bigr)\bigl(
\a(e^{-x/2}+1)+2\l\bigr) \biggr] \ .
\end{align*}

Therefore
\begin{align*}
& Q_\d (x) \\
& = 
\frac{\a}{2(\a e^{-x/2}-\l)^2 (\a-\l)^2} \biggl\{
\frac{ihx-1}{x^2} \biggl[ 
\Psi_1(x)(\a e^{-x/2}-\l)(\a-\l) 
-\a\l (1-e^{-x/2})^2 \\
& \qquad
+ 2\biggl( \frac{e^{-2\d x} -1}{x/2} +2\d e^{-2\d x}+2\d \biggr)
\frac{(\a-\l)^2(\a e^{-x/2} -\l)}{\a} - \bigl( e^{-2\d x} -1+2\d x\bigr) 
e^{-x/2} (\a-\l)^2 \\
&\qquad
-2\d x\biggl(\frac{1-e^{-x/2}}{x/2} -1\biggr) 
\bigl(\a e^{-x/2}-\l\bigr)  
(\a-\l)
+\d x\bigl(e^{-x/2} -1\bigr)\l(\a-\l)\biggr] 
 \\
& \quad
+ \frac1x \biggl[
2(h^2-\d^2)\biggl( \bigl( e^{-2\d x} -1\bigr) 
\frac{(\a -\l)^2 (\a e^{-x/2}-\l)}{\a} + \bigl(1-e^{-x/2}\bigr)
(\a-\l)(\a e^{-x/2} -\l) \biggr) \\
& \qquad
+ \d\biggl( e^{-x/2} \bigl(e^{-2\d x}-1\bigr)(\a-\l)^2  
+ \bigl(e^{-x/2}-1\bigr)(\a-\l)^2  
+\a\bigl(1-e^{-x/2}\bigr)\bigl(
\a(e^{-x/2}+1)+2\l\bigr) \biggr) \biggr] \biggr\} \\
& = \frac{\a}{2(\a e^{-x/2}-\l)^2 (\a-\l)^2} 
\biggl[ 
(\a e^{-x/2}-\l)(\a-\l) G_1(x) + (\a-\l)^2 G_2(x) \\
& \qquad\qquad
+\frac{(\a-\l)^2(\a e^{-x/2}-\l)
}{\a} G_3(x)
+\a\l G_4(x)+\l(\a-\l) G_5(x) +\a^2 G_6(x) \biggr] \, ,
\end{align*}
where $G_1, \dots, G_6$ are entire functions of $x$, independent of $\a$
and $\l$. 
Explicitly, recalling the dependence on the choice of
the sign in front of $\d$, they are
 given by
\begin{equation}\label{G1G6}
\begin{array}{l}
G^\pm_1(x) ={\displaystyle 
\frac{ihx-1}{x^2} \biggl( 
2\frac{1-e^{-x/2}}{x/2} - e^{-x/2} -1\biggr)
-2\d x\biggl(\frac{1-e^{-x/2}}{x/2} -1\biggr) 
 +(h^2 -\d^2)
  \frac{1-e^{-x/2}}{x/2} }\medskip\\
G^\pm_2(x) = {\displaystyle 
\frac{ihx-1}{x^2} 
\biggl[ -e^{-x/2}\bigl( e^{\mp2\d x} -1\pm 2\d x\bigr) \biggr]
\pm\d e^{-x/2} \biggl( \frac{e^{\mp2\d x} -1}{x/2} + \frac{e^{-x/2} -1}{x}
\biggr) }\medskip\\
G^\pm_3(x) = {\displaystyle 
2\frac{ihx-1}{x^2} \biggl( 
\frac{e^{\mp2\d x} -1}{x/2} \pm 2\d e^{-2\d x} \pm2\d\biggr) +(h^2-\d^2)
\frac{e^{\mp2\d x} -1}{x/2}  }\medskip\\
G^\pm_4(x) = {\displaystyle 
-\frac{ihx-1}{x^2} \bigl(1-e^{-x/2}\bigr)^2 \pm 2\d e^{-x/2}
\frac{1-e^{-x/2}}{x} }\medskip\\
G^\pm_5(x) = {\displaystyle 
\mp2\d \frac{ihx-1}{x^2} \bigl(e^{-x/2}-1\bigr) }\medskip\\
G^\pm_6(x) = {\displaystyle 
 \pm\d \frac{e^{-x}-1}{x} }
\, .
\end{array}
\end{equation}
That completes the proof of Lemma \ref{lemma-Q}.
\bigskip
\endpf

\section{Concluding Remarks}\label{concluding-remarks}
\vspace*{.12in}

We have provided a rather complete analysis of the Bergman kernel for
the non-smooth worms $D_\b$ and $D'_\b$.  Our work has been
facilitated by the fact that the boundaries of these domains are Levi
flat, hence they are qualitatively like product domains (and our
resulting formulas reflect this structure).  Another way to look at
this is that the domains $D_\b$ and $D'_\b$ possess many more
symmetries with respect to the smooth worm $\mathcal W_\b$, since
the translations in the real part of $z_1$, $(z_1,z_2)\ra (z_1+a,z_2)$
with $a\in\RR$, are automorphisms of both 
$D_\b$ and $D'_\b$.  

It is a matter of considerable interest to perform the analogous
analysis on the smooth worm $\mathcal W_\b$.  That work
will be of a different nature, for the boundary of $\mathcal W_\b$
is strongly pseudoconvex at all points except those on a singular
annulus in the boundary.  There certainly is no ``product structure''.
We hope to complete such an analysis in a future paper.

The results that we present provide a new way to view the
negative results of Kiselman, Barrett, and Christ about mapping
properties of the Bergman projection on worm domains.
We hope that they will serve as a stepping stone to future studies.


\bigskip


\begin{thebibliography}{CheW}


\bibitem[A]{APF}  L. Apfel, Localization properties and boundary
behavior of the Bergman kernel, thesis, Washington University, 2003.

\bibitem[Bar1]{BAR1}  D. Barrett, Regularity of the Bergman projection
on domains with transverse symmetries, {\it Math.\ Ann.} 258(1982),
441--446.

\bibitem[Ba2]{BAR2}  D. Barrett, Behavior of the Bergman projection
on the Diederich-Forn\ae ss worm, {\it Acta Math.} 168(1992), 1--10.

\bibitem[Be1]{BEL1}  S. R. Bell, Biholomorphic mappings and the
$\overline{\partial}$-problem, {\it Ann.\ Math.} 114(1981), 103--113.

\bibitem[Be2]{BEL2}  S. R. Bell, Local boundary behavior of
proper holomorphic mappings, {\it Proceedings of Symposia in Pure Math.},
American Math.\ Society, Providence, RI, 41, 1984, 1--7.

\bibitem[BoS1]{BOS1}  H. Boas and E. Straube, Complete Hartogs domains
in $\CC^2$ have regular Bergman and Szeg\"{o} projections,
{\it Math.\ Zeit.} 201(1989), 441--454.

\bibitem[BoS2]{BOS2}  H. Boas and E. Straube, Equivalence of regularity
for the Bergman projection and the $\overline{\partial}$-Neumann
operator, {\it Manuscripta Math.} 67(1990), 25--33.

\bibitem[Che1]{CHE1}  S.-C. Chen, Characterization of the automorphisms
on the Barrett and the Diederich-Fornaess worm domains, {\it Trans.\
Amer.\ Math.\ Soc.} 338(1993), 431--440.

\bibitem[Che2]{CHE2}  S.-C. Chen, A counterexample to the differentiability
of the Bergman kernel, {\it Proc.\ AMS} 124(1996), 1807--1810.

\bibitem[CheS]{CHS}  S.-C. Chen and M.-C. Shaw, {\it Partial Differential
Equations in Several Complex Variables}, American Mathematical Society,
Providence, RI, 2001.

\bibitem[Chr]{CHR}  M. Christ, Global $C^\infty$ irregularity of the
$\overline{\partial}$-Neumann problem for worm domains, {\it J. of the
Amer.\ Math. Soc.} 9(1996), 1171--1185.

\bibitem[DF]{DF}  K. Diederich and J. E. Forn\ae ss, Pseudoconvex
domains:  An example with nontrivial Nebenh\"{u}lle, {\it Math. Ann.}
225(1977), 275--292.

\bibitem[FR]{FORR} F. Forelli and W. Rudin, Projections on
spaces of holomorphic functions in balls, {\it Indiana Univ.\ Math.\
J.} 24(1974/5), 593--602.

\bibitem[Ke]{KER}  The Bergman kernel function.  Differentiability
at the boundary, {\em Math. Ann.} 195(1972), 149-158.

\bibitem[Ki]{KIS}  C. Kiselman, A study of the Bergman projection
in certain Hartogs domains, {\it Proc.\ Symposia in Pure Math.},
American Math.\ Society, Providence, RI, 52, III(1991), 219--231.

\bibitem[Kr1]{KRA1}  S. G. Krantz, {\it Function Theory of Several
Complex Variables}, $2^{\rm nd}$ ed., American Mathematical Society,
Providence, RI, 2000.

\bibitem[Kr2]{KRA2}  S. G. Krantz, {\it Partial Differential Equations
and Complex Analysis}, CRC Press, Boca Raton, FL, 1992.

\bibitem[KrP]{KP2} S. G. Krantz, M. M. Peloso,
Analysis on the worm domains, preprint 2007.

\bibitem[L]{LIG}  E. Ligocka, Remarks on the Bergman kernel function
of a worm domain, {\it Studia Mathematica} 130(1998), 109--113.

\bibitem[Sa]{SAD}  C. Sadosky, {\it Interpolation of Operators
and Singular Integrals}, Marcel Dekker, New York, 1979.

\bibitem[Si]{SIU}  Y.-T. Siu, Non-H\"{o}lder property of Bergman
projection of smooth worm domain, {\it Aspects of Mathematics---Algebra,
Geometry, and Several Complex Variables}, N. Mok (ed.), University 
of Hong Kong, 1996, 264--304.

\bibitem[Str]{Str} E. M. Straube, MR\,1149863 and personal
  communication. 

\bibitem[Z]{Zhu} K. Zhu, {\em Spaces of holomorphic functions in the
    unit ball}, Graduate Texts in Mathematics 226, Springer-Verlag,
    New York, 2005. 
\end{thebibliography}
\end{document}